\documentclass{elsarticle}

\usepackage[margin=1in]{geometry}

\usepackage{bbm}
\usepackage{graphicx}
\usepackage{pstool}
\usepackage{wrapfig}
\usepackage{caption}
\usepackage{subcaption}
\usepackage[section]{placeins} 

\usepackage{fancyhdr} 
\usepackage{lastpage} 
\usepackage{extramarks} 

\usepackage{xcolor} 
\usepackage{enumerate} 
\usepackage{paralist} 
\usepackage{amsmath, amsthm, amssymb, mathtools}
\usepackage{mathabx, pifont, stmaryrd} 
\usepackage[explicit]{titlesec} 
\usepackage{etoolbox} 
\usepackage{bibentry} 
\makeatletter\let\saved@bibitem\@bibitem\makeatother 
\usepackage[colorlinks, bookmarksopen, bookmarksnumbered,
            citecolor=red,urlcolor=red]{hyperref} 
\makeatletter\let\@bibitem\saved@bibitem\makeatother 
\usepackage{rotating}

\usepackage{algorithm}
\usepackage{algorithmic}

\newtheorem{remark}{Remark}

\theoremstyle{definition}

\usepackage{bm} 
\usepackage{amsfonts} 

\DeclareMathOperator*{\argmin}{arg\,min}

\newcommand{\ds}[1]{\ensuremath{\displaystyle{#1}}}
\newcommand{\func}[3]{\ensuremath{#1 : #2 \rightarrow #3}}

\newcommand{\norm}[1]{\ensuremath{\left\| #1 \right\|}}

\newcommand{\suchthat}{\mathrel{}\middle|\mathrel{}}

\newcommand{\optconOne}[3]{
\begin{aligned}
& \underset{#1}{\text{minimize}}
& & #2 \\
& \text{subject to} & & #3
\end{aligned}}

\newcommand{\pder}[2]{\ensuremath{\frac{\partial #1}{\partial #2}}}

\newcommand{\KKTMat}{\ensuremath{\Abm}}
\newcommand{\KKTb}{\ensuremath{\bbm}}

\newcommand{\KKTC}{\ensuremath{\tilde{\Abm}_{C}}}
\newcommand{\Ju}{\ensuremath{\Jbm_\ubm}}
\newcommand{\Jx}{\ensuremath{\Jbm_\xbm}}

\newcommand{\Prec}{\ensuremath{\tilde{\Abm}}}
\newcommand{\PrecInv}{\ensuremath{\Prec^{-1}}}
\newcommand{\PrecV}[1]{\ensuremath{\Prec_{\text{#1}}}}
\newcommand{\JBJ}{\ensuremath{\tilde{\Jbm}_{\text{BJ}}}}

\newcommand{\JBILU}{\ensuremath{\tilde{\Jbm}_{\text{BILU}}}}

\newcommand{\ABJ}{\PrecV{BJ}}
\newcommand{\ABJp}{\PrecV{BJp0}}
\newcommand{\ABJilu}{\PrecV{BJ/ilu}}

\newcommand{\ABILU}{\PrecV{BILU}}
\newcommand{\ABILUp}{\PrecV{BILUp0}}
\newcommand{\ABILUilu}{\PrecV{BILU/ilu}}

\newcommand{\Ao}{\PrecV{0}}
\newcommand{\Aop}{\PrecV{0p0}}

\newcommand{\Prol}{\ensuremath{\Pbm}}
\newcommand{\Rest}{\ensuremath{\Qbm}}

\newcommand{\sumEl}{\ensuremath{\sum^{\vert \Ecal_h \vert }_{e=1}}}

\newcommand{\ecyl}{\texttt{cylinder}}
\newcommand{\edmnd}{\texttt{diamond}}

\newcommand{\Ecal}{\ensuremath{\mathcal{E}}}

\newcommand{\Gcal}{\ensuremath{\mathcal{G}}}
\newcommand{\Hcal}{\ensuremath{\mathcal{H}}}

\newcommand{\Lcal}{\ensuremath{\mathcal{L}}}

\newcommand{\Ncal}{\ensuremath{\mathcal{N}}}

\newcommand{\Pcal}{\ensuremath{\mathcal{P}}}

\newcommand{\Vcal}{\ensuremath{\mathcal{V}}}
\newcommand{\Wcal}{\ensuremath{\mathcal{W}}}

\newcommand{\Gbb}{\ensuremath{\mathbb{G}}}

\newcommand{\Nbb}{\ensuremath{\mathbb{N} }}

\newcommand{\Rbb}{\ensuremath{\mathbb{R} }}

\newcommand\Abm{{\ensuremath{\bm{A}}}}
\newcommand\Bbm{{\ensuremath{\bm{B}}}}
\newcommand\Cbm{{\ensuremath{\bm{C}}}}
\newcommand\Dbm{{\ensuremath{\bm{D}}}}

\newcommand\Fbm{{\ensuremath{\bm{F}}}}
\newcommand\Gbm{{\ensuremath{\bm{G}}}}

\newcommand\Ibm{{\ensuremath{\bm{I}}}}
\newcommand\Jbm{{\ensuremath{\bm{J}}}}

\newcommand\Lbm{{\ensuremath{\bm{L}}}}
\newcommand\Mbm{{\ensuremath{\bm{M}}}}

\newcommand\Pbm{{\ensuremath{\bm{P}}}}
\newcommand\Qbm{{\ensuremath{\bm{Q}}}}
\newcommand\Rbm{{\ensuremath{\bm{R}}}}

\newcommand\Ubm{{\ensuremath{\bm{U}}}}

\newcommand\bbm{{\ensuremath{\bm{b}}}}

\newcommand\gbm{{\ensuremath{\bm{g}}}}

\newcommand\rbm{{\ensuremath{\bm{r}}}}
\newcommand\sbm{{\ensuremath{\bm{s}}}}

\newcommand\ubm{{\ensuremath{\bm{u}}}}
\newcommand\vbm{{\ensuremath{\bm{v}}}}
\newcommand\wbm{{\ensuremath{\bm{w}}}}
\newcommand\xbm{{\ensuremath{\bm{x}}}}
\newcommand\ybm{{\ensuremath{\bm{y}}}}
\newcommand\zbm{{\ensuremath{\bm{z}}}}

\newcommand\lambdabold{{\ensuremath{\boldsymbol{\lambda}}}}

\newcommand\etabold{{\ensuremath{\boldsymbol{\eta}}}}

\newcommand\phibold{{\ensuremath{\boldsymbol{\phi}}}}

\newcommand\zerobold{\ensuremath{\mathbf{0}}}

\usepackage{tikz}
\usepackage{pgfplots}
\usepackage{pgfplotstable, filecontents, booktabs}
\pgfplotsset{compat=1.9}

\usetikzlibrary{pgfplots.groupplots}
\usepgfplotslibrary{fillbetween}
\usetikzlibrary{calc,fit,matrix,arrows,automata,positioning,shapes}
\usetikzlibrary{arrows.meta}

\pgfplotsset{select coords between index/.style 2 args={
    x filter/.code={
        \ifnum\coordindex<#1\fi
        \ifnum\coordindex>#2\fi
    }
}}

\tikzset{
 invisible/.style={opacity=0},
 visible on/.style={alt={#1{}{invisible}}},
 alt/.code args={<#1>#2#3}{%
   \alt<#1>{\pgfkeysalso{#2}}{\pgfkeysalso{#3}}
 },
}

\newcommand{\colorbarMatlabParula}[5]{
\begin{tikzpicture}
\begin{axis}[
   hide axis, scale only axis,
   height=0pt, width=0pt,
   colormap={parula}{rgb255=(62,38,168) rgb255=(62,39,172) rgb255=(63,40,175) rgb255=(63,41,178) rgb255=(64,42,180) rgb255=(64,43,183) rgb255=(65,44,186) rgb255=(65,45,189) rgb255=(66,46,191) rgb255=(66,47,194) rgb255=(67,48,197) rgb255=(67,49,200) rgb255=(67,50,202) rgb255=(68,51,205) rgb255=(68,52,208) rgb255=(69,53,210) rgb255=(69,55,213) rgb255=(69,56,215) rgb255=(70,57,217) rgb255=(70,58,220) rgb255=(70,59,222) rgb255=(70,61,224) rgb255=(71,62,225) rgb255=(71,63,227) rgb255=(71,65,229) rgb255=(71,66,230) rgb255=(71,68,232) rgb255=(71,69,233) rgb255=(71,70,235) rgb255=(72,72,236) rgb255=(72,73,237) rgb255=(72,75,238) rgb255=(72,76,240) rgb255=(72,78,241) rgb255=(72,79,242) rgb255=(72,80,243) rgb255=(72,82,244) rgb255=(72,83,245) rgb255=(72,84,246) rgb255=(71,86,247) rgb255=(71,87,247) rgb255=(71,89,248) rgb255=(71,90,249) rgb255=(71,91,250) rgb255=(71,93,250) rgb255=(70,94,251) rgb255=(70,96,251) rgb255=(70,97,252) rgb255=(69,98,252) rgb255=(69,100,253) rgb255=(68,101,253) rgb255=(67,103,253) rgb255=(67,104,254) rgb255=(66,106,254) rgb255=(65,107,254) rgb255=(64,109,254) rgb255=(63,110,255) rgb255=(62,112,255) rgb255=(60,113,255) rgb255=(59,115,255) rgb255=(57,116,255) rgb255=(56,118,254) rgb255=(54,119,254) rgb255=(53,121,253) rgb255=(51,122,253) rgb255=(50,124,252) rgb255=(49,125,252) rgb255=(48,127,251) rgb255=(47,128,250) rgb255=(47,130,250) rgb255=(46,131,249) rgb255=(46,132,248) rgb255=(46,134,248) rgb255=(46,135,247) rgb255=(45,136,246) rgb255=(45,138,245) rgb255=(45,139,244) rgb255=(45,140,243) rgb255=(45,142,242) rgb255=(44,143,241) rgb255=(44,144,240) rgb255=(43,145,239) rgb255=(42,147,238) rgb255=(41,148,237) rgb255=(40,149,236) rgb255=(39,151,235) rgb255=(39,152,234) rgb255=(38,153,233) rgb255=(38,154,232) rgb255=(37,155,232) rgb255=(37,156,231) rgb255=(36,158,230) rgb255=(36,159,229) rgb255=(35,160,229) rgb255=(35,161,228) rgb255=(34,162,228) rgb255=(33,163,227) rgb255=(32,165,227) rgb255=(31,166,226) rgb255=(30,167,225) rgb255=(29,168,225) rgb255=(29,169,224) rgb255=(28,170,223) rgb255=(27,171,222) rgb255=(26,172,221) rgb255=(25,173,220) rgb255=(23,174,218) rgb255=(22,175,217) rgb255=(20,176,216) rgb255=(18,177,214) rgb255=(16,178,213) rgb255=(14,179,212) rgb255=(11,179,210) rgb255=(8,180,209) rgb255=(6,181,207) rgb255=(4,182,206) rgb255=(2,183,204) rgb255=(1,183,202) rgb255=(0,184,201) rgb255=(0,185,199) rgb255=(0,186,198) rgb255=(1,186,196) rgb255=(2,187,194) rgb255=(4,187,193) rgb255=(6,188,191) rgb255=(9,189,189) rgb255=(13,189,188) rgb255=(16,190,186) rgb255=(20,190,184) rgb255=(23,191,182) rgb255=(26,192,181) rgb255=(29,192,179) rgb255=(32,193,177) rgb255=(35,193,175) rgb255=(37,194,174) rgb255=(39,194,172) rgb255=(41,195,170) rgb255=(43,195,168) rgb255=(44,196,166) rgb255=(46,196,165) rgb255=(47,197,163) rgb255=(49,197,161) rgb255=(50,198,159) rgb255=(51,199,157) rgb255=(53,199,155) rgb255=(54,200,153) rgb255=(56,200,150) rgb255=(57,201,148) rgb255=(59,201,146) rgb255=(61,202,144) rgb255=(64,202,141) rgb255=(66,202,139) rgb255=(69,203,137) rgb255=(72,203,134) rgb255=(75,203,132) rgb255=(78,204,129) rgb255=(81,204,127) rgb255=(84,204,124) rgb255=(87,204,122) rgb255=(90,204,119) rgb255=(94,205,116) rgb255=(97,205,114) rgb255=(100,205,111) rgb255=(103,205,108) rgb255=(107,205,105) rgb255=(110,205,102) rgb255=(114,205,100) rgb255=(118,204,97) rgb255=(121,204,94) rgb255=(125,204,91) rgb255=(129,204,89) rgb255=(132,204,86) rgb255=(136,203,83) rgb255=(139,203,81) rgb255=(143,203,78) rgb255=(147,202,75) rgb255=(150,202,72) rgb255=(154,201,70) rgb255=(157,201,67) rgb255=(161,200,64) rgb255=(164,200,62) rgb255=(167,199,59) rgb255=(171,199,57) rgb255=(174,198,55) rgb255=(178,198,53) rgb255=(181,197,51) rgb255=(184,196,49) rgb255=(187,196,47) rgb255=(190,195,45) rgb255=(194,195,44) rgb255=(197,194,42) rgb255=(200,193,41) rgb255=(203,193,40) rgb255=(206,192,39) rgb255=(208,191,39) rgb255=(211,191,39) rgb255=(214,190,39) rgb255=(217,190,40) rgb255=(219,189,40) rgb255=(222,188,41) rgb255=(225,188,42) rgb255=(227,188,43) rgb255=(230,187,45) rgb255=(232,187,46) rgb255=(234,186,48) rgb255=(236,186,50) rgb255=(239,186,53) rgb255=(241,186,55) rgb255=(243,186,57) rgb255=(245,186,59) rgb255=(247,186,61) rgb255=(249,186,62) rgb255=(251,187,62) rgb255=(252,188,62) rgb255=(254,189,61) rgb255=(254,190,60) rgb255=(254,192,59) rgb255=(254,193,58) rgb255=(254,194,57) rgb255=(254,196,56) rgb255=(254,197,55) rgb255=(254,199,53) rgb255=(254,200,52) rgb255=(254,202,51) rgb255=(253,203,50) rgb255=(253,205,49) rgb255=(253,206,49) rgb255=(252,208,48) rgb255=(251,210,47) rgb255=(251,211,46) rgb255=(250,213,46) rgb255=(249,214,45) rgb255=(249,216,44) rgb255=(248,217,43) rgb255=(247,219,42) rgb255=(247,221,42) rgb255=(246,222,41) rgb255=(246,224,40) rgb255=(245,225,40) rgb255=(245,227,39) rgb255=(245,229,38) rgb255=(245,230,38) rgb255=(245,232,37) rgb255=(245,233,36) rgb255=(245,235,35) rgb255=(245,236,34) rgb255=(245,238,33) rgb255=(246,239,32) rgb255=(246,241,31) rgb255=(246,242,30) rgb255=(247,244,28) rgb255=(247,245,27) rgb255=(248,247,26) rgb255=(248,248,24) rgb255=(249,249,22) rgb255=(249,251,21) },
   colorbar horizontal,
   point meta min=#1, point meta max=#5,
   colorbar style={width=10cm, xtick={#1,#2,#3,#4,#5}}
]
\addplot [draw=none] coordinates {(0,0)};
\end{axis}
\end{tikzpicture}
}

\usepackage{setspace}

\newbool{fastcompile}
\setbool{fastcompile}{false}

\begin{document}
	\title{Preconditioned iterative solvers for constrained \\ high-order implicit shock tracking methods}
	
	\author[1]{Jakob Vandergrift\corref{cor}%
		\fnref{fn1}}
	\ead{vandergrift@fdy.tu-darmstadt.de}
	
	\author[rvt1]{Matthew J. Zahr\fnref{fn2}\corref{cor1}}
	\ead{mzahr@nd.edu}
	
	\affiliation[1]{organization={TU Darmstadt, Chair of Fluid Dynamics},\\ 
		addressline={Otto-Berndt-Str. 2 },
		postcode={64287}, 
		city={Darmstadt}, 
		country={Germany}}
	
	\address[rvt1]{Department of Aerospace and Mechanical Engineering, University
		of Notre Dame, Notre Dame, IN 46556, United States}
	\cortext[cor1]{Corresponding author}
	
	\fntext[fn1]{PhD Student, Department of Mechanical Engineering, Technical University of Darmstadt}
	\fntext[fn2]{Assistant Professor, Department of Aerospace and Mechanical
		Engineering, University of Notre Dame}
	
	\begin{keyword} 
		Shock fitting, high-order methods, discontinuous Galerkin, constrained optimization,
		preconditioners, iterative solvers
	\end{keyword}
	
	\begin{abstract}
		High-order implicit shock tracking (fitting) is a class of high-order numerical methods that use numerical optimization
		to simultaneously compute a high-order approximation to a conservation law solution and align elements
		of the computational mesh with non-smooth features. This alignment ensures that non-smooth features
		are perfectly represented by inter-element jumps and high-order basis functions approximate smooth regions of the solution
		without nonlinear stabilization, which leads to accurate approximations on traditionally coarse meshes. In this work, we
		devise a family of preconditioners for the saddle point linear system that defines the step toward optimality at each iteration
		of the optimization solver so Krylov solvers can be effectively used. Our preconditioners integrate standard preconditioners
		from constrained optimization with popular preconditioners for discontinuous Galerkin discretizations such as block Jacobi,
		block incomplete LU factorizations with minimum discarded fill reordering, and $p$-multigrid. Thorough studies are performed
		using two inviscid compressible flow problems to evaluate the effectivity of each preconditioner in this family and their sensitivity
		to critical shock tracking parameters such as the mesh and Hessian regularization, linearization state, and resolution of the
		solution space.
	\end{abstract}
	
	\maketitle
	
	
	\section{Introduction}
	\label{sec:intro}
	
	Accurate and robust simulation of shock-dominated flows remains a significant challenge
for modern computational fluid dynamics methods. High-order methods, such as discontinuous
Galerkin (DG) methods \cite{cockburn_rungekutta_2001, hesthaven_nodal_2008}, have
received considerable attention because they are highly accurate per degree of
freedom, introduce minimal dissipation, provide geometric flexibility, and exhibit
excellent parallel scalability \cite{2013_wang_highordercfd}. Despite these advantages,
high-order methods are known to lack robustness for shock-dominated flows because
high-order approximation of shocks and contact discontinuities leads to spurious oscillations
that cause a breakdown of numerical solvers.

A new class of numerical methods, known as high-order implicit shock tracking (fitting) \cite{zahr2018shktrk,corrigan2019moving},
has emerged that uses numerical optimization to simultaneously compute a high-order approximation
to a conservation law solution and align elements of the computational mesh with the non-smooth
features. This ensures non-smooth features are perfectly represented by inter-element jumps and
high-order basis functions approximate smooth regions of the solution without nonlinear stabilization,
which leads to accurate approximations on traditionally coarse meshes. These techniques have been used
to resolve steady and unsteady, inert and reacting shock-dominated flows in the transonic, supersonic, and
hypersonic regimes.
To this point, implicit shock tracking research has focused on the variational formulation
\cite{zahr_implicit_2020,corrigan2019moving,kercher2020moving,kercher2020least},
proper choice of objective and constraint functions \cite{zahr2018shktrk,zahr_implicit_2020},
robust solvers for the optimization problems \cite{huangRobustHighorderImplicit2022a}, and various
applications \cite{huang_high-order_2023}. To this point, little-to-no attention has been given to solvers
for the linearized optimality system that defines the search direction at each optimization iteration.

In this work, we propose a family of preconditioners for the linearized optimality systems that arise
in sequential quadratic programming (SQP) solvers for constrained implicit shock tracking methods. 
For concreteness, we focus on the High-Order Implicit Shock Tracking (HOIST) method 
\cite{zahr_implicit_2020,huangRobustHighorderImplicit2022a} that
uses an enriched residual as the objective function, although our preconditioners generalize to other
objective functions such as the Rankine-Hugoniot conditions \cite{corrigan2019moving}. In SQP methods, the step
toward optimality is the solution of the linearized Karush-Kuhn-Tucker (KKT) conditions of the original
constrained optimization problem. The proposed preconditioners are built on a class of \textit{constrained
preconditioners} \cite{kellerConstraintPreconditioningIndefinite2000} that mimic the structure of the original
saddle point problems, which have been successfully combined with conjugate gradient methods
\cite{colemanPreconditionedConjugateGradient,gouldSolutionEqualityConstrained2001} and other Krylov subspace methods
\cite{luksanNumericalExperienceIterative2001,dollarConstraintStylePreconditionersRegularized2007} to solve nonlinear
programming problems. We build a family of cost-effective constrained preconditioners
by approximating the constraint Jacobian with standard preconditioners from the DG community (e.g., block Jacobi
and block incomplete LU factorization with minimum discarded fill reordering) \cite{persson_newton-gmres_2008},
dropping some blocks of the Lagrangian
Hessian, and using standard preconditioners to approximate other blocks. Similar approaches that approximate the constraint
Jacobian and neglect blocks of the Lagrangian Hessian have been used to develop matrix-based preconditioners for optimal control
problems \cite{biros_parallel_2000}. A two-level $p$-multigrid acceleration strategy is defined that can be used in combination with
any preconditioner in the proposed family. Thorough studies are performed using two inviscid compressible flow problems to evaluate the 
effectivity of each preconditioner in this family and their sensitivity to critical shock tracking parameters such as the mesh and Hessian 
regularization, linearization state, and resolution of the solution space.

	The remainder of the  paper is organized as follows: Section \ref{sec:govern} introduces the transformed system of conservation laws and its high-order DG discretization. Section \ref{sec:hoist} presents the HOIST formulation and details the sparsity structure of the linearized optimality system. Section \ref{sec:linsolv} discusses popular preconditioners for DG discretizations and uses these to derive specialized matrix-based preconditioners for the implicit shock tracking linearized optimality system. Section \ref{sec:mulgrd:hoist} presents a two-level $p$-multigrid method, which is developed and integrated with each preconditioner proposed. Extensive experimentations with all preconditioners proposed, highlighting their dependence on several crucial optimization solver parameters, are presented and analyzed in Section~\ref{sec:numexp}. Finally, Section~\ref{sec:conclude} offers conclusions and identifies relevant avenues for future research.

	\section{Governing equations and high-order discretization}
	\label{sec:govern}
	In this section, we introduce the governing partial differential
	equations, specifically a system of steady inviscid conservation laws
	and its transformation to a reference
	domain so that domain deformations appear explicitly in the governing
	equations (Section~\ref{sec:govern:claw}). Secondly, we present its discretization
	using a high-order DG method (high-order with respect to both the
	solution and geometry) (Section~\ref{sec:govern:dg}). Lastly, we examine the sparsity structure of discrete operators for DG methods (Section~\ref{sec:govern:sparse}) which are needed in applications employing nonlinear solvers (Section~\ref{sec:govern:nlsolve}).
	
	\subsection{Transformed system of conservation laws}
	\label{sec:govern:claw}
	Consider a general system of $m$ inviscid conservation laws, defined on the
	fixed domain $\Omega \subset \Rbb^d$ and subject to appropriate boundary
	conditions,
	\begin{equation} \label{eqn:claw-phys}
		\nabla\cdot F(U) = S(U) \quad \text{in}~~\Omega,
	\end{equation}
	where $\func{U}{\Omega}{\Rbb^m}$ is the solution of the system of
	conservation laws, $\func{F}{\Rbb^m}{\Rbb^{m\times d}}$ is the flux
	function, $\func{S}{\Rbb^m}{\Rbb^m}$ is the source term,
	$\ds{\nabla \coloneqq (\partial_{x_1},\dots,\partial_{x_d})}$
	is the gradient operator in the physical domain, and the boundary of
	the domain $\partial\Omega$ has outward unit normal
	$\func{n}{\partial\Omega}{\Rbb^d}$. In general, the solution $U(x)$ may contain discontinuities, in which case, the conservation laws
	(\ref{eqn:claw-phys}) hold away from these
	and the Rankine-Hugoniot conditions \cite{majda2012compressible} 
	hold at the discontinuities.
	
	Before discretizing equation (\ref{eqn:claw-phys}), it is advantageous to explicitly handle deformations to the conservation law domain $\Omega$. These deformations, which will eventually occur due to mesh adjustments as nodal coordinates move to track discontinuities, can be managed by transforming the problem to a fixed reference domain, denoted as $\Omega_0 \subset \mathbb{R}^d$. Let $\mathcal{G}$ represent the set of diffeomorphisms from the reference domain $\Omega_0$ to the physical domain $\Omega$ defined as
	\begin{equation} \label{eqn:dom-map}
		\Gbb :=  \left\{\func{\Gcal}{\Omega_0}{\Omega} \suchthat 
		\Gcal : X \mapsto \Gcal(X)\right\} .
	\end{equation}
	For any $\Gcal\in\Gbb$, the conservation law on the physical domain $\Omega$ is transformed to a conservation law on the reference domain $\Omega_0$ 
	\begin{equation} \label{eqn:claw-ref}
		\bar\nabla\cdot\bar{F}(\bar{U};G) =
		\bar{S}(\bar{U};g) \quad \text{in}~~\Omega_0.
	\end{equation}
	Here we denote by $\func{\bar{U}}{\Omega_0}{\Rbb^m}$ the solution of the transformed conservation law, by $\func{\bar{F}}{\Rbb^m\times\Rbb^{d\times d}}{\Rbb^{m\times d}}$ the
	transformed flux function, by $\bar\nabla \coloneqq (\partial_{X_1},\dots,\partial_{X_d})$ the gradient operator on the reference domain, by $\func{G}{\Omega_0}{\Rbb^{d\times d}}$ the deformation gradient and by $\func{g}{\Omega_0}{\Rbb}$ the mapping Jacobian. The latter are defined as
	\begin{equation}
		G = \bar\nabla \Gcal, \qquad g = \det G.
	\end{equation}
	The unit outward normal to the reference domain is denoted $\func{N}{\partial\Omega_0}{\Rbb^d}$ and the following relation to the unit normal in the physical domain holds 
	\begin{equation} \label{eqn:transf-normal}
		n \circ \Gcal = \frac{g G^{-T}N}{\norm{g G^{-T}N}}.
	\end{equation}
	For any $X\in\Omega_0$, the transformed and physical solution are related by
	\begin{equation}
		\bar{U}(X) = U(\Gcal(X)),
	\end{equation}
	whereas the transformed flux and source term are defined as
	\begin{equation}
		\bar{F} : (\bar{W}; \Theta) \mapsto (\det\Theta) F(\bar{W}) \Theta^{-T},
		\qquad
		\bar{S} : (\bar{W}; q) \mapsto q S(\bar{W}).
	\end{equation}
	\begin{remark} \label{rem:refphys}
		In general, the reference domain can be defined such that it maps to the
		physical domain under the action of a smooth, invertible mapping
		$\func{\hat{\Gcal}}{\Rbb^d}{\Rbb^d}$, i.e., $\Omega_0 = \hat{\Gcal}^{-1}(\Omega)$.
		In this work, we take the reference and physical domains to be the same
		set, i.e., $\hat{\Gcal}=\mathrm{Id}$.
	\end{remark}
	\subsection{Discontinuous Galerkin discretization}
	\label{sec:govern:dg}
	We employ a nodal discontinuous Galerkin method \cite{cockburn_rungekutta_2001, hesthaven_nodal_2008} to discretize the transformed conservation law (\ref{eqn:claw-ref}). Here $\mathcal{E}_h$ represents a discretization of the reference domain $\Omega_0$ into distinct, possibly curved, non-overlapping computational elements. To establish the finite-dimensional DG formulation, we introduce the DG approximation space, consisting of discontinuous piecewise polynomials associated with the mesh $\mathcal{E}_h$
	\begin{equation}
		\Vcal_h^p := \left\{v \in [L^2(\Omega_0)]^m \suchthat
		\left.v\right|_K \in [\Pcal_p(K)]^m,
		~\forall K \in \Ecal_h\right\},
	\end{equation}
	where $\Pcal_p(K)$ is the space of polynomial functions of degree at most
	$p \geq 1$ on the element $K$. Furthermore, we define the space of globally
	continuous piecewise polynomials of degree $q$ associated with the mesh
	$\Ecal_h$ as
	\begin{equation} \label{eqn:gfcnsp}
		\Wcal_h := \left\{v \in C^0(\Omega_0) \suchthat
		\left.v\right|_K \in \Pcal_q(K),~\forall K \in \Ecal_h\right\}
	\end{equation}
	and discretize the domain mapping with the corresponding vector-valued
	space $[\Wcal_h]^d$.
	
	Considering the DG test space as $\mathcal{V}_h^{p'}$, where $p' \geq p$, the DG formulation can be expressed as follows: given $\Gcal_h\in[\Wcal_h]^d$, find
	$\bar{U}_h\in\Vcal_h^p$ such that for all $\bar\psi_h\in\Vcal_h^{p'}$ the following condition holds:
	\begin{equation} \label{eqn:claw-weak-elem}
		\int_{\partial K} \bar\psi_h^+ \cdot
		\bar\Hcal(\bar{U}_h^+,\bar{U}_h^-,N_h;\bar\nabla\Gcal_h) \, dS -
		\int_K \bar{F}(\bar{U}_h; \bar\nabla\Gcal_h):\bar\nabla \bar\psi_h \, dV =
		\int_K \bar\psi_h \cdot \bar{S}(\bar{U}_h; \det(\bar\nabla\Gcal_h)) \, dV,
	\end{equation}
	where $\func{N_h}{\partial K}{\Rbb^d}$ is the unit outward normal to
	element $K\in\Ecal_h$, $\bar{W}_h^+$ ($\bar{W}_h^-$) denotes the interior
	(exterior) trace of $\bar{W}_h$ to the element $K$ for $\bar{W}_h\in\Vcal_h^s$
	for any $s\in\Nbb$ (for $X\in\partial K\cap\partial\Omega_0$,
	$\bar{U}_h^-$ is a boundary state constructed to enforce the appropriate
	boundary condition). Moreover, $\func{\bar\Hcal}{\Rbb^m\times\Rbb^m\times\Rbb^d\times\Rbb^{d\times d}}{\Rbb^m}$ is the numerical flux function linked with the reference inviscid flux $\bar{F}$. This function is crucial as it ensures that the surface integrand remains single-valued and can be designed to guarantee that the DG discretization maintains properties of consistency, conservativeness, and stability \cite{hesthaven_nodal_2008}. The expression for the reference numerical flux function can be derived from the standard physical numerical flux function \cite{zahr_implicit_2020}. The residual form of the
	DG equation in (\ref{eqn:claw-weak-elem}) is given by
	$\func{r_h^{p',p}}{\Vcal_h^{p'}\times\Vcal_h^p\times[\Wcal_h]^d}{\Rbb}$
	\begin{equation}
		r_h^{p',p} : (\bar\psi_h,\bar{W}_h,\Gcal_h) \mapsto
		\sum_{K\in\Ecal_h} r_K^{p',p}(\bar\psi_h,\bar{W}_h,\Gcal_h),
	\end{equation}
	where the elemental DG form is given by
	$\func{r_K^{p',p}}{\Vcal_h^{p'}\times\Vcal_h^p\times[\Wcal_h]^d}{\Rbb}$
	\begin{equation}
		\begin{aligned}
			r_K^{p',p} : (\bar\psi_h,\bar{W}_h,\Gcal_h) \mapsto
			&\int_{\partial K} \bar\psi_h^+ \cdot
			\bar\Hcal(\bar{W}_h^+,\bar{W}_h^-,N_h;\bar\nabla\Gcal_h) \, dS \\
			&-\int_K \bar{F}(\bar{W}_h; \bar\nabla\Gcal_h):\bar\nabla \bar\psi_h \, dV \\
			&-\int_K \bar\psi_h \cdot \bar{S}(\bar{W}_h; \det(\bar\nabla\Gcal_h)) \, dV.
		\end{aligned}
	\end{equation}
	
	Next, we introduce a (nodal) basis for the test space ($\Vcal_h^{p'}$),
	trial space ($\Vcal_h^p$), and domain mapping space ($[\Wcal_h]^d$).
	This enables the transformation of the weak formulation into a system of nonlinear algebraic equations in residual form.
	In the case where $p'=p$, we denote the algebraic residual
	\begin{equation}
		\func{\rbm}{\Rbb^{N_\ubm}\times\Rbb^{N_\xbm}}{\Rbb^{N_\ubm}}, \qquad
		\rbm : (\ubm,\xbm) \mapsto \rbm(\ubm,\xbm),
	\end{equation}
	where $N_\ubm=\dim \Vcal_h^p$ and $N_\xbm=\dim[\Wcal_h]^d$. In this notation, a standard DG discretization in algebraic form can be expressed as follows: given $\ubm\in\Rbb^{N_\ubm}$ such that $\rbm(\ubm,\xbm)=\zerobold$, where $\ubm$ are the DG solution coefficients and $\xbm$ are the coefficients of the domain mapping (nodal coordinates). Typically, $\xbm$ is predetermined during mesh generation and remains fixed. However, in this work, it will be determined through optimization to ensure that the mesh accurately tracks and aligns the element faces with all flow discontinuities. Finally, we define the algebraic \textit{enriched residual}
	\begin{equation}
		\func{\Rbm}{\Rbb^{N_\ubm}\times\Rbb^{N_\xbm}}{\Rbb^{N_\ubm'}}, \qquad
		\Rbm : (\ubm,\xbm) \mapsto \Rbm(\ubm,\xbm)
	\end{equation}
	associated with a test space of degree $p'$, where $N_\ubm'=\dim\Vcal_h^{p'}$. The enriched residual will be used to construct the implicit shock tracking objective function and in this work we take $p'=p+1$. 
	
	\subsection{Sparsity of discrete operators}
	\label{sec:govern:sparse}
	Next, we examine the sparsity structure of the Jacobians of the DG residuals $\rbm$ and $\Rbm$ with respect to the variables $\ubm$ and $\xbm$, as they will be central to the shock tracking optimization method. For any element $K_e\in\Ecal_h$, let $\ubm_e \in \Rbb^{N_p}$ and $\xbm_e \in \Rbb^{N^e_\xbm}$ denote the degrees of freedom (DOFs) of $\ubm$ and $\xbm$, respectively, associated with element $K_e$, where $ N_p =m \dim[\Pcal_p(K_e)]$ and $N^e_\xbm=d \dim[\Pcal_q(K_e)]$. The elemental DOFs are related to the global DOFs via the selection matrices, $\Pbm_e\in\{0,1\}^{N_\ubm\times N_p}$ and $\Qbm_e\in\{0,1\}^{N_\xbm\times N_\xbm^e}$, which are subsets of the identity matrix that extract selected rows from $N_\ubm$- and $N_\xbm$-vectors, respectively,
	\begin{equation}
	\ubm_e = \Pbm_e^T\ubm, \qquad \xbm_e = \Qbm_e^T\xbm.
	\end{equation}
	Furthermore, denote the DOFs corresponding to the neighbors of element $K_e$ as $\hat\ubm_e\in \Rbb^{\hat N^e_p}$, where $\hat N^e_p =m\dim[\Pcal_p(\mathcal{N}_e)]$, and $\hat\Pbm_e\in\{0,1\}^{N_\ubm\times \hat{N}_p^e}$ as the corresponding selection matrix such that
	\begin{equation}
	 \hat\ubm_e = \hat\Pbm_e^T\ubm.
	\end{equation}
	Here, $\Ncal_e \subset \Ecal_h$ is the collection of elements neighboring (i.e., sharing a face with) element $K_e$.
	
\begin{figure}
			\centering
				\ifbool{fastcompile}{}{
			\includegraphics[height=55mm]{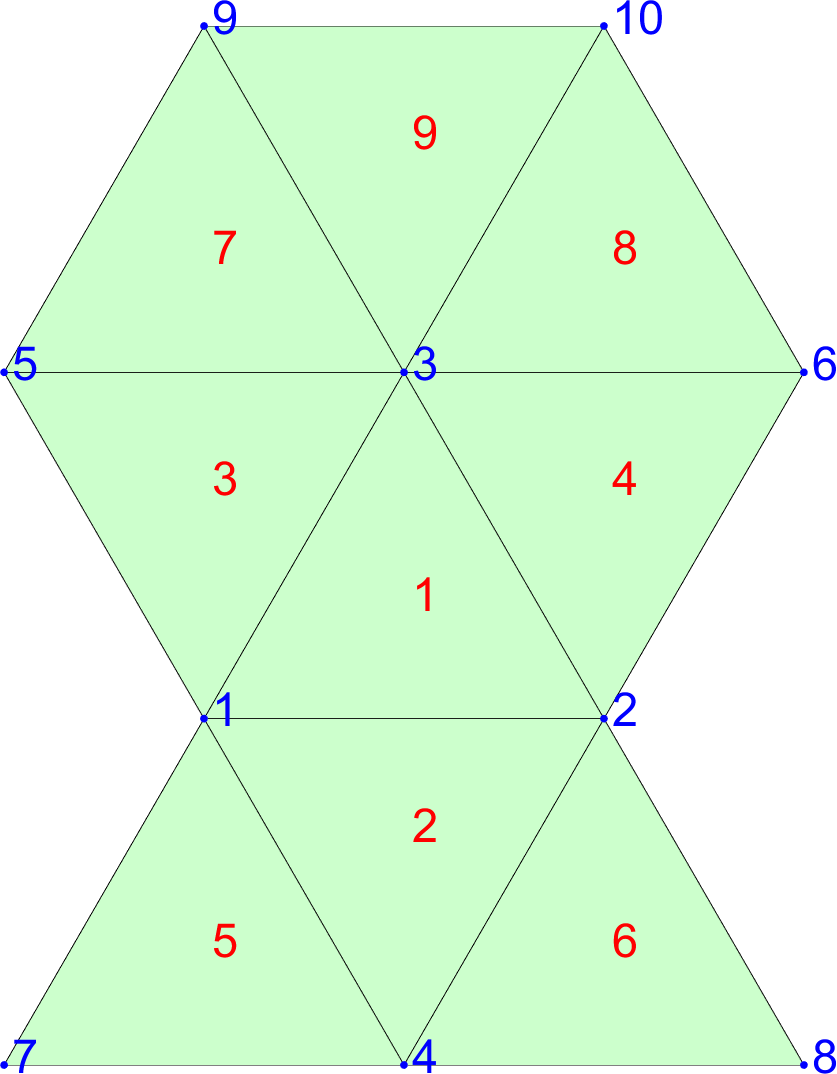} \quad \includegraphics[height=55mm]{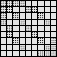}
				}
			\caption{Example two-dimensional mesh (\textit{left}) ($10$ nodes and $9$ elements) and corresponding sparsity structure of $\Ju$ (\textit{right}) for a polynomial degree of $p=1$ and a single conservation law ($m=1$). This choice leads to $9$ blocks of size $3 \times 3$ for \Ju.}
			\label{fig:Ju_sprs}
		\end{figure}

\begin{figure}
			\centering
				\ifbool{fastcompile}{}{
			\includegraphics[height=55mm]{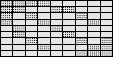} \qquad \qquad
				}
			\caption{Sparsity structure of $(\partial\Rbm/ \partial\ubm)^T$ for the mesh in Figure \ref{fig:Ju_sprs}, polynomial degrees of $p=1,~p'=2$ and a single conservation law ($m=1$). This choice results in $9$ blocks of size $6 \times 3$ for $\partial\Rbm/\partial\ubm$.}
			\label{fig:dR_sprs}
		\end{figure}

	With this notation, the elemental DG residuals, $r_{K_e}^{p,p}$ and $r_{K_e}^{p',p}$, can be written algebraically as
	\begin{eqnarray}
		\func{\rbm_e}{\Rbb^{N_p}\times\Rbb^{\hat{N}_p^e}\times\Rbb^{N^e_\xbm}}{\Rbb^{N_p}}, \qquad
		\rbm_e :  (\ubm_e,\hat\ubm_e,\xbm_e) \mapsto \rbm_e (\ubm_e,\hat\ubm_e,\xbm_e) \\
		\func{\Rbm_e}{\Rbb^{N_p}\times\Rbb^{\hat{N}_p^e}\times\Rbb^{N^e_\xbm}}{\Rbb^{N_{p'}}}, \qquad
		\Rbm_e : (\ubm_e,\hat\ubm_e,\xbm_e) \mapsto \Rbm_e(\ubm_e,\hat\ubm_e,\xbm_e).
	\end{eqnarray}
	The global residuals are formed by summing over all elements and assembling into the appropriate DOF as
	\begin{equation} \label{eqn:globres}
	 \rbm(\ubm,\xbm)= \sumEl \Pbm_e \rbm_e(\ubm_e, \hat\ubm_e, \xbm_e)= \sumEl \Pbm_e \rbm_e(\Pbm_e^T \ubm, \hat\Pbm_e^T \ubm, \Qbm_e^T \xbm).
	\end{equation}
	Direct differentiation leads to an expression for the Jacobian $\Ju(\ubm,\xbm)\in\Rbb^{N_\ubm \times N_\ubm}$ that exposes its block structure
	\begin{equation}
		\label{eqn:drdU-spars}
		\Ju(\ubm,\xbm) \coloneqq \pder{\rbm}{\ubm}(\ubm,\xbm) = \sumEl\Pbm_e\left(\pder{\rbm_e}{\ubm_e}(\ubm_e, \hat\ubm_e, \xbm_e)\Pbm_e^T + \pder{\rbm_e}{\hat\ubm_e}(\ubm_e, \hat\ubm_e, \xbm_e)\hat\Pbm_e^T\right),
	\end{equation}
	where $\pder{\rbm_e}{\ubm_{e}}(\ubm_e, \hat\ubm_e, \xbm_e)\in \Rbb^{N_p \times N_p}$,  $\pder{\rbm_e}{\hat\ubm_{e}}(\ubm_e, \hat\ubm_e, \xbm_e)\in \Rbb^{N_p \times \hat N^e_\ubm}$ are its matrix blocks. The matrix $\Ju$ is a $\vert \Ecal_h \vert \times \vert \Ecal_h \vert$ block matrix with blocks of size $N_p \times N_p$. For an example two-dimensional mesh consisting of $9$ elements (Figure \ref{fig:Ju_sprs}), a polynomial degree of $p=1$, and a single conservation law ($m=1$), the sparsity of $\Ju$ is shown in Figure \ref{fig:Ju_sprs}.				
	Using the same arguments, one can derive the sparsity pattern for the Jacobian $\partial \Rbm/ {\partial \ubm}$, the only difference being that the blocks $\pder{\Rbm_e}{\ubm_{e'}}\in \Rbb^{N_{p'}\times N_p}$ have a different size to account for the additional constraints. The sparsity pattern of $(\partial\Rbm /\partial{\ubm})^T$ for the exemplary mesh in Figure \ref{fig:Ju_sprs} is shown in Figure \ref{fig:dR_sprs} with $9$ blocks of size $3 \times 6$ coming from the enriched polynomial degree $p'=2$.
	\begin{remark}
		The elemental residuals do not depend on the neighboring nodes $\hat \xbm_e$ as the coupling to the neighboring elements is only due to trace values on the boundaries. For a more detailed discussion the reader is referred to the work by Wen et al. \cite{wenGloballyConvergentMethod2023}.
	\end{remark}
	
	\subsection{Nonlinear solvers}
	\label{sec:govern:nlsolve}
	Typically, in the context of a DG method without mesh adaptation (fixed mesh), one aims to solve the algebraic equation $\rbm(\ubm,\xbm)=0$ for $\ubm\in\Rbb^{N_\ubm}$ for fixed $\xbm\in\Rbb^{N_\xbm}$. This is a nonlinear system of equations, which is usually solved using a nonlinear iterative method, i.e., Newton's method or pseudo-transient continuation (PTC). For each iteration, the equations are linearized around the current iterate $\ubm_k$ and a linear system $\Abm \sbm =\bbm$ must be solved, where $\sbm=\Delta \ubm$, $\bbm = -\rbm(\ubm_k,\xbm)$, and $\Abm =\Ju(\ubm_k,\xbm)$ (Newton's method) or $\Abm =\Ju(\ubm_k,\xbm)+ (1/\Delta t) \Mbm$ (PTC), where $\Mbm$ is the mass matrix and $\Delta t\in\Rbb_{>0}$ is a pseudo time-step. For scenarios with many DOFs, direct solvers are not practical because of their large computational cost, memory footprint, and poor parallel scaling. For these situations, iterative linear solvers combined with efficient preconditioners are commonly used. In Section \ref{sec:linsolv:dg} we review matrix-based preconditioners that exploit the DG block structure and have been successfully used for linear systems arising from DG discretizations. They will serve as a basis for the preconditioners introduced in Section \ref{sec:linsolv:hoist} for implicit shock tracking.
	\begin{remark}
		The sparsity patterns of both the $\pder{\rbm}{\ubm}$ and $\pder{\Rbm}{\ubm}$ sub-blocks are influenced by the specific choice of the DG basis. If a nodal basis is used, dense blocks appear along the diagonal (representing the $K_e$-$K_e$ interaction) and sparse blocks off the diagonal (representing the $K_e$-$K_{e'}$ interaction). This sparsity arises because the elemental residuals are sensitive only to changes that occur at the common face of the $K_e$ and $K_{e'}$ elements.
	\end{remark}
	\section{High-order implicit shock tracking formulation and solver}
	\label{sec:hoist}
	In this section, we review the optimization formulation (Section \ref{sec:hoist:form}) and the sequential quadratic programming (SQP) solver (Sections \ref{sec:hoist:solve}-\ref{sec:hoist:hessian}) \cite{huangRobustHighorderImplicit2022a} on which the HOIST method is based, including the linear system that defines the SQP step. Finally, we investigate the sparsity structure of the SQP linear system (Section \ref{sec:hoist:sparsity}).
	
	\subsection{Formulation}
	\label{sec:hoist:form}	
	The HOIST method, as described in \cite{zahr_implicit_2020,huangRobustHighorderImplicit2022a}, is a high-order technique that simultaneously computes both the discrete solution of the conservation law and the nodal coordinates of the mesh, aligning element faces with discontinuities. This process is accomplished using a fully discrete, full-space optimization formulation, where the optimization variables consist of the discrete flow solution and the nodal coordinates of the mesh. 
	We begin the description of the HOIST method by introducing a boundary-preserving parameterization of the physical nodes (details of its construction can be found in \cite{huangRobustHighorderImplicit2022a})
	\begin{equation} \label{eqn:mapparam}
		\func{\phibold}{\Rbb^{N_\ybm}}{\Rbb^{N_\xbm}}, \qquad
		\phibold: \ybm \mapsto \phibold(\ybm),
	\end{equation}
	such that $\Gcal_h(\Ecal_h;\phibold(\ybm))$ conforms to $\partial\Omega$
	for any $\ybm\in\Rbb^{N_\ybm}$ that does not cause element inversion. With this parameterization of the mesh motion the HOIST method is formulated as
	\begin{equation} \label{eqn:pde-opt}
		(\ubm^\star,\ybm^\star) \coloneqq
		\argmin_{\ubm\in\Rbb^{N_\ubm},\ybm\in\Rbb^{N_\ybm}} f(\ubm,\ybm) \quad \text{subject to:} \quad \rbm(\ubm,\phibold(\ybm)) = \zerobold,
	\end{equation}
	where $\func{f}{\Rbb^{N_\ubm}\times\Rbb^{N_\ybm}}{\Rbb}$ is the objective
	function and the nodal coordinates of the aligned mesh are
	$\xbm^\star = \phibold(\ybm^\star)$. The objective function is composed
	of two terms as
	\begin{equation}\label{eqn:obj0}
		f : (\ubm,\ybm) \mapsto f_\text{err}(\ubm,\ybm) + \kappa^2 f_\text{msh}(\ybm),
	\end{equation}
	balancing alignment of the mesh with non-smooth features and the quality
	of the elements. Here, $\kappa\in\Rbb_{\geq 0}$ is an adaptively chosen mesh penalty parameter to weight the two terms such that the first term is prioritized \cite{huangRobustHighorderImplicit2022a}. The mesh alignment term, $\func{f_\text{err}}{\Rbb^{N_\ubm}\times\Rbb^{N_\ybm}}{\Rbb}$, is taken
	to be the norm of the enriched DG residual
	\begin{equation}\label{eqn:obj1}
		f_\text{err} : (\ubm,\ybm) \mapsto \frac{1}{2}\norm{\Rbm(\ubm,\phibold(\ybm))}_2^2.
	\end{equation}
	We also want to ensure that the elements of the
	discontinuity-aligned mesh are of high quality, which leads to the definition of the mesh distortion
	term, $\func{f_\text{msh}}{\Rbb^{N_\ybm}}{\Rbb}$, as
	\begin{equation}\label{eqn:obj2}
		f_\text{msh} : \ybm \mapsto \frac{1}{2}\norm{\Rbm_\text{msh}(\phibold(\ybm))}_2^2,
	\end{equation}
	where $\func{\Rbm_\text{msh}}{\Rbb^{N_\ybm}}{\Rbb^{|\Ecal_h|}}$ is the
	elemental mesh distortion with respect to an ideal element
	\cite{zahr_implicit_2020,knupp_algebraic_2001,roca_defining_2012}. 
	
	To obtain the first-order optimality system of the implicit shock
	tracking formulation (\ref{eqn:pde-opt}), we introduce the corresponding
	Lagrangian,
	$\func{\Lcal}{\Rbb^{N_\ubm}\times\Rbb^{N_\ybm}\times\Rbb^{N_\ubm}}{\Rbb}$,
	defined as
	\begin{equation}
		\Lcal : (\ubm,\ybm,\lambdabold) \mapsto f(\ubm,\ybm)-\lambdabold^T\rbm(\ubm,\phibold(\ybm)).
	\end{equation}
	Then, the first-order optimality, or Karush-Kuhn-Tucker (KKT), conditions state
	that $(\ubm^\star,\ybm^\star)\in\Rbb^{N_\ubm}\times\Rbb^{N_\ybm}$ is a
	first-order solution of the optimization problem in (\ref{eqn:pde-opt}) if there
	exists $\lambdabold^\star\in\Rbb^{N_\ubm}$ such that the Lagrangian is
	stationary, i.e.
	\begin{equation}\label{eqn:kkt0}
		\nabla\Lcal(\ubm^\star,\ybm^\star,\lambdabold^\star) =0.
	\end{equation}
	\subsection{Sequential quadratic programming solver}
	\label{sec:hoist:solve}
	Next, we briefly describe the SQP solver \cite{huangRobustHighorderImplicit2022a} for the optimization problem in (\ref{eqn:pde-opt}). It is a full-space approach that aims to converge the DG solution and the mesh to their optimal values simultaneously. To this end, we define a new variable $\zbm\in\Rbb^{N_\zbm}$ ($N_\zbm=N_\ubm+N_\ybm$) that combines the DG solution $\ubm$ and the unconstrained mesh coordinates $\ybm$ as
	\begin{equation}
		\zbm = (\ubm, \ybm),
	\end{equation}
	and use $\zbm$ interchangeably with $(\ubm,\ybm)$.
	For brevity, we introduce the following notation for the
	derivatives of the objective function,
	$\func{\gbm}{\Rbb^{N_\zbm}}{\Rbb^{N_\zbm}}$, and the DG residual,
	$\func{\Jbm}{\Rbb^{N_\zbm}}{\Rbb^{N_\ubm}\times\Rbb^{N_\zbm}}$, as
	\begin{equation}
		\gbm : \zbm \mapsto
		\begin{bmatrix}
			\ds{\pder{f}{\ubm}(\ubm,\ybm)^T} \\[1em]
			\ds{\pder{f}{\ybm}(\ubm,\ybm)^T}
		\end{bmatrix}, \quad
		\Jbm : \zbm \mapsto
		\begin{bmatrix}
			\ds{\pder{\rbm}{\ubm}(\ubm,\phibold(\ybm))} &
			\ds{\pder{\rbm}{\xbm}(\ubm,\phibold(\ybm))\pder{\phibold}{\ybm}(\ybm)}
		\end{bmatrix}.
	\end{equation}
	The SQP method in \cite{huangRobustHighorderImplicit2022a} produces a sequence of iterates
	$\{\zbm_k\}_{k=0}^\infty$ such that
	$\zbm_k=(\ubm_k,\ybm_k)\rightarrow\zbm^\star=(\ubm^\star,\ybm^\star)$,
	where $(\ubm^\star,\ybm^\star)$ satisfies the first-order
	optimality conditions in (\ref{eqn:kkt0}). The sequence of iterates
	is generated as
	\begin{equation} \label{eqn:sqp_step0}
		\zbm_{k+1} = \zbm_k + \alpha_k\Delta \zbm_k,
	\end{equation}
	where the search direction $\Delta\zbm_k\in\Rbb^{N_\zbm}$
	is computed as the solution of the following quadratic program
	\begin{equation}\label{eqn:qp0}
		\optconOne{\Delta\zbm \in \Rbb^{N_\zbm}}
		{\gbm_k^T\Delta\zbm+\frac{1}{2}\Delta\zbm^T\Bbm_k\Delta\zbm}
		{\rbm_k+\Jbm_k\Delta\zbm = \zerobold},
	\end{equation}
	$\gbm_k\in\Rbb^{N_\zbm}$, $\rbm_k\in\Rbb^{N_\ubm}$, and
	$\Jbm_k\in\Rbb^{N_\ubm\times N_\zbm}$ are the objective gradient,
	residual, and residual Jacobian, respectively, evaluated at $\zbm_k$
	\begin{equation}
		\rbm_k \coloneqq \rbm(\ubm_k,\phibold(\ybm_k)), \qquad
		\gbm_k \coloneqq \gbm(\zbm_k), \qquad
		\Jbm_k \coloneqq \Jbm(\zbm_k),
	\end{equation}
	$\Bbm_k\in\Rbb^{N_\zbm \times N_\zbm}$ is a symmetric positive
	definite (SPD) approximation to the Hessian of the Lagrangian at $\zbm_k$,
	and $\alpha_k\in\Rbb_{>0}$ is the step length. The latter is computed by an inexact line search employing a standard $l_1$-merit function \cite{nocedal2006numerical} and the first-order
	optimality conditions of the quadratic program lead to the following linear system
	of equations
	\begin{equation} \label{eqn:sqp_sys0}
		\begin{bmatrix}
			\Bbm_k & \Jbm_k^T \\
			\Jbm_k & \zerobold
		\end{bmatrix}
		\begin{bmatrix}
			\Delta\zbm_k \\ \etabold_k
		\end{bmatrix}
		=
		-\begin{bmatrix}
			\gbm_k \\ \rbm_k
		\end{bmatrix},
	\end{equation}
	where $\etabold_k\in\Rbb^{N_\ubm}$ are the Lagrange multipliers associated
	with the linearized constraint in (\ref{eqn:qp0}). This linear system of size $2N_\ubm + N_\xbm$ must to be solved at each iteration $k$ to compute the step $\Delta \zbm_k$ to update the DG solution and mesh (\ref{eqn:sqp_step0}). For large-scale applications with many DOFs, direct solvers are not a viable option as these systems are larger than standard DG system (size: $N_\ubm$). In the next two sections, we explore the structure of the linear system in (\ref{eqn:sqp_sys0}) that will facilitate the development of efficient preconditioners and iterative linear solvers in Section~\ref{sec:linsolv:hoist}.
	
\begin{remark}
This SQP method \cite{zahr_implicit_2020} proved not to be robust enough to handle complex problems, such as high Mach number flows with complex discontinuities, so several robustness measures were introduced \cite{huangRobustHighorderImplicit2022a}. These measures manipulate the state $\zbm_{k+1}$ after the SQP update only for a fixed number of iterations $M>0$ (to ensure SQP convergence in the limit) and include 1) boundary-preserving, shock-aware element removal, 2) geometric curvature removal from inverted or ill-conditioned elements, and 3) elemental solution reinitialization; see \cite{huangRobustHighorderImplicit2022a} for details. These operations have a small positive impact on linear solvers as they locally reduce sources of ill-conditioning, which can lead to abrupt (positive) changes in the performance of iterative solvers when comparing between different states (Section~\ref{sec:numexp}).
\end{remark}
		
	\subsection{Hessian approximation}
	\label{sec:hoist:hessian}
	Implicit shock tracking methods employ a Levenberg-Marquardt Hessian approximation introduced in \cite{corrigan2019moving, zahr_implicit_2020} to define $\Bbm_k$. To this end, $\Bbm_k$ is expanded as
	\begin{equation} \label{eqn:hess0}
		\Bbm_k =
		\begin{bmatrix}
			\Bbm_{\ubm\ubm,k} & \Bbm_{\ubm\ybm,k} \\[0.5em]
			\Bbm_{\ubm\ybm,k}^T & \Bbm_{\ybm\ybm,k}
		\end{bmatrix},
	\end{equation}
	where the individual components
	$\Bbm_{\ubm\ubm,k}\in\Rbb^{N_\ubm\times N_\ubm}$,
	$\Bbm_{\ubm\ybm,k}\in\Rbb^{N_\ubm\times N_\ybm}$, and
	$\Bbm_{\ybm\ybm,k}\in\Rbb^{N_\ybm\times N_\ybm}$
	are defined as
	\begin{equation} \label{eqn:hess1}
		\begin{aligned}
			\Bbm_{\ubm\ubm,k} &\coloneqq \pder{\Fbm}{\ubm}(\zbm_k)^T\pder{\Fbm}{\ubm}(\zbm_k) \\
			\Bbm_{\ubm\ybm,k} &\coloneqq \pder{\Fbm}{\ubm}(\zbm_k)^T\pder{\Fbm}{\ybm}(\zbm_k) \\
			\Bbm_{\ybm\ybm,k} &\coloneqq \pder{\Fbm}{\ybm}(\zbm_k)^T\pder{\Fbm}{\ybm}(\zbm_k)+
			\gamma_k\pder{\phibold}{\ybm}(\ybm_k)^T\Dbm_k\pder{\phibold}{\ybm}(\ybm_k),
		\end{aligned}
	\end{equation}
	where $\func{\Fbm}{\Rbb^{N_\ubm}\times\Rbb^{N_\ybm}}{\Rbb^{N_\ubm'+|\Ecal_h|}}$ is the residual function 
	\begin{equation} \label{eqn:residualF}
		\Fbm : (\ubm,\ybm) \mapsto \begin{bmatrix}\Rbm(\ubm,\phibold(\ybm)) \\ \kappa \Rbm_\text{msh}(\phibold(\ybm)) \end{bmatrix},
	\end{equation}
and $\Dbm_k\in\Rbb^{N_\xbm\times N_\xbm}$ is a SPD matrix constructed to
	regularize the mesh motion. The regularization
	parameter $\gamma_k\in\Rbb_{\geq 0}$ is chosen adaptively during the optimization process and has a strong impact on the number of iterations needed for an iterative solver, which will be observed in the numerical experiments in Section~\ref{sec:numexp:results:gamma}. 
		
	\subsection{Sparsity of the linearized optimality system}
	\label{sec:hoist:sparsity}
\begin{figure}
			\centering
				\ifbool{fastcompile}{}{
			\includegraphics[width=55mm]{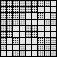} \quad \includegraphics[height=55mm]{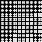} 
				}
			\caption{Sparsity structure of $\Bbm_{\ubm\ubm}$ (\textit{left}) and $\Bbm_{\ybm\ybm}$ (\textit{right}) (assuming no boundary constraints, i.e.  $\phibold(\ybm)=\ybm$) for mesh depicted in Figure \ref{fig:Ju_sprs} with polynomial degrees of $p=1,~p'=2$ and a single conservation law ($m=1$).}
			\label{fig:Buu_sprs}
		\end{figure}

	In this section, we detail the sparsity of the linear system \eqref{eqn:sqp_sys0} as it has significant implications for the design requirements of efficient preconditioners. From this point forward, we fix the state $\zbm_k$ and drop the $k$ subscript on all terms.
	The sparsity of $\Ju$ was already examined in Section~\ref{sec:govern:sparse} so we begin with $\Bbm_{\ubm\ubm}$. First, recall that  $\partial\Rbm_{\text{msh}}/{\partial \ubm}=0$ by construction. Therefore, we can derive the block sparsity structure using the elemental decomposition of the DG Jacobian as follows
	\begin{eqnarray}
		\Bbm_{\ubm\ubm} = \pder{\Fbm}{\ubm}^T \pder{\Fbm}{\ubm}=\pder{\Rbm}{\ubm}^T \pder{\Rbm}{\ubm} = &\left(\sum_{e=1}^{|\Ecal_h|} \Pbm_e'\left(\pder{\Rbm_e}{\ubm_e}\Pbm_e^T + \pder{\Rbm_e}{\hat\ubm_e}\hat\Pbm_e^T\right)\right)^T \left(\sum_{E=1}^{|\Ecal_h|} \Pbm_E'\left(\pder{\Rbm_E}{\ubm_E}\Pbm_E^T + \pder{\Rbm_E}{\hat\ubm_E}\hat\Pbm_E^T\right)\right) \\
		= &\sum_{e,E=1}^{|\Ecal_h|} \left(\left(\Pbm_e\pder{\Rbm_e}{\ubm_e}^T + \hat\Pbm_e\pder{\Rbm_e}{\hat\ubm_e}^T\right)(\Pbm_e')^T \Pbm_E'\left(\pder{\Rbm_E}{\ubm_E}\Pbm_E^T + \pder{\Rbm_E}{\hat\ubm_E}\hat\Pbm_E^T\right)\right). \nonumber
	\end{eqnarray}
	Due to the fact that $(\Prol_e')^T\Prol_E' = \Ibm \delta_{eE}$ (because $\Rbm$ is a DG residual), we finally obtain 
	\begin{equation}
		\Bbm_{\ubm\ubm} = \sumEl \left(\Pbm_e\pder{\Rbm_e}{\ubm_e}^T\pder{\Rbm_e}{\ubm_e}\Pbm_e^T+ \Pbm_e\pder{\Rbm_e}{\ubm_e}^T\pder{\Rbm_e}{\hat\ubm_e}\hat\Pbm_e^T + \hat\Pbm_e\pder{\Rbm_e}{\hat\ubm_e}^T\pder{\Rbm_e}{\ubm_e}\Pbm_e^T + \hat\Pbm_e\pder{\Rbm_e}{\hat\ubm_e}^T\pder{\Rbm_e}{\hat\ubm_e}\hat\Pbm_e^T\right).
	\end{equation}
		From this identity it can be deduced that $\Bbm_{\ubm\ubm}$ has an element based block structure like $\Ju$, but with an extended (denser) sparsity pattern. The difference lies in the additional non-zero blocks due to the neighbor-neighbor interaction $\hat\Pbm_e\pder{\Rbm_e}{\hat\ubm_e}^T\pder{\Rbm_e}{\hat\ubm_e}\hat\Pbm_e^T$, which does not exist for $\Ju$. This amounts to an increase of non-zero blocks by a factor of $1+d(d+1)/(d+2)$ for a simplicial mesh (see Remark~\ref{remark:factorBuu}), which is expensive and memory-intensive (Table~\ref{tab:mul}) and requires parallel communication of the blocks to form the product. Therefore, in the next section, we will avoid preconditioners that require explicitly forming $\Bbm_{\ubm\ubm}$. On the other hand, matrix-vector products of the form $\Bbm_{\ubm\ubm}\vbm$ can be performed efficiently as $\pder{\Rbm}{\ubm}^T(\pder{\Rbm}{\ubm} \vbm)$, making it well-suited for use with an iterative (Krylov) solver. For illustrative purposes, the sparsity of $\Bbm_{\ubm\ubm}$ is shown in Figure \ref{fig:Buu_sprs} (\textit{left}) for the same exemplary mesh as in Figure \ref{fig:Ju_sprs}, and the substantial decrease in sparsity can be observed. 

	\begin{remark}\label{remark:factorBuu}
	Let us quantify the sparsity of $\Ju$ relative to $\Bbm_{\ubm\ubm}$ for simplicial grids. Let $m_1$ and $m_2$ denote the number of non-zero blocks per row of $\Ju$ and $\Bbm_{\ubm\ubm}$, respectively. For simplicity, we consider a row corresponding to an element sufficiently far from a boundary to avoid enumerating special cases. Because the $e$th block row of $\Ju$ has a non-zero for each element neighboring $K_e$, we have $m_1 = d+2$ (the block diagonal plus $d+1$ neighbors). On the other hand, the $e$th block row of $\Bbm_{\ubm\ubm}$ has a non-zero for all neighbors of $K_e$ and all neighbors of neighbors of $K_e$, which gives $m_2 = m_1 + d(d+1)$ (each of the $d+1$ neighbors of $K_e$ adds at most $d$ new neighbors). Thus, the ratio of non-zero blocks in $\Bbm_{\ubm\ubm}$ to those in $\Ju$ is
	\begin{equation} \label{eqn:mul}
		 \frac{m_2}{m_1} = 1 + \frac{d(d+1)}{d+2},
	\end{equation}
which is a significant factor (Table~\ref{tab:mul}), especially considering the DG Jacobians themselves are already memory-intensive to form and store. This motivates our decision to avoid explicitly forming $\Bbm_{\ubm\ubm}$ in the proposed preconditioners in Section~\ref{sec:linsolv:hoist}.
	\begin{table}
	\caption{Growth of block sparsity structure of $\Bbm_{\ubm\ubm}$ ($m_2$) relative to $\Ju$ ($m_1$).}
	\label{tab:mul}
	\centering
	\begin{tabular}{c|cccc}
		$d$ & $1$ & $2$ & $3$ & $4$ \\\hline
		$m_2/m_1$ & $1.67$ & $2.5$ & $3.4$ & $4.33$ \\
	\end{tabular}
\end{table}
	\end{remark}
	
\begin{figure}
			\centering
				\ifbool{fastcompile}{}{
			\includegraphics[width=120mm]{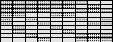} \qquad \quad
				}
			\caption{Sparsity structure of $(\partial\Rbm /{\partial\xbm})^T$ for mesh depicted in Figure \ref{fig:Ju_sprs} with polynomial degrees $p=1,~p'=2$ and a single conservation law ($m=1$).}
			\label{fig:dRdx_sprs}
		\end{figure}

	Next, we examine the sparsity of $\Bbm_{\ubm\ybm}$ and $\Bbm_{\ybm \ybm}$. First, we build the DG Jacobian with respect to $\ybm$ using the Jacobian of the $\phibold$ mapping
	\begin{equation}
		\label{eqn:dRdY-chainrule}
		\pder{}{\ybm}\Rbm(\ubm,\phibold(\ybm))=\pder{\Rbm}{\xbm}(\ubm,\phibold(\ybm)) \pder{\phibold}{\ybm}(\ybm)=\pder{\Rbm}{\xbm}(\ubm,\xbm) \pder{\phibold}{\ybm}(\ybm).
	\end{equation}
	Direct differentiation of (\ref{eqn:globres}), replacing $\rbm$ with $\Rbm$ and with respect to $\xbm$ exposes the assembled block structure of the DG residual
	\begin{equation}
		\label{eqn:dRdY-spars}
		\pder{\Rbm}{\xbm}(\ubm,\xbm) = \sumEl\Pbm_e'\pder{\Rbm_e}{\xbm_e}(\Pbm_e^T\ubm, \hat\Pbm_e^T\ubm, \Qbm_e^T\xbm) \Qbm_e^T=\sumEl \Pbm_e'\pder{\Rbm_e}{\xbm_e}(\ubm_e, \hat\ubm_e, \xbm_e) \Qbm_e^T,
	\end{equation}
	where $\pder{\Rbm_e}{\xbm_e} \in \Rbb^{N_{p'} \times N_\xbm^e}$. Note that because most nodes $\xbm_e$ are shared between two elements, it is not possible to obtain an elemental block structure for the columns (elemental block structure does exist for the rows). For illustrative purposes only, we assume no boundary constraints ($\phibold(\ybm)=\ybm$) and refer to Figure~\ref{fig:Ju_sprs} for our exemplary mesh. The sparsity pattern of $(\partial\Rbm/\partial \xbm)^T$ is illustrated in Figure \ref{fig:dRdx_sprs}, revealing $9$ block rows. 
	
	Going further, we obtain the following sparsity-block structure for $\Bbm_{\ubm\ybm}$ from the following identity
	\begin{align}
		\label{eqn:BUY-spars-deriv}
		\Bbm_{\ubm\ybm}&=\pder{\Fbm}{\ubm}^T\pder{\Fbm}{\ybm}=\pder{\Rbm}{\ubm}^T \pder{\Rbm}{\xbm}\pder{\phibold}{\ybm} = \left(\sumEl\Pbm_e'\left(\pder{\Rbm_e}{\ubm_e}\Pbm_e^T + \pder{\Rbm_e}{\hat\ubm_e}\hat\Pbm_e^T\right)\right)^T\left(\sum_{E=1}^{|\Ecal_h|}\Pbm_E'\pder{\Rbm_E}{\xbm_E}\Qbm_E^T\right) \pder{\phibold}{\ybm} \\&= \sumEl \left(\Pbm_e\pder{\Rbm_e}{\ubm_e}^T + \hat\Pbm_e\pder{\Rbm_e}{\hat\ubm_e}^T\right)\left(\pder{\Rbm_e}{\xbm_e}\Qbm_e^T\right)\pder{\phibold}{\ybm}. \nonumber
	\end{align}
	The structure of $\partial \rbm/\partial \ybm$ is identical to that of $\partial \Rbm/ \partial \ybm$ by repeating the above derivation. Because this is a (rectangular) off-diagonal term, the proposed preconditioners and linear solvers only require products with $\Bbm_{\ubm\ybm}$, which can be computed as $\pder{\Rbm}{\ubm}^T(\pder{\Rbm}{\xbm}(\pder{\phibold}{\ybm} \vbm))$ for any vector $\vbm$, so $\Bbm_{\ubm\ybm}$ never needs to be explicitly computed.
	
Lastly, we consider the structure of $\Bbm_{\ybm\ybm}$. From a simple application of the chain rule, we have
\begin{equation}
	\Bbm_{\ybm\ybm} = \pder{\phibold}{\ybm}^T \Bbm_{\xbm\xbm} \pder{\phibold}{\ybm},
\end{equation}
where
\begin{equation}
 \Bbm_{\xbm\xbm}= \pder{\Rbm}{\xbm}^T\pder{\Rbm}{\xbm}+ \kappa^2\pder{\Rbm_{\text{msh}}}{\xbm}^T\pder{\Rbm_{\text{msh}}}{\xbm} + \gamma \Dbm.
\end{equation}
Furthermore, from (\ref{eqn:dRdY-spars}), we have
	\begin{equation}
		\label{eqn:BYY-spars-deriv}
		\pder{\Rbm}{\xbm}^T\pder{\Rbm}{\xbm} = \left(\sumEl\Pbm_e'\pder{\Rbm_e}{\xbm_e}\Qbm_e^T\right)^T\left(\sum_{E=1}^{|\Ecal_h|}\Pbm_E'\pder{\Rbm_E}{\xbm_E}\Qbm_E^T\right) =\sumEl \left(\Qbm_e\pder{\Rbm_e}{\xbm_e}^T\pder{\Rbm_e}{\xbm_e}\Qbm_e^T\right),
	\end{equation}
which loses block structure once assembled because of the overlapping entries in $\Qbm_e$ for different elements.
The sparsity of $\pder{\Rbm_{\text{msh}}}{\xbm}^T\pder{\Rbm_{\text{msh}}}{\xbm}$ is a subset of $\pder{\Rbm}{\xbm}^T\pder{\Rbm}{\xbm}$ because each entry of the mesh distortion $\Rbm_{\text{msh}}$ is defined individually for each element $K_e$ and solely depends on the element nodes $\xbm_e$. The sparsity of the regularization matrix, $\Dbm$, depends solely on its specific choice. We choose $\Dbm$ as the linear elasticity (isotropic) stiffness matrix, with the elasticity modulus being inversely proportional to the volume of elements in the reference mesh. Thus, the sparsity of $\Dbm$ is a subset of $\pder{\Rbm}{\xbm}^T\pder{\Rbm}{\xbm}$, as it originates from the continuous finite element discretization of the elasticity equations. Finally, the mapping $\phibold$ determines the final structure of $\Bbm_{\ybm \ybm}$.

	\section{Iterative linear solvers and preconditioners}
	\label{sec:linsolv}
In this section, we introduce preconditioners for implicit shock tracking linearized systems, which are derived from successful preconditioners utilized for DG methods. We begin with a brief overview of Krylov iterative solvers (Section~\ref{sec:linsolv:krylov}) and review commonly used preconditioners for DG discretization (Section~\ref{sec:linsolv:dg}). Finally, we present the novel preconditioners for implicit shock tracking (Section~\ref{sec:linsolv:hoist}).

	\subsection{Krylov solvers and preconditioning}
	\label{sec:linsolv:krylov}
	In this work, we consider Krylov subspace methods for solving the linear system $\Abm\sbm=\bbm$. Krylov methods only require the action of the matrix $\Abm$ on vectors, not the entire matrix itself, which minimizes storage cost. This is particularly advantageous for implicit shock tracking because it allows us to avoid explicitly forming all blocks of $\Bbm$ (the regularization Lagrangian Hessian approximation).
	
On the other hand, Krylov methods rely on preconditioning, i.e., transformation of the system $\Abm \sbm=\bbm$ to enhance its spectral properties. Left preconditioning is achieved by multiplying the linear system on the left by some non-singular matrix $\tilde\Abm^{-1}$ to yield
\begin{equation}
	\PrecInv \Abm \sbm=\PrecInv \bbm,
\end{equation}
which has the same solution as the original system. Here, $\Prec \approx \Abm$ is the preconditioner and must be inexpensive to apply its inverse to a vector $\Prec^{-1} \vbm$ to be practical. Generally, as the preconditioner approaches the original matrix $\Abm$, the number of Krylov iterations decreases while the associated costs increase (only one iteration is necessary if $\Prec =\Abm$). Finding a suitable preconditioner that balances the need for fewer Krylov iterations with increased costs per iteration requires a specialized solution tailored to the matrix structure, discretization method, and equations at hand. The most effective preconditioners require all or part of the matrix $\Abm$, which partially neutralizes the matrix-free benefits of Krylov methods. In Section~\ref{sec:linsolv:hoist}, we will develop matrix-based preconditioners for implicit shock tracking that (1) build on established preconditioners for the DG system (Section~\ref{sec:linsolv:dg}) and (2) avoid forming the entire $\Bbm$ matrix.

	\subsection{Preconditioners for discontinuous Galerkin methods}
	\label{sec:linsolv:dg}
	Two established matrix-based preconditioners for the DG system ($\Abm=\Ju,\bbm=-\rbm(\ubm)$) introduced in \cite{persson_newton-gmres_2008} are the block Jacobi preconditioner (Section \ref{sec:linsolv:dg:ilu0}) and the block incomplete LU (ILU) preconditioner with minimum discarded fill (MDF) (Section \ref{sec:linsolv:dg:ilu0}). Both preconditioners utilize the block structure of the Jacobian matrix $\Ju$ and are efficient in terms of computational cost and memory to form and apply. They will be building blocks for HOIST preconditioners.
	
	\subsubsection{Block Jacobi}
	\label{sec:linsolv:dg:jacobi}
	The block Jacobi preconditioner is obtained by setting all blocks of the original matrix $\Abm$ off the diagonal to zero, which can be written compactly as
	\begin{equation}
		\JBJ \coloneqq \sumEl\Pbm_e\pder{\rbm_e}{\ubm_e}\Pbm_e^T.
	\end{equation}
	This block diagonal preconditioner can be easily formed from $\Ju$ and its inverse can be explicitly formed by inverting each $N_p\times N_p$ block as
	\begin{equation}
		\JBJ^{-1} = \sumEl\Pbm_e\left(\pder{\rbm_e}{\ubm_e}\right)^{-1}\Pbm_e^T.
	\end{equation}
	Because the size of each block is relatively small, a direct solver can be used. According to \cite{persson_newton-gmres_2008}, this preconditioner shows good performance in specific cases, but loses effectivity as the Reynolds number or timestep increases, and in the low Mach limit.
	
		\subsubsection{Block incomplete LU preconditioning with minimum discarded fill reordering}
	\label{sec:linsolv:dg:ilu0}
	A more advanced preconditioner is the Incomplete Block LU Factorization (BILU) with Minimum Discarded Fill (MDF), which is achieved by performing an ILU0 factorization of the matrix $\Abm$ on the block level. This procedure involves limiting a standard LU factorization to maintain the sparsity structure of $\Abm$, i.e., any operation that would introduce new non-zero blocks (known as ``fill in'') are skipped. To optimize the performance of an ILU, it is augmented with an initial re-ordering of the matrix block rows to minimize fill-in. Readers are referred to \cite{persson_newton-gmres_2008} for the complete algorithm and implementation details.
	
	The preconditioner $\JBILU$ is formed as $\tilde\Pbm\JBILU=\tilde \Lbm \tilde \Ubm$, where $\tilde\Pbm$ is the MDF reordering permutation,  $\tilde\Lbm$ is a lower block-triangular matrix with the identity matrix along the diagonal, and $\tilde\Ubm$ is an upper block-triangular matrix; both $\tilde\Lbm$ and $\tilde\Ubm$ that share the same sparsity pattern as $\Ju$. Because of the complementary structure of $\tilde\Lbm$ and $\tilde\Ubm$, the matrix $\Abm$ can be mutated in-place into $\tilde\Lbm$ (strict lower block triangle) and $\tilde\Ubm$ (upper block triangle). To apply the inverse of $\JBILU$ to a vector $\wbm$ ($\JBILU^{-1}\wbm$), we must solve the system $\JBILU \vbm = \wbm$. First, we multiply this equation by the permutation and substitute the ILU factorization
	\begin{equation}
	 \tilde\Pbm\JBILU\vbm = \tilde\Lbm\tilde\Ubm\vbm = \tilde\Pbm\wbm.
	\end{equation}
Then, we apply the usual forward-backward substitution process to solve for $\vbm$: first solve $\tilde\Lbm\tilde\vbm = \tilde\Pbm\wbm$ for $\tilde\vbm$ using block forward substitution, then solve $\tilde\Ubm \vbm = \tilde\vbm$ for $\vbm$ using block backward substitution. Because the block diagonal of $\tilde\Lbm$ are identity matrices, forward substitution only requires matrix-vector products at the element level. On the other hand, backward substitution requires solving linear systems of size $N_p\times N_p$, which is usually performed with a direct solver because of the relatively small size. According to \cite{persson_newton-gmres_2008}, this preconditioner works effectively for a wide range of problems, particularly when combined with $p$-multigrid.

	\subsection{Preconditioners for implicit shock tracking}
	\label{sec:linsolv:hoist}
	In this section, we introduce matrix-based preconditioners tailored for the HOIST linearized system in \eqref{eqn:sqp_sys0} (Section~\ref{sec:linsolv:hoist:blkconprec}). These preconditioners are derived from \textit{constrained preconditioners}, commonly employed for linear systems encountered in constrained optimization (Section~\ref{sec:linsolv:hoist:conprec}). We close the section with a summary of all preconditioners proposed and studied in this work (Section~\ref{sec:linsolv:hoist:consideredprec}). We are interested in efficient preconditioners that do not require formation of $\Bbm_{\ubm\ubm}$ or involve the inverse of $\Ju$ or $\Bbm_{\ybm\ybm}$; however, we consider a suite of preconditioners to study what is lost by these requirements.
	
	\subsubsection{Constrained preconditioners}
	\label{sec:linsolv:hoist:conprec}
	The system matrix, which must be solved at every iteration of the HOIST method, repeated here for reference
		\begin{equation} \label{eqn:sqp_mat}
		\Abm=\begin{pmatrix}
			\Bbm & \Jbm^T \\
			\Jbm & \zerobold
		\end{pmatrix},
	\end{equation}
is a symmetric saddle-point matrix. Typically, matrices of this type are known to suffer from bad condition numbers and there exists a wide variety of preconditioners tailored to the specific scenarios where they arise \cite{benzi_numerical_2005}.  In the realm of constrained optimization, where these saddle-point systems naturally emerge from first-order optimality conditions, a class of popular preconditioners known as \textit{constrained preconditioners}, denoted $\KKTC$, is commonly used
\begin{equation}
	\label{eqn:approx_con_prec}
	\KKTC=\begin{pmatrix}
		\Gbm  &\tilde{\Jbm}^T \\
		\tilde{\Jbm} &0
	\end{pmatrix}.
\end{equation}
Here, $\Gbm\approx \Bbm$ and $\tilde{\Jbm} \approx \Jbm$ are approximations to the Hessian and constraint matrices. If $\tilde{\Jbm}=\Jbm$, the preconditioner is the coefficient matrix for a modified saddle-point problem with the same linearized constraint.  Furthermore, $\Gbm$ and $\tilde{\Jbm}$ are generally chosen as such that $\KKTC$ and $\Gbm$ are invertible and $\Gbm^{-1}$, $\tilde{\Jbm} \Gbm^{-1} \tilde{\Jbm}^T$ are easy to compute. In this case, the inverse of $\KKTC$ can be explicitly computed as
			\begin{equation}
	\label{equation:CP_inv}
	\begin{pmatrix}
		\Gbm  &\tilde{\Jbm}^T \\
		\tilde{\Jbm} &0
	\end{pmatrix}^{-1} = \begin{pmatrix}
		\Ibm & -\Gbm^{-1}\tilde{\Jbm}^T \\
		0 & \Ibm
	\end{pmatrix}\begin{pmatrix}
		\Gbm^{-1} & 0 \\
		0 & -(\tilde{\Jbm} \Gbm^{-1} \tilde{\Jbm}^T)^{-1}
	\end{pmatrix}\begin{pmatrix}
		\Ibm & 0 \\
		-\tilde{\Jbm} \Gbm^{-1} & \Ibm
	\end{pmatrix}.
\end{equation}
However, we consider more restrictive approximations because of our desire to avoid formation of $\Bbm_{\ubm\ubm}$ and inverses of $\Ju$ and $\Bbm_{\ybm\ybm}$.

\subsubsection{Block anti-triangular constrained preconditioner}
\label{sec:linsolv:hoist:blkconprec}
We propose a class of preconditioners for the HOIST linearized system with
\begin{equation} \label{equation:G_block_triag}
\Gbm =
	\begin{pmatrix}
		0 & 0  \\
		0 & \tilde\Bbm_{\ybm \ybm}
	\end{pmatrix}, \qquad
\tilde\Jbm = \begin{pmatrix} \tilde{\Jbm}_\ubm & \Jx \end{pmatrix},
\end{equation}
where $\tilde\Bbm_{\ybm\ybm}$ is an approximation to $\Bbm_{\ybm\ybm}$ and
$\tilde{\Jbm}_\ubm$ is an approximation to $\Ju$. Substitution of these choices
into the constrained preconditioner leads to a lower block anti-triangular matrix, denoted
$\PrecV{AT}$,
\begin{equation}
	\label{equation:BTCP}
	\PrecV{AT} = \begin{pmatrix}
		0 & 0 &\tilde{\Ju}^T \\
		0 & \tilde{\Bbm}_{\ybm \ybm} & \Jx^T \\
		\tilde{\Ju} &  \Jx &0
	\end{pmatrix}
\end{equation}
that will be referred to as the \textit{approximate block anti-triangular constrained preconditioner} in the remainder. The
inverse of $\PrecV{AT}$ is
\begin{equation}
	\PrecV{AT}^{-1}= \begin{pmatrix}
		\tilde\Cbm\tilde\Bbm_{\ybm\ybm}^{-1}\tilde\Cbm^T & \tilde\Cbm\tilde\Bbm_{\ybm\ybm}^{-1} &\tilde\Ju^{-1} \\
		\tilde\Bbm_{\ybm\ybm}^{-1}\tilde\Cbm^T & \tilde\Bbm_{\ybm\ybm}^{-1} & 0\\
		\tilde\Ju^{-T} &  0&0
	\end{pmatrix},
\end{equation}
where $\tilde\Cbm := -\tilde\Ju^{-1}\Jx$. Furthermore, the action of $\PrecV{AT}^{-1}$ on a vector
$\vbm=\begin{pmatrix} \vbm_1 & \vbm_2 & \vbm_3\end{pmatrix}^T$, is
		\begin{equation} \label{eqn:ATtimesb}
	\PrecV{AT}^{-1}\begin{pmatrix} \vbm_1 \\\vbm_2 \\\vbm_3\end{pmatrix} =\begin{pmatrix}
		\tilde\Ju^{-1}(-\Jx\tilde\Bbm_{\ybm\ybm}^{-1}(-\Jx^T\tilde\Ju^{-T} \vbm_1+ \vbm_2) +\vbm_3)\\
		\tilde\Bbm_{\ybm\ybm}^{-1}(-\Jx^T\tilde\Ju^{-T} \vbm_1+ \vbm_2)\\
		\tilde\Ju^{-T} \vbm_1
	\end{pmatrix}.
\end{equation}
The current arrangement of the $\PrecV{AT}^{-1}\vbm$ shows the following sequence of operation is required to compute the product: 1) a linear solve of the form $\tilde\Ju^T \wbm_1 = \vbm_1$, 2) matrix product of the form $\tilde\wbm_2 = \Jx^T\wbm_1$, 3) a linear solve of the form $\tilde\Bbm_{\ybm\ybm} \wbm_2 = -\tilde\wbm_2 + \vbm_2$, 4) matrix product of the form $\tilde\wbm_3 = \Jx \wbm_2$, and 5) a linear solve of the form $\tilde\Ju\wbm_3 = -\tilde\wbm_3 + \vbm_3$. Hence, the product $\PrecV{AT}^{-1}\vbm$ requires three linear system solves with the matrices $\tilde\Ju$, $\tilde\Bbm_{\ybm\ybm}$, and $\tilde\Ju^T$, and two matrix products with $\Jx$ and $\Jx^T$. Thus, preconditioners of this form entirely circumvent the need to form $\Bbm_{\ubm\ubm}$ or invert $\Ju^T\Ju$, $\Jx^T\Jx$, $\Bbm_{\ubm\ubm}$. 

\subsubsection{Considered preconditioners}
\label{sec:linsolv:hoist:consideredprec}
The effectivity and cost of the anti-triangular constrained preconditioner are determined by the approximations $\tilde\Ju$ and $\tilde\Bbm_{\ybm\ybm}$. In this work, we consider three choices of $\tilde\Ju$, including the standard DG preconditioners:
(1) $\tilde\Ju = \Ju$, (2) $\tilde\Ju = \JBJ$, and (3) $\tilde\Ju = \JBILU$. We also study
three choices for $\tilde\Bbm_{\ybm\ybm}$, including standard preconditioners for general sparse matrices: (1) $\tilde\Bbm_{\ybm\ybm} = \Bbm_{\ybm\ybm}$, (2) $\tilde\Bbm_{\ybm\ybm} = \text{diag}(\Bbm_{\ybm\ybm})$ (point Jacobi), and (3) $\tilde\Bbm_{\ybm\ybm} = \text{ilu}(\Bbm_{\ybm\ybm})$ (point ILU0), where $\text{diag}(\Bbm)$ extracts the diagonal of $\Bbm$ and $\text{ilu}(\Bbm)$ is the ILU0 factorization of $\Bbm$. The combinations of these choices studied in this work are summarized in Table~\ref{tab:all_preconditioners}, including preconditioners combined with $p$-multigrid (Section~\ref{sec:mulgrd:hoist}). The preconditioner $\Abm_0$ that uses $\tilde\Ju = \Ju$ and $\tilde\Bbm_{\ybm\ybm} = \Bbm_{\ybm\ybm}$ is not practical as the action of the preconditioner inverse to a vector will involve linear solves with $\Ju$, $\Ju^T$, and $\Bbm_{\ybm\ybm}$; however, it is included in our study as a benchmark for comparison, representing a best-case scenario in terms of iterative solver iterations.

\subsection{$p$-Multigrid for implicit shock tracking}
\label{sec:mulgrd:hoist}

In the context of the DG method, several studies have utilized $p$-multigrid techniques. These techniques are employed either as stand-alone methods to iteratively solve the linear system $\Abm \sbm=\bbm$ \cite{fidkowski_p-multigrid_2005} or as preconditioners for iterative solvers like GMRES \cite{persson_newton-gmres_2008}. The term \textit{$p$-multigrid} refers to a multi-level approach combined with a smoother $\tilde \Abm$ where the high-order linear system (for instance, when $p>2$) and the current iterate $\sbm$ are projected onto spaces of lower polynomial order. On the fine levels, the solution is smoothed using an operation of the form $\sbm \leftarrow \sbm + \tilde\Abm^{-1}(\bbm-\Abm \sbm)$. On the coarsest level (typically with $p=0$ or $p=1$), the linear system is solved exactly. 

Following \cite{persson_newton-gmres_2008}, we employ a two-level $p$-multigrid strategy for the linearized HOIST system $\KKTMat \sbm =\bbm$. In this approach, we first restrict the state variables $(\ubm,\ybm)$ to the coarse scale. Specifically, $\ubm$ will be restricted to a piecewise constant solution ($p=0$), and $\ybm$ is constrained to a mesh with straight-sided elements ($q=1$). Upon returning to the finer level through prolongation, a smoothing operation $\tilde \Abm$ is applied; for this, one of the preconditioners outlined in Section \ref{sec:linsolv:hoist:consideredprec}. This entire process is interpreted as an operator $\tilde\Abm_{p0}^{-1}$ approximating $ \Abm^{-1}$ and is employed as a preconditioner for a Krylov solver.

\begin{figure}
		\centering
			\ifbool{fastcompile}{}{
		\includegraphics[width=50mm]{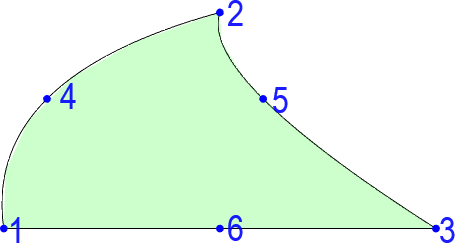} 
		\quad
		\includegraphics[width=50mm]{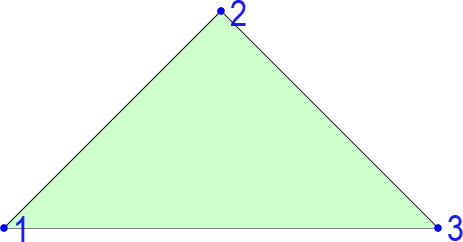}
		\quad
		\includegraphics[width=50mm]{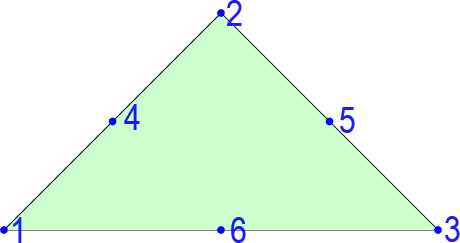}
	}
		\caption{Example of mesh restriction/prolongation for a second order mesh ($q=2$) with one element (\textit{left}). The original element is restricted to $q=1$ (\textit{middle}) removing the high order nodes $4,5,6$. Prolongation (\textit{right}) is performed by inserting high order nodes interpolating the low order element.}
		\label{fig:MeshProlRest}
	\end{figure}

The prolongation process involves utilizing a linear operator $\Prol$ represented as
\begin{equation}
	\label{eqn:full_prol}
	\Prol =\begin{pmatrix}
		\Prol_\ubm  & 0& 0\\
		0& \Prol_\ybm & 0\\
		0&0 & \Prol_\ubm
	\end{pmatrix}
\end{equation}
 which transfers a solution from the coarse level $\tilde \sbm^{(0)}$ to the fine level $\tilde \sbm $ via $\tilde\sbm=\Pbm \tilde \sbm^{(0)}$. Here $\Pbm_\ubm$ represents the prolongation operator for both the DG coefficients $\ubm$ and the Lagrange multiplier $\lambdabold$. Detailed information about its construction can be found in \cite{fidkowski_p-multigrid_2005}. Additionally, the prolongation $\Pbm_\ybm$ for the mesh $\ybm$ involves inserting high-order nodes into the linear elements, as illustrated in Figure \ref{fig:MeshProlRest}. 
Similarly, a linear restriction operator $\Rest$ is applied, defined as
\begin{equation}
	\label{eqn:full_rest}
	\Rest =\begin{pmatrix}
	\Prol_\ubm^T  & 0& 0\\
	0& \Rest_\ybm & 0\\
	0&0 & \Prol_\ubm^T
\end{pmatrix}.
\end{equation}
This operator projects a fine level solution $\tilde \sbm $ to the coarse level $\tilde \sbm^{(0)}$ via $\tilde\sbm^{(0)}=\Qbm \tilde \sbm$. Here $\Pbm_\ubm^T$ is used for the restriction of both the DG coefficients $\ubm$ and the Lagrange multiplier $\lambdabold$. Furthermore, for the mesh variables $\ybm$, a distinct restriction operator $\Qbm_\ybm$ is employed, which effectively functions as a selection operator, eliminating all high-order nodes ($q>1$). For a single element, this process is depicted in Figure \ref{fig:MeshProlRest}.
\begin{remark}
	\label{remark:mgrd}
	An alternate $p$-multigrid strategy does not restrict/prolongate the mesh nodes $\ybm$, i.e., $\Qbm_\ybm=\Pbm_\ybm=\Ibm$. In a direct comparison with the presented approach, without the restriction/prolongation we obtained slightly better results in terms iterative solver iterations. However, when choosing $\Qbm_\ybm=\Pbm_\ybm=\Ibm$, the coarse level matrix $\Abm^{(0)}$ is significantly bigger in size compared to the presented approach because $\Bbm_{\ybm\ybm}$ is left untouched in the restriction process, which significantly increases the cost per iteration of the multigrid-based preconditioners. Because of the substantial additional cost per iteration and marginal reduction in overall iterations, this approach is not competitive.
\end{remark}

\begin{algorithm}
	\caption{Two level $p$-multigrid}
	\label{alg:mulgrd}
	\begin{algorithmic}
		\REQUIRE KKT matrix $\KKTMat$, right-hand-side $\KKTb$, precomputed prolongation operator $\Prol$, restriction operator $\Rest$, smoother $\PrecV{}$, and coarse matrix $\KKTMat^{(0)}=\Rest \KKTMat \Prol$
		\ENSURE approximate solution $\tilde{\sbm}$ to $\KKTMat \sbm =\KKTb$
		\STATE \textbf{Restrict right-hand side}: $\KKTb^{(0)}=\Rest \KKTb$
		\STATE \textbf{Solve coarse problem}: $\KKTMat^{(0)}\tilde{\sbm}^{(0)} =\KKTb^{(0)}$ 
		\STATE \textbf{Prolongate solution to fine level}: $\tilde{\sbm}=\Prol\tilde{\sbm}^{(0)}$
		\STATE \textbf{Apply smoother}: $\tilde{\sbm} = \tilde{\sbm} + \PrecV{}^{-1}(\KKTb -\KKTMat \tilde{\sbm})$
	\end{algorithmic}
\end{algorithm}

The entire algorithm (Algorithm \ref{alg:mulgrd}) is described as follows: Given the coarse matrix $\KKTMat^{(0)} := \Rest \KKTMat \Prol$, written as
\begin{equation}
	\label{eqn:coarseA}
	\Abm^{(0)}:= \Rest \Abm \Prol =\begin{pmatrix}
		\Prol_\ubm^T \Bbm_{\ubm \ubm} \Prol_\ubm  & 	\Prol_\ubm^T\Bbm_{\ubm\ybm}\Prol_\ybm& 	\Prol_\ubm^T \Ju^T \Prol_\ubm \\
		\Rest_\ybm \Bbm_{\ubm \ybm}^T \Prol_\ubm &\Rest_\ybm \Bbm_{\ybm \ybm} \Prol_\ybm & \Rest_\ybm \Jx^T \Prol_\ubm\\
		\Prol_\ubm^T \Ju \Prol_\ubm&\Prol_\ubm^T \Jx \Prol_\ybm & 0
	\end{pmatrix},
\end{equation} 
the right-hand-side is restricted to the coarse level: $\KKTb^{(0)} = \Rest \KKTb$. Then, the coarse problem is solved: $\KKTMat^{(0)}\tilde{\sbm}^{(0)} = \KKTb^{(0)}$ using a direct sparse solve. Subsequently, the solution is prolonged back to the fine level as $\tilde{\sbm} = \Prol\tilde{\sbm}^{(0)}$ and an iterative smoothing process is applied as $\tilde{\sbm} = \tilde{\sbm} + \PrecV{}^{-1}(\KKTb -\KKTMat \tilde{\sbm})$. 

We close this section by summarizing all eight preconditioners that will be studied in Section~\ref{sec:numexp:results} in Table \ref{tab:all_preconditioners}. As we did not observe a significant benefit in preliminary studies, we do not combine $p$-multigrid with the preconditioners where $\tilde\Bbm_{\ybm \ybm} = \text{ilu}(\Bbm_{\ybm \ybm})$.
\begin{table}[htbp]
	\centering
	\renewcommand{\arraystretch}{1.5} 
	\caption{Summary of all HOIST preconditioners studied.}
	\label{tab:all_preconditioners}
	\begin{tabular}{|c|c|c|c|}
		\hline
		\textbf{Preconditioner} &$\tilde\Bbm_{\ybm \ybm}\approx \Bbm_{\ybm \ybm} $ & $\tilde \Ju \approx \Ju$ & $p$-\text{multigrid}\\ \hline
		$\PrecV{0}$ & $\Bbm_{\ybm \ybm}$ & $\Ju$ & \text{no}\\ \hline
		$\ABJ$ & $\text{diag}(\Bbm_{\ybm \ybm})$ & $\tilde{\Jbm}_{\text{BJ}}$ & \text{no}\\ \hline
		$\ABILU$ & $\text{diag}(\Bbm_{\ybm \ybm})$ & $\tilde{\Jbm}_{\text{BILU}}$ & \text{no}\\ \hline
		$\ABJilu$  & $\text{ilu}(\Bbm_{\ybm \ybm})$ & $\tilde{\Jbm}_{\text{BJ}}$ & \text{no}\\ \hline
		$\ABILUilu$ & $\text{ilu}(\Bbm_{\ybm \ybm})$ & $\tilde{\Jbm}_{\text{BILU}}$ & \text{no}\\ \hline
		$\PrecV{0p0}$ & $\Bbm_{\ybm \ybm}$ & $\Ju$ & \text{yes}\\ \hline
		$\ABJp$ & $\text{diag}(\Bbm_{\ybm \ybm})$ & $\tilde{\Jbm}_{\text{BJ}}$ & \text{yes}\\ \hline
		$\ABILUp$ & $\text{diag}(\Bbm_{\ybm \ybm})$ & $\tilde{\Jbm}_{\text{BILU}}$ & \text{yes}\\ \hline
	\end{tabular}
\end{table}
	\section{Numerical experiments}
	\label{sec:numexp}
In this section, we present a series of numerical experiments designed to evaluate the performance of the introduced preconditioners. First, we define the metrics employed to measure their effectiveness (Section~\ref{sec:numexp:metrics}) and describe two shock-dominated flow benchmarks (Euler equations) that will be used to study the preconditioners (Section~\ref{sec:numexp:cases}). Finally, we present and analyze the results from numerical experiments, focusing on various key HOIST parameters (Section~\ref{sec:numexp:results}). We solely consider the generalized minimum residual (GMRES) Krylov solver in all studies because the preconditioned system does not have special structure that would allow us to use a more specialized solver. 

In some applications, GMRES is used with a restart technique, to mitigate the linear growth (with respect to iteration numbers) of memory needed. The idea is to limit the number of basis vectors stored by restarting the GMRES algorithm after a fixed number of iterations $m'$. However, restarting can sometimes lead to slower convergence because valuable information from previous iterations is discarded and the choice of $m'$ can significantly impact performance. In order to not introduce another parameter ($m'$) to study, we choose to use GMRES without restart for the present work and leave the investigation for future research.
\subsection{Description of metrics}
	\label{sec:numexp:metrics}
We assess the performance of the preconditioners based on the number of GMRES iterations required to achieve a convergence criterion. In practical applications, this involves monitoring the relative residual norm of the preconditioned system and stopping at the first iteration where
	\begin{equation}
		\label{eqn:gmres-convcrit-res}
		 \frac{\Vert \tilde \Abm^{-1} \Abm \sbm - \tilde \Abm^{-1}\bbm \Vert}{\Vert\tilde \Abm^{-1}\bbm \Vert} < \text{tol},
	\end{equation}
	where $\text{tol}>0$ is a specified tolerance. It is important to note that this convergence criterion is preconditioner-dependent. To ensure a fair comparison, we opt for a convergence criterion based on the exact solution $\sbm_{\text{ex}}$ satisfying $\Abm \sbm_{\text{ex}}=\bbm$:
	\begin{equation}
	\label{eqn:gmres-convcrit-ex}
		\frac{\Vert \sbm_{\text{ex}}-\sbm \Vert}{\Vert \sbm_{\text{ex}} \Vert}<\text{tol},
	\end{equation}
	with $\text{tol}=10^{-3}$. We also set the maximal GMRES iterations to be $1000$.
	
	The parameter space influencing the effectiveness of preconditioners for the HOIST method is vast and multifaceted. It encompasses choices related to the equations, specific problem formulations, the number of elements $\vert \Ecal_h \vert$ utilized, the polynomial degrees $p$ and $q$, the state $\zbm_k$ around which the system is linearized, and finally, the selection of $\gamma$ and $\kappa$, which significantly affect the condition number of the system. Studying all these dimensions collectively is infeasible. Consequently, we will conduct separate investigations to gauge the relative impact of each of these parameters.
	\begin{table}
		\centering
		\caption{Legend for plots comparing (\# GMRES Iterations) for each preconditioner}
		\label{tab:legend}
		\begin{tabular}{c|c|c|c|c|c|c|c}
			$\PrecV{0}$ &$\PrecV{0$p0$}$ &$\ABILU$ &$\PrecV{BILU$p0$}$&$\ABJ$ &$\PrecV{BJ$p0$}$&$\ABILUilu$&$\ABJilu$ \\\hline
			(\ref{line:A0}) & (\ref{line:A0p0}) & (\ref{line:BILU}) &
			(\ref{line:BILUp0}) & (\ref{line:BJ}) & (\ref{line:BJp0}) &
			(\ref{line:BILU_ILU}) & (\ref{line:BJ_ILU})
		\end{tabular}
	\end{table}

	\subsection{Description of examined cases}
	\label{sec:numexp:cases}
	In this work, we focus exclusively on experiments related to the steady, inviscid two-dimensional Euler equations (Section \ref{sec:numexp:euler}). Specifically, we consider two problems with unique solution features: supersonic flow around a cylinder (Section \ref{sec:numexp:euler:cylinder}) and supersonic flow over a diamond-shaped obstacle in a tunnel (Section \ref{sec:numexp:euler:diamond}).
	
	\subsubsection{Inviscid Euler equations}
	\label{sec:numexp:euler}
	Compressible, inviscid flow through the domain $\Omega \in \Rbb^d$ is modeled using the Euler equations of gasdynamics
			\begin{equation} \label{eqn:euler}
		\begin{split}
			\pder{}{x_j}\left(\rho(x) v_j(x)\right) &= 0 \\
			\pder{}{x_j}\left(\rho(x) v_i(x)v_j(x)+P(x)\delta_{ij}\right) &= 0 \\
			\pder{}{x_j}\left(\left[\rho(x)E(x)+P(x)\right]v_j(x)\right) &= 0
		\end{split}
	\end{equation}
for all $x\in\Omega$ and where $i=1,\dots,d$ and summation is implied over
the repeated index $j=1,\dots,d$. The density of the fluid $\func{\rho}{\Omega}{\Rbb_{>0}}$, the velocity of the fluid $\func{v_i}{\Omega}{\Rbb}$ in $x_i$
direction for $i=1,\dots,d$, the total energy of the fluid $\func{E}{\Omega}{\Rbb_{>0}}$, the pressure of the fluid $\func{P}{\Omega}{\Rbb_{>0}}$ are implicitly defined as the solution of (\ref{eqn:euler}). The equations are closed by introducing an equation of state, the ideal gas law in this work,
		\begin{equation}
	P = (\gamma-1)\left(\rho E - \frac{\rho v_i v_i}{2}\right),
\end{equation}
where $\gamma\in\Rbb_{>0}$ is the ratio of specific heats, typically  $\gamma=1.4$ for air at standard conditions.

	\subsubsection{Supersonic flow over a cylinder}
	\label{sec:numexp:euler:cylinder}
		In our first problem, we explore a supersonic flow (Mach 2) over a cylinder to demonstrate the preconditioners performance for problems with curved shocks. 		
	 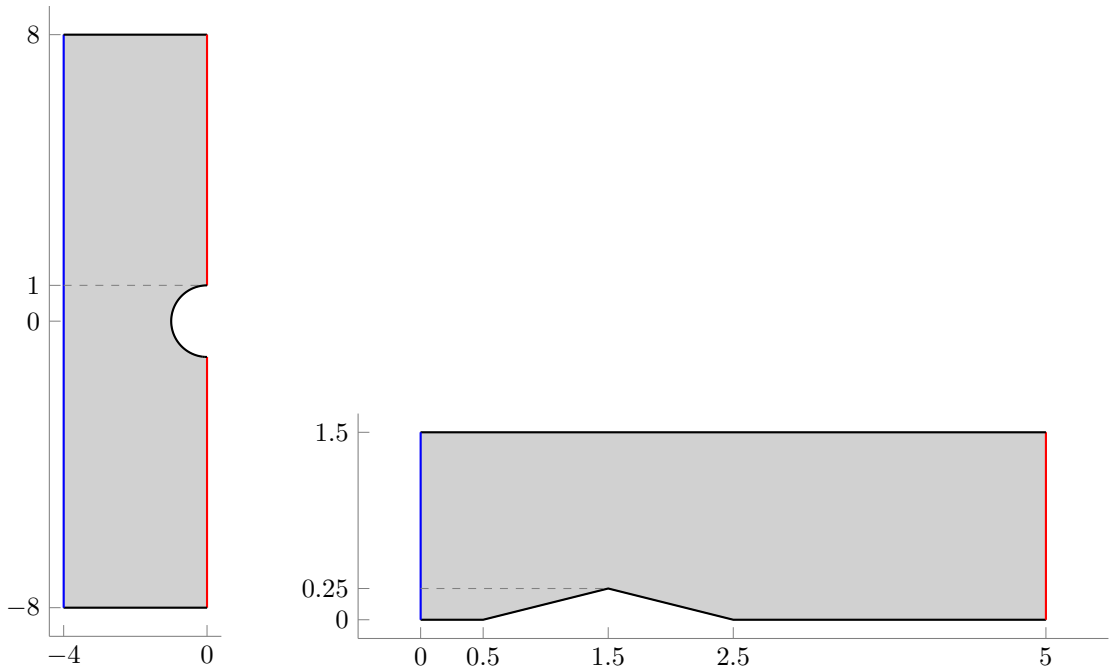
\begin{figure}[!htbp]
		\centering
		\ifbool{fastcompile}{}{
			\begin{tikzpicture}
\begin{axis}[
axis equal image,
axis line style={gray},
axis x line*=bottom,
axis y line*=left,
width=0.7\textwidth,
xtick={-4.0,0.0},
ytick={-8.0,0.0,1.0, 8.0},
ymax=8.8,
xmax=0.4000000000000001,
xmin=-4.4,
ymin=-8.8]
\addplot [opacity=0.6, fill=black!30!white, forget plot]
coordinates {
( 0.00000000e+00, -1.00000000e+00)
( 0.00000000e+00, -8.00000000e+00)
(-4.00000000e+00, -8.00000000e+00)
(-4.00000000e+00,  8.00000000e+00)
( 0.00000000e+00,  8.00000000e+00)
( 0.00000000e+00,  1.00000000e+00)
( 6.12323400e-17,  1.00000000e+00)
(-6.27905195e-02,  9.98026728e-01)
(-1.25333234e-01,  9.92114701e-01)
(-1.87381315e-01,  9.82287251e-01)
(-2.48689887e-01,  9.68583161e-01)
(-3.09016994e-01,  9.51056516e-01)
(-3.68124553e-01,  9.29776486e-01)
(-4.25779292e-01,  9.04827052e-01)
(-4.81753674e-01,  8.76306680e-01)
(-5.35826795e-01,  8.44327926e-01)
(-5.87785252e-01,  8.09016994e-01)
(-6.37423990e-01,  7.70513243e-01)
(-6.84547106e-01,  7.28968627e-01)
(-7.28968627e-01,  6.84547106e-01)
(-7.70513243e-01,  6.37423990e-01)
(-8.09016994e-01,  5.87785252e-01)
(-8.44327926e-01,  5.35826795e-01)
(-8.76306680e-01,  4.81753674e-01)
(-9.04827052e-01,  4.25779292e-01)
(-9.29776486e-01,  3.68124553e-01)
(-9.51056516e-01,  3.09016994e-01)
(-9.68583161e-01,  2.48689887e-01)
(-9.82287251e-01,  1.87381315e-01)
(-9.92114701e-01,  1.25333234e-01)
(-9.98026728e-01,  6.27905195e-02)
(-1.00000000e+00,  1.22464680e-16)
(-9.98026728e-01, -6.27905195e-02)
(-9.92114701e-01, -1.25333234e-01)
(-9.82287251e-01, -1.87381315e-01)
(-9.68583161e-01, -2.48689887e-01)
(-9.51056516e-01, -3.09016994e-01)
(-9.29776486e-01, -3.68124553e-01)
(-9.04827052e-01, -4.25779292e-01)
(-8.76306680e-01, -4.81753674e-01)
(-8.44327926e-01, -5.35826795e-01)
(-8.09016994e-01, -5.87785252e-01)
(-7.70513243e-01, -6.37423990e-01)
(-7.28968627e-01, -6.84547106e-01)
(-6.84547106e-01, -7.28968627e-01)
(-6.37423990e-01, -7.70513243e-01)
(-5.87785252e-01, -8.09016994e-01)
(-5.35826795e-01, -8.44327926e-01)
(-4.81753674e-01, -8.76306680e-01)
(-4.25779292e-01, -9.04827052e-01)
(-3.68124553e-01, -9.29776486e-01)
(-3.09016994e-01, -9.51056516e-01)
(-2.48689887e-01, -9.68583161e-01)
(-1.87381315e-01, -9.82287251e-01)
(-1.25333234e-01, -9.92114701e-01)
(-6.27905195e-02, -9.98026728e-01)};

\addplot [thick, color=red]
coordinates {
( 0.00000000e+00, -1.00000000e+00)
( 0.00000000e+00, -8.00000000e+00)};\label{line:cyl:outlet}

\addplot [thick, color=black]
coordinates {
( 0.00000000e+00, -8.00000000e+00)
(-4.00000000e+00, -8.00000000e+00)};\label{line:cyl:wall}

\addplot [thick, color=blue]
coordinates {
(-4.00000000e+00, -8.00000000e+00)
(-4.00000000e+00,  8.00000000e+00)};\label{line:cyl:inlet}

\addplot [thick, color=black]
coordinates {
(-4.00000000e+00,  8.00000000e+00)
( 0.00000000e+00,  8.00000000e+00)};\label{line:cyl:wall}

\addplot [thick, color=red]
coordinates {
( 0.00000000e+00,  8.00000000e+00)
( 0.00000000e+00,  1.00000000e+00)};\label{line:cyl:outlet}

\addplot [dashed, color=gray, forget plot]
coordinates {
(-4.00000000e+00,  1.00000000e+00)
( 0.00000000e+00,  1.00000000e+00)};

\addplot [thick, color=black]
coordinates {
( 0.00000000e+00,  1.00000000e+00)
( 6.12323400e-17,  1.00000000e+00)
(-6.27905195e-02,  9.98026728e-01)
(-1.25333234e-01,  9.92114701e-01)
(-1.87381315e-01,  9.82287251e-01)
(-2.48689887e-01,  9.68583161e-01)
(-3.09016994e-01,  9.51056516e-01)
(-3.68124553e-01,  9.29776486e-01)
(-4.25779292e-01,  9.04827052e-01)
(-4.81753674e-01,  8.76306680e-01)
(-5.35826795e-01,  8.44327926e-01)
(-5.87785252e-01,  8.09016994e-01)
(-6.37423990e-01,  7.70513243e-01)
(-6.84547106e-01,  7.28968627e-01)
(-7.28968627e-01,  6.84547106e-01)
(-7.70513243e-01,  6.37423990e-01)
(-8.09016994e-01,  5.87785252e-01)
(-8.44327926e-01,  5.35826795e-01)
(-8.76306680e-01,  4.81753674e-01)
(-9.04827052e-01,  4.25779292e-01)
(-9.29776486e-01,  3.68124553e-01)
(-9.51056516e-01,  3.09016994e-01)
(-9.68583161e-01,  2.48689887e-01)
(-9.82287251e-01,  1.87381315e-01)
(-9.92114701e-01,  1.25333234e-01)
(-9.98026728e-01,  6.27905195e-02)
(-1.00000000e+00,  1.22464680e-16)
(-9.98026728e-01, -6.27905195e-02)
(-9.92114701e-01, -1.25333234e-01)
(-9.82287251e-01, -1.87381315e-01)
(-9.68583161e-01, -2.48689887e-01)
(-9.51056516e-01, -3.09016994e-01)
(-9.29776486e-01, -3.68124553e-01)
(-9.04827052e-01, -4.25779292e-01)
(-8.76306680e-01, -4.81753674e-01)
(-8.44327926e-01, -5.35826795e-01)
(-8.09016994e-01, -5.87785252e-01)
(-7.70513243e-01, -6.37423990e-01)
(-7.28968627e-01, -6.84547106e-01)
(-6.84547106e-01, -7.28968627e-01)
(-6.37423990e-01, -7.70513243e-01)
(-5.87785252e-01, -8.09016994e-01)
(-5.35826795e-01, -8.44327926e-01)
(-4.81753674e-01, -8.76306680e-01)
(-4.25779292e-01, -9.04827052e-01)
(-3.68124553e-01, -9.29776486e-01)
(-3.09016994e-01, -9.51056516e-01)
(-2.48689887e-01, -9.68583161e-01)
(-1.87381315e-01, -9.82287251e-01)
(-1.25333234e-01, -9.92114701e-01)};\label{line:cyl:wall}

\addplot [thick, color=black]
coordinates {
(-6.27905195e-02, -9.98026728e-01)
(-1.25333234e-01, -9.92114701e-01)};\label{line:cyl:wall}

\addplot [thick, color=black]
coordinates {
( 0.00000000e+00, -1.00000000e+00)
(-6.27905195e-02, -9.98026728e-01)};\label{line:cyl:wall}

\end{axis}
\end{tikzpicture} \qquad
			\begin{tikzpicture}
\begin{axis}[
axis equal image,
axis line style={gray},
axis x line*=bottom,
axis y line*=left,
width=0.7\textwidth,
xtick={0.0, 0.5, 1.5, 2.5, 5.0},
ytick={0.0, 0.25, 1.5},
ymax=1.65,
xmax=5.5,
xmin=-0.5,
ymin=-0.15]
\addplot [opacity=0.6, fill=black!30!white, opacity=0.6, forget plot]
coordinates {
( 0.00000000e+00,  0.00000000e+00)
( 5.00000000e-01,  0.00000000e+00)
( 1.50000000e+00,  2.50000000e-01)
( 2.50000000e+00,  0.00000000e+00)
( 5.00000000e+00,  0.00000000e+00)
( 5.00000000e+00,  1.50000000e+00)
( 0.00000000e+00,  1.50000000e+00)
( 0.00000000e+00,  0.00000000e+00)};

\addplot [thick, color=black]
coordinates {
( 0.00000000e+00,  0.00000000e+00)
( 5.00000000e-01,  0.00000000e+00)
( 1.50000000e+00,  2.50000000e-01)
( 2.50000000e+00,  0.00000000e+00)
( 5.00000000e+00,  0.00000000e+00)};\label{line:diamond:wall}

\addplot [thick, color=red]
coordinates {
( 5.00000000e+00,  0.00000000e+00)
( 5.00000000e+00,  1.50000000e+00)};\label{line:diamond:outlet}

\addplot [thick, color=black]
coordinates {
( 5.00000000e+00,  1.50000000e+00)
( 0.00000000e+00,  1.50000000e+00)};\label{line:diamond:wall}

\addplot [thick, color=blue]
coordinates {
( 0.00000000e+00,  1.50000000e+00)
( 0.00000000e+00,  0.00000000e+00)};\label{line:diamond:inlet}

\addplot [dashed, color=gray, forget plot]
coordinates {
( 0.00000000e+00,  2.50000000e-01)
( 1.50000000e+00,  2.50000000e-01)};

\end{axis}
\end{tikzpicture}
		}
		\caption{Geometry and boundary conditions for the \ecyl~(\textit{left}) and \texttt{diamond}~(\textit{right}) test cases.
			Boundary conditions: slip walls (\ref{line:cyl:wall}), Mach 2 supersonic inflow
			(\ref{line:cyl:inlet}), and supersonic outflow (\ref{line:cyl:outlet}).}
		\label{fig:euler:cyl:geom}
	\end{figure}
 
\ifbool{fastcompile}{}{
	\begin{figure}[!htbp]
		\centering
		\begin{tikzpicture}
			\begin{groupplot}[
				group style={
					group size=6 by 1,
					horizontal sep=0.4cm,
					vertical sep=0.4cm
				},
				width=0.7\textwidth,
				axis equal image,
				xticklabels={,,},
				yticklabels={,,},
				xmin=-4, xmax=0,
				ymin=-8, ymax=8
				]
				\nextgroupplot
				\addplot graphics [xmin=-4, xmax=0, ymin=-8, ymax=8] {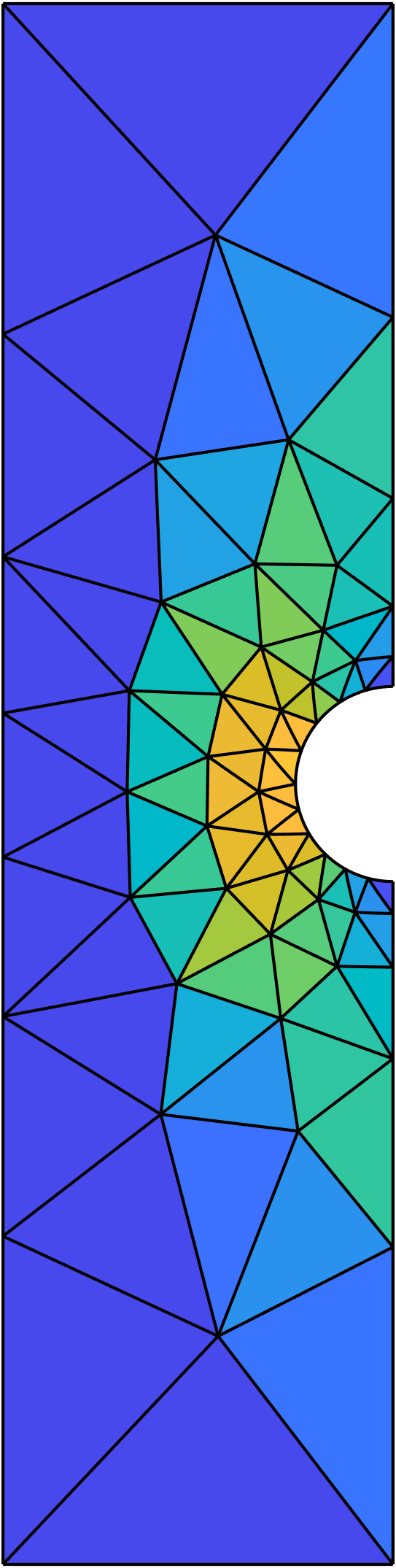};
				\nextgroupplot
				\addplot graphics [xmin=-4, xmax=0, ymin=-8, ymax=8] {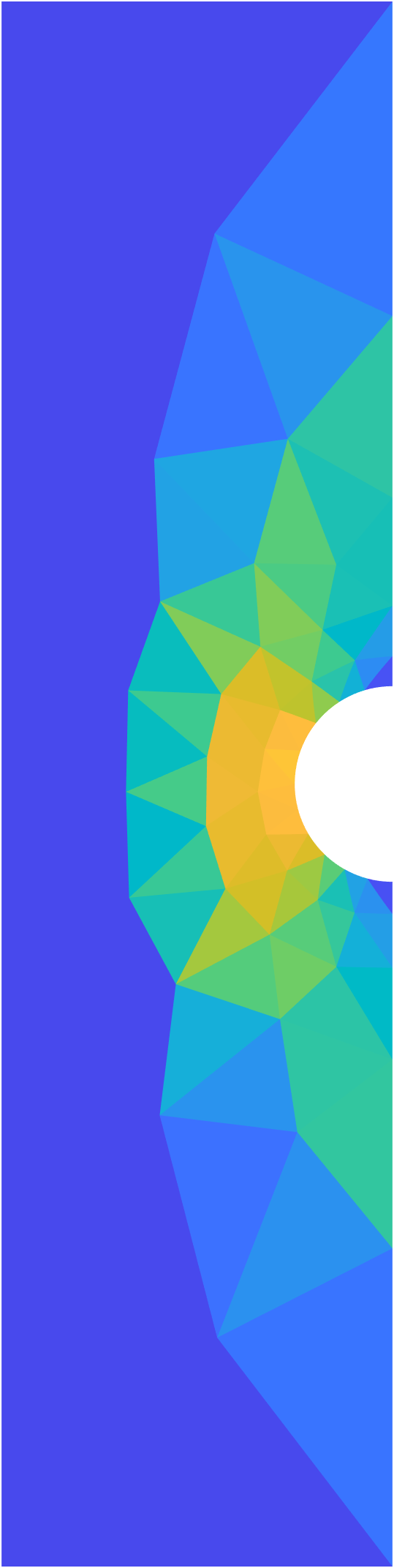};
				\nextgroupplot
				\addplot graphics [xmin=-4, xmax=0, ymin=-8, ymax=8] {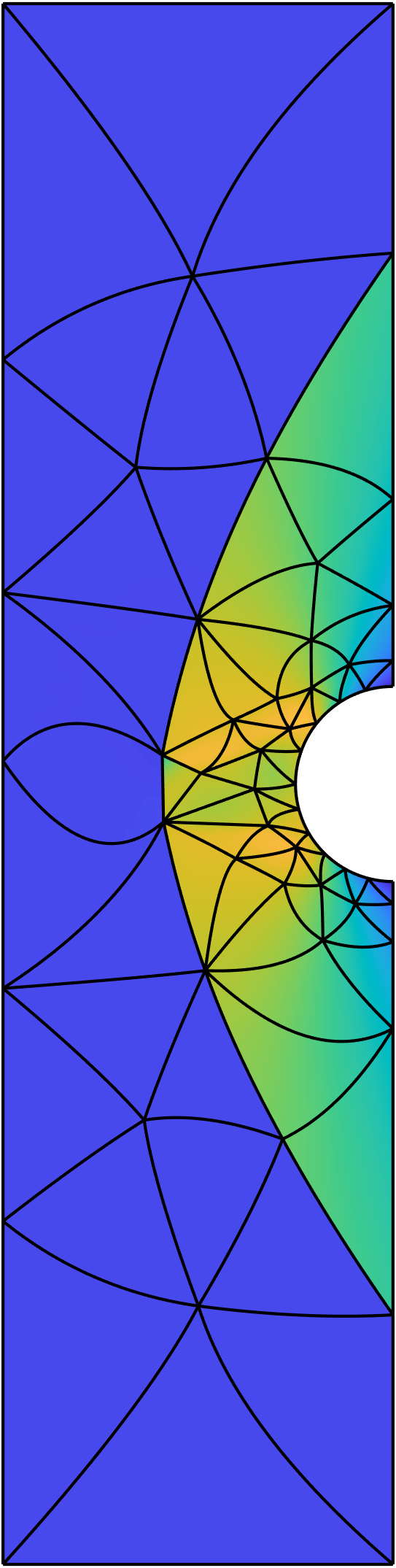};
				\nextgroupplot
				\addplot graphics [xmin=-4, xmax=0, ymin=-8, ymax=8] {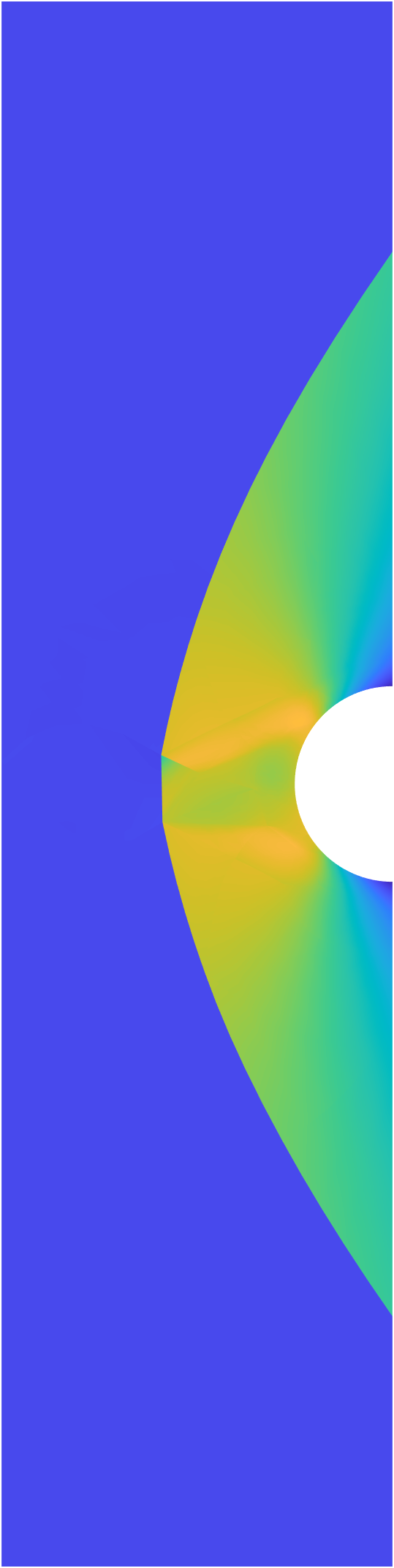};
				\nextgroupplot
				\addplot graphics [xmin=-4, xmax=0, ymin=-8, ymax=8] {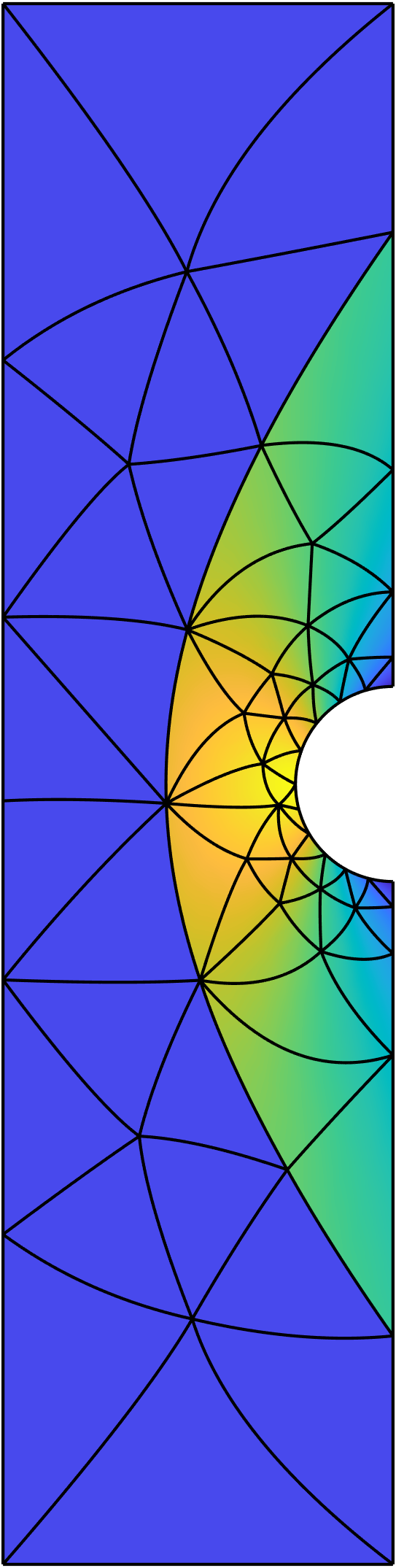};
				\nextgroupplot
				\addplot graphics [xmin=-4, xmax=0, ymin=-8, ymax=8] {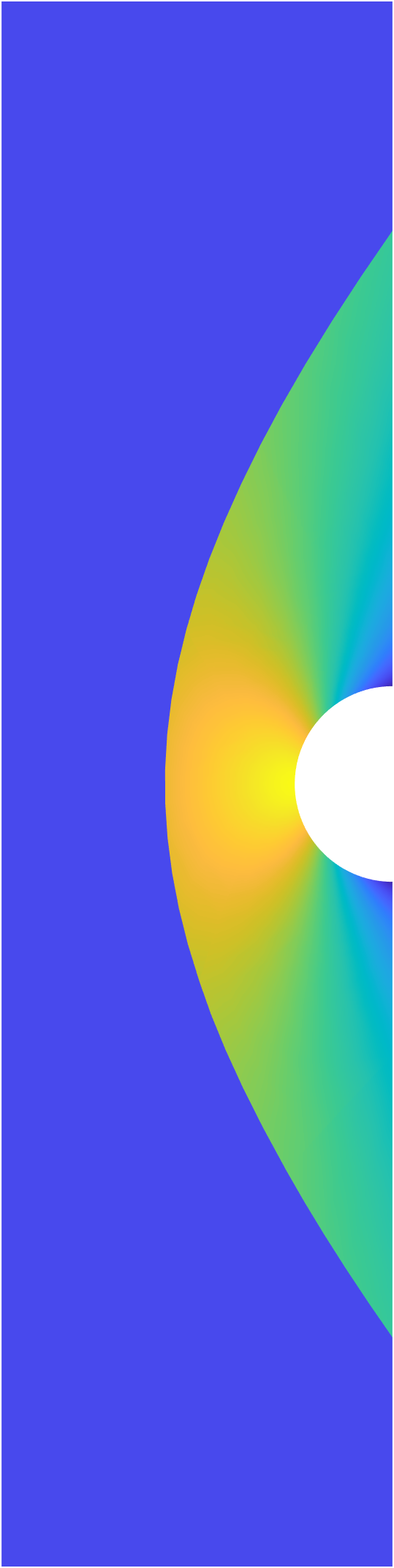};
			\end{groupplot}
		\end{tikzpicture}
		\colorbarMatlabParula{1.0}{1}{2}{3}{4.4}
		\caption{Selected $p=q=2$ HOIST iterations $k\in \{1,50,100\}$ (\textit{left-to-right}) for the \ecyl~test case (density) shown with and without mesh edges.}
		\label{fig:euler:cyl:sltn}
	\end{figure}
}
	  The domain (Figure \ref{fig:euler:cyl:geom}, problem: \ecyl) is discretized using a coarse unstructured triangular mesh with $90$ elements and 
	 throughout most of our investigations, we utilize a third-order approximation for the flow variables and the geometry ($p=q=2$), allowing the HOIST method to iterate to $k=100$. A first-order finite volume solution is used for initialization of the method and it converges to a mesh that tracks the shock. Figure \ref{fig:euler:cyl:sltn} displays the density for selected iterations $k=1,50,100$ obtained for this configuration. For the upcoming studies, the corresponding linear systems derived from the states $\zbm_1,\zbm_{50},\zbm_{100}$ are utilized to evaluate the performance of the preconditioners.
	  
	 	\begin{figure}[!htbp]
	 		\centering
	 			 \ifbool{fastcompile}{}{
	 		\begin{tikzpicture}
	 			\begin{groupplot}[
	 				group style={
	 					group size=6 by 1,
	 					horizontal sep=0.4cm,
	 					vertical sep=0.4cm
	 				},
	 				width=0.7\textwidth,
	 				axis equal image,
	 				xticklabels={,,},
	 				yticklabels={,,},
	 				xmin=-4, xmax=0,
	 				ymin=-8, ymax=8
	 				]
	 				\nextgroupplot
	 				\addplot graphics [xmin=-4, xmax=0, ymin=-8, ymax=8] {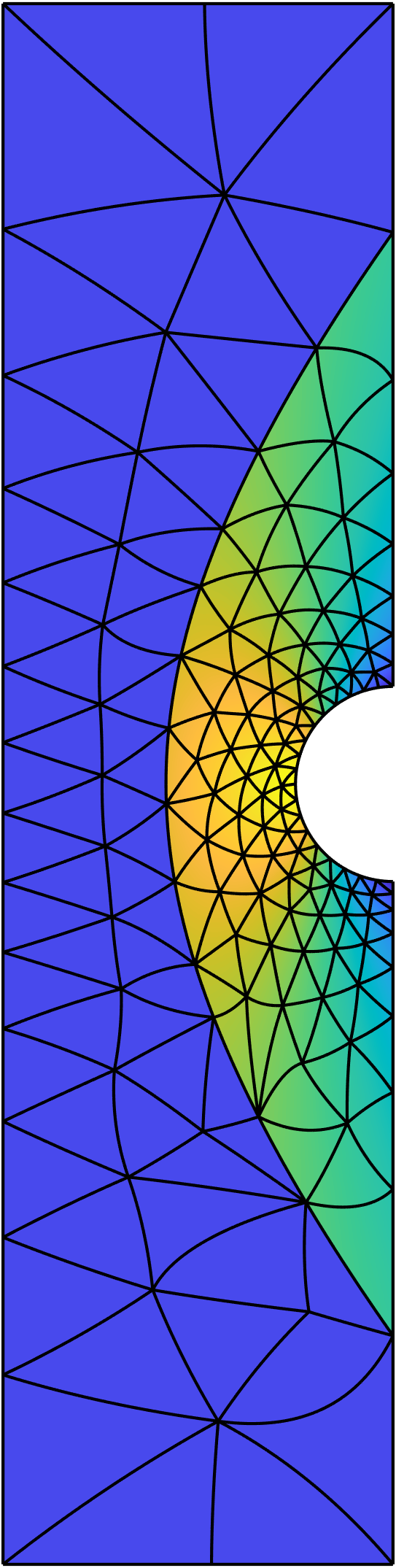};
	 				\nextgroupplot
	 				\addplot graphics [xmin=-4, xmax=0, ymin=-8, ymax=8] {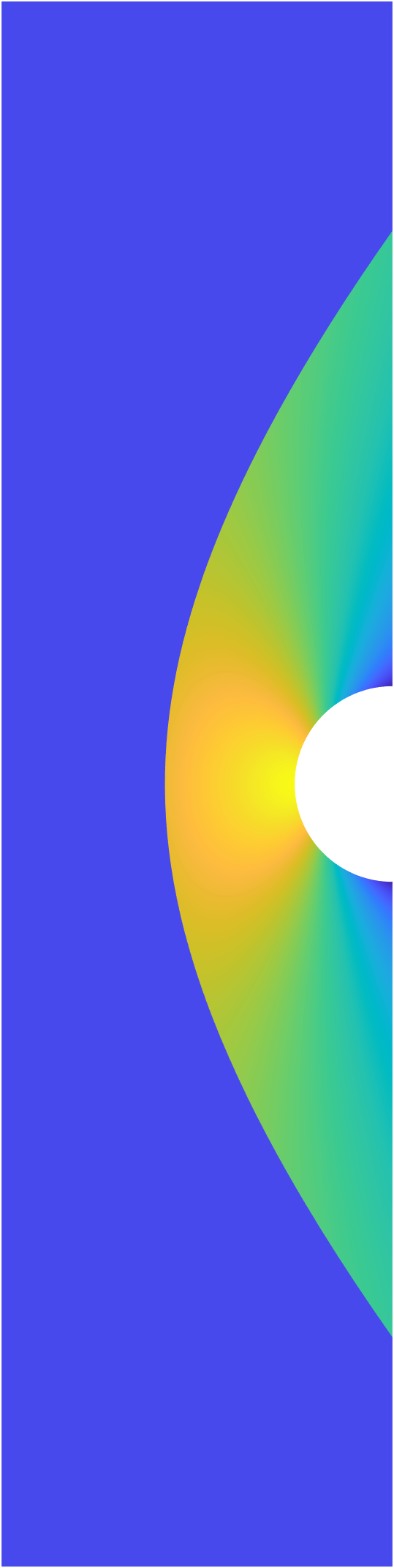};
	 				\nextgroupplot
	 				\addplot graphics [xmin=-4, xmax=0, ymin=-8, ymax=8] {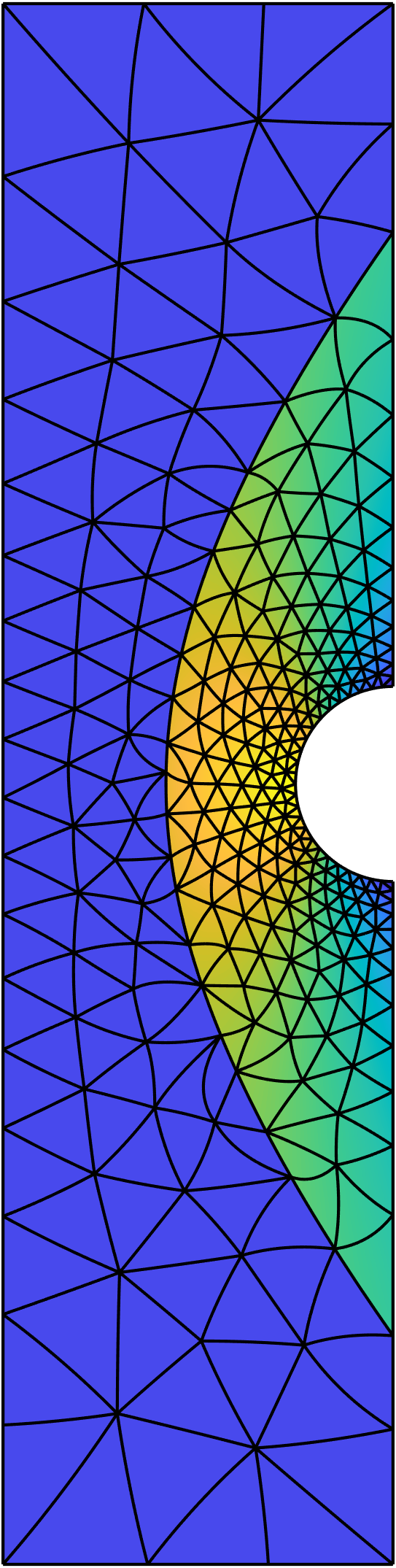};
	 				\nextgroupplot
	 				\addplot graphics [xmin=-4, xmax=0, ymin=-8, ymax=8] {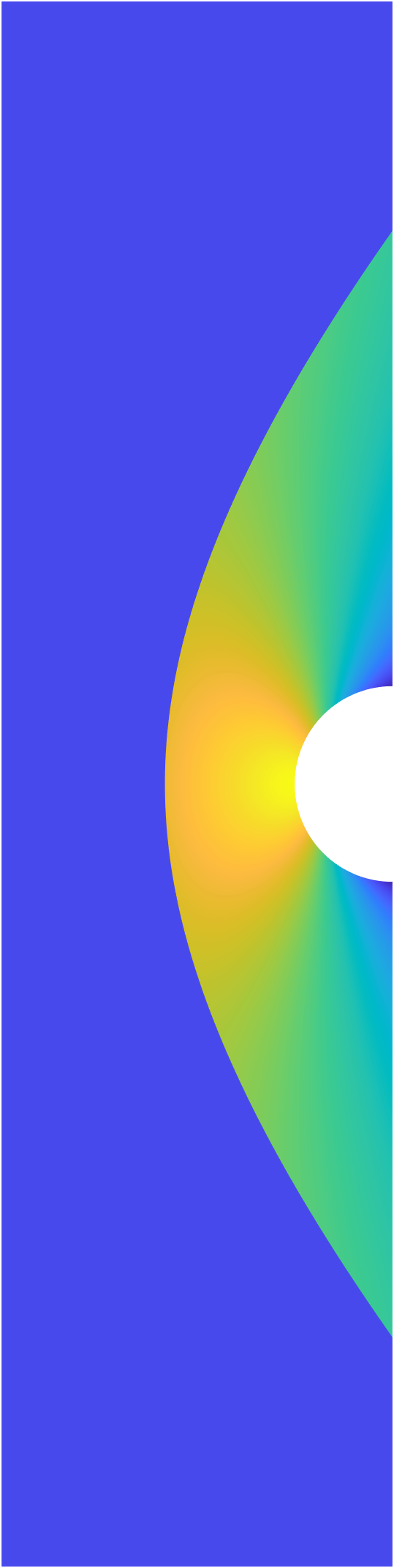};
	 				\nextgroupplot
	 				\addplot graphics [xmin=-4, xmax=0, ymin=-8, ymax=8] {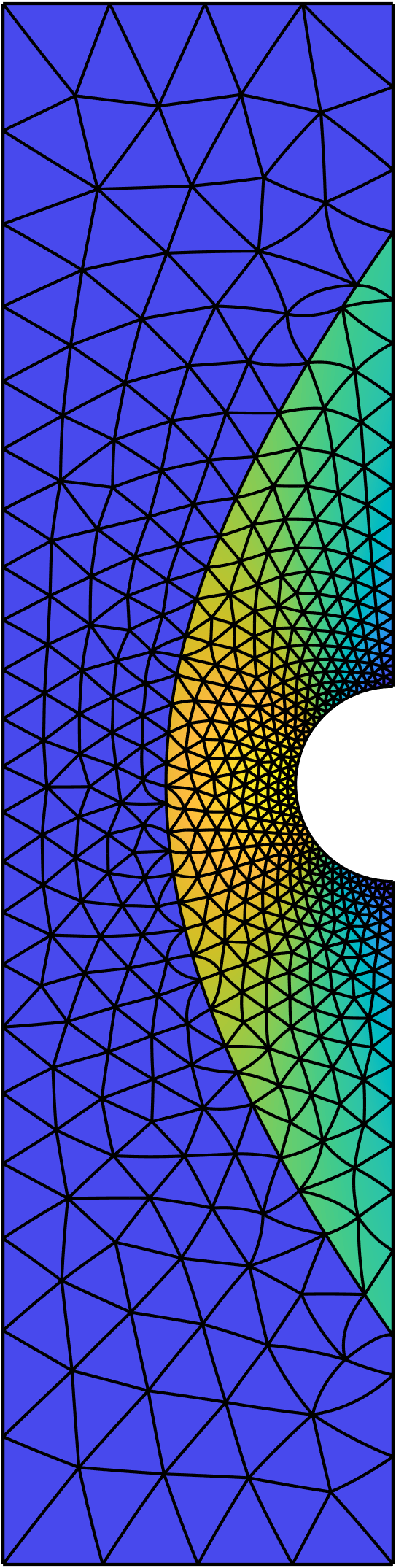};
	 				\nextgroupplot
	 				\addplot graphics [xmin=-4, xmax=0, ymin=-8, ymax=8] {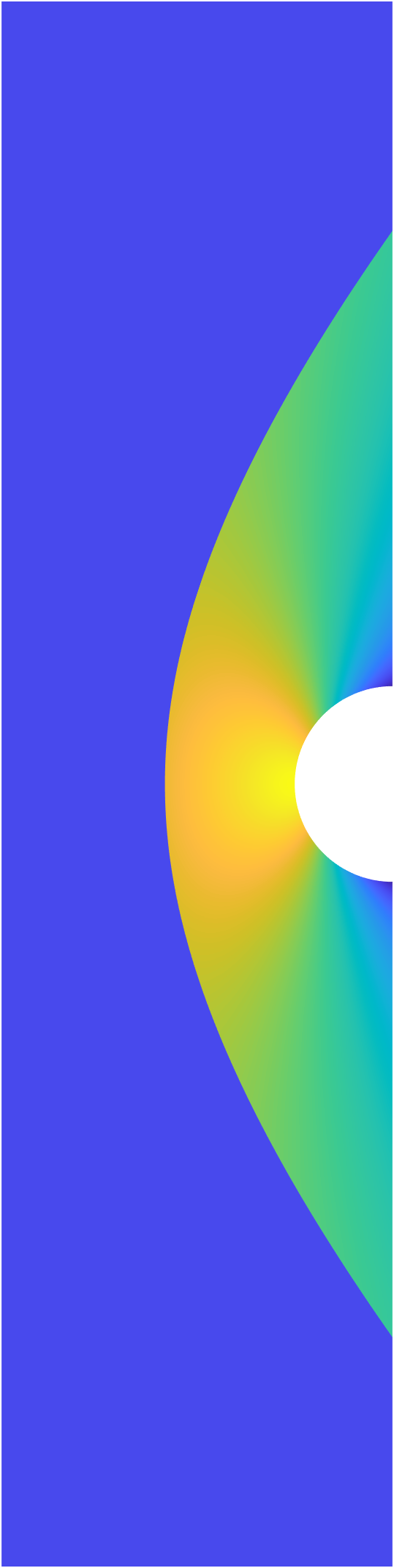};
	 			\end{groupplot}
	 		\end{tikzpicture}
	 		\colorbarMatlabParula{1.0}{1}{2}{3}{4.4}
	 			 }
	 		\caption{Selected $p=q=2$ HOIST iterations $k=100$ and different refinement levels $n_\text{ref}=1,2,3$ (\textit{left-to-right}) for the \ecyl~test case (density) shown with and without mesh edges.}
	 		\label{fig:euler:cyl:refinement}
	 	\end{figure}

\begin{figure}[!htbp]
			\centering
				\ifbool{fastcompile}{}{
			\begin{tikzpicture}
				\begin{groupplot}[
					group style={
						group size=6 by 1,
						horizontal sep=0.4cm,
						vertical sep=0.4cm
					},
					width=0.7\textwidth,
					axis equal image,
					xticklabels={,,},
					yticklabels={,,},
					xmin=-4, xmax=0,
					ymin=-8, ymax=8
					]
					\nextgroupplot
					\addplot graphics [xmin=-4, xmax=0, ymin=-8, ymax=8] {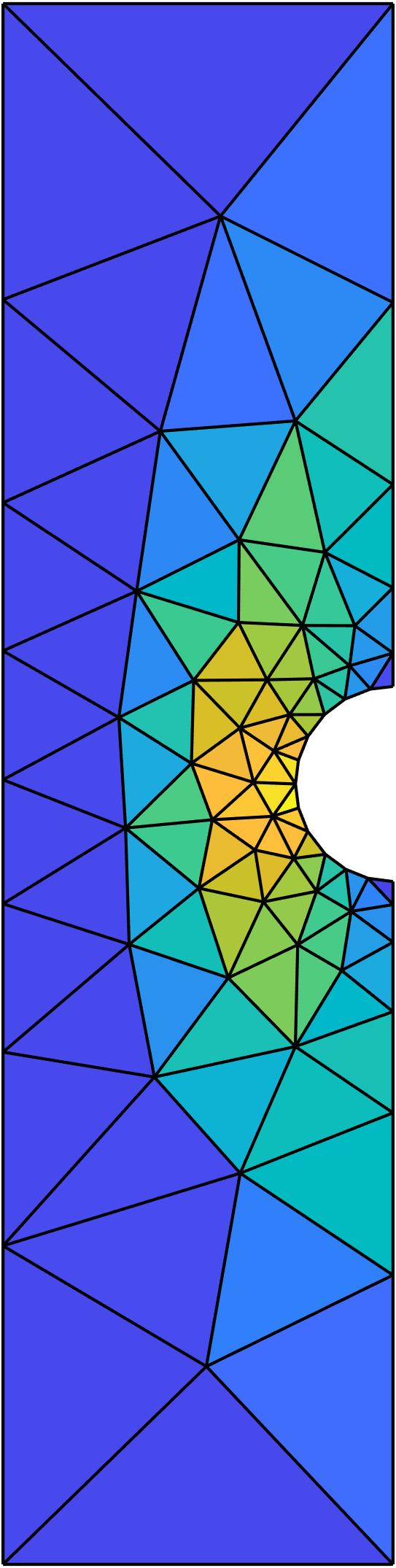};
					\nextgroupplot
					\addplot graphics [xmin=-4, xmax=0, ymin=-8, ymax=8] {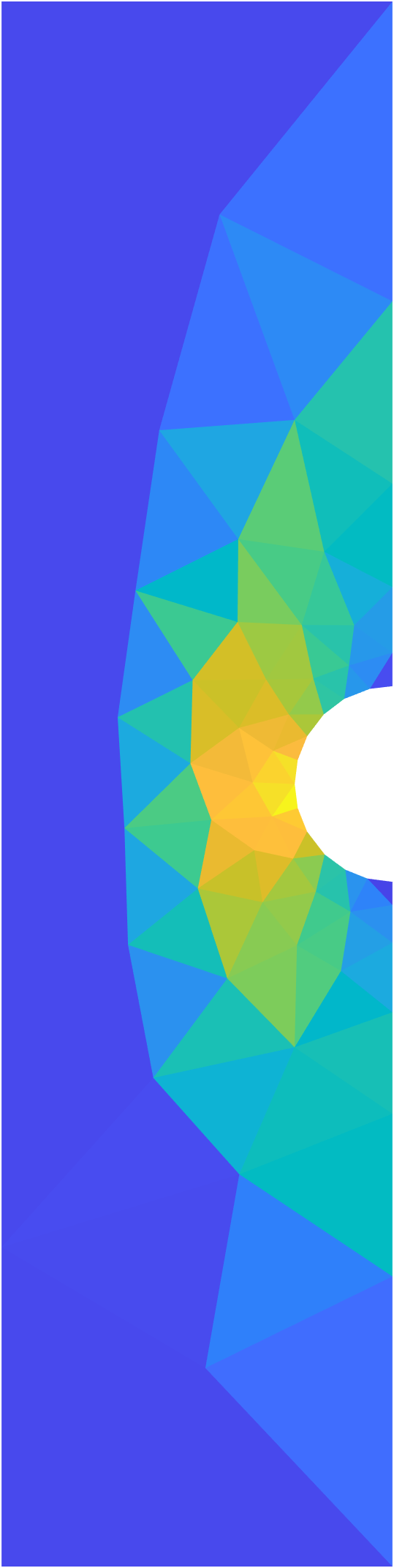};
					\nextgroupplot
					\addplot graphics [xmin=-4, xmax=0, ymin=-8, ymax=8] {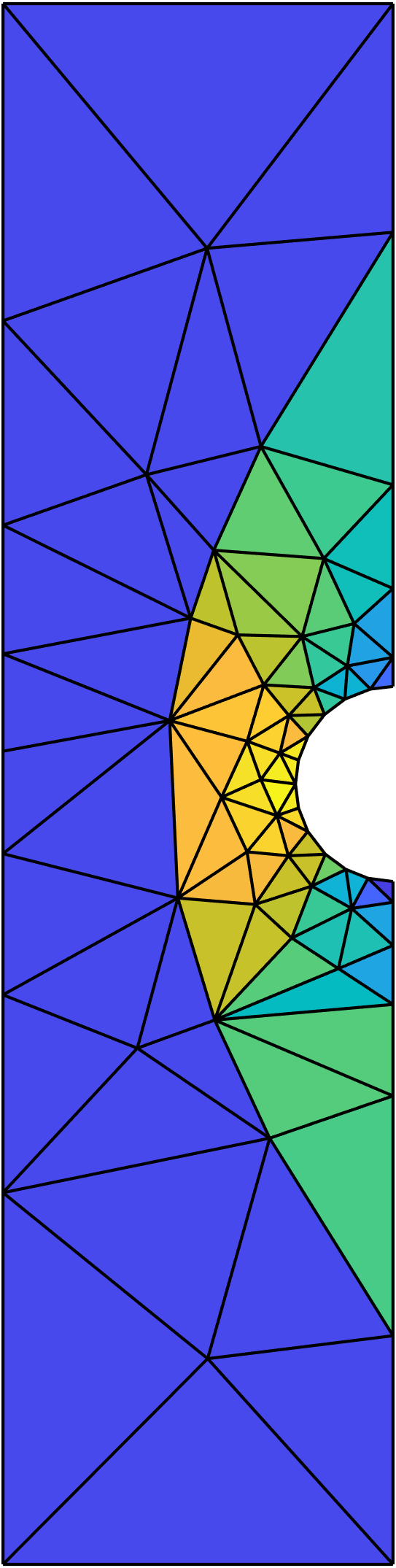};
					\nextgroupplot
					\addplot graphics [xmin=-4, xmax=0, ymin=-8, ymax=8] {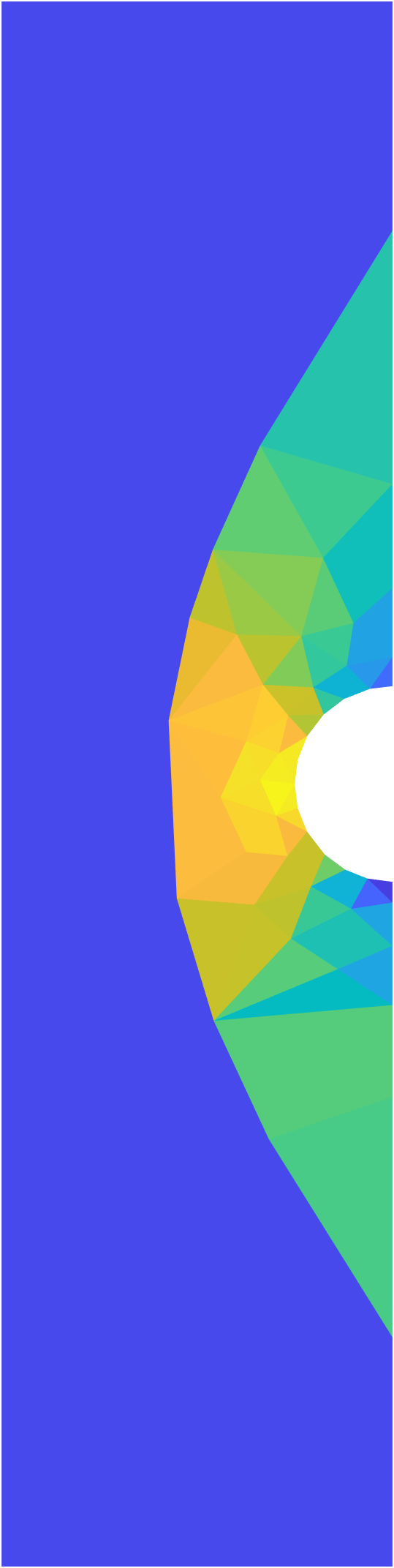};
					\nextgroupplot
					\addplot graphics [xmin=-4, xmax=0, ymin=-8, ymax=8] {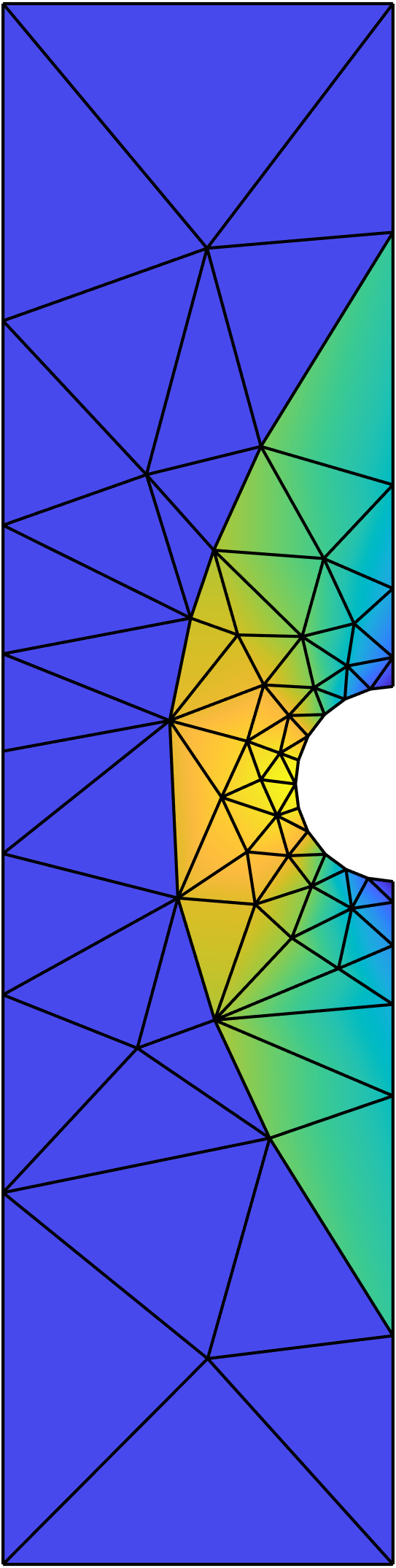};
					\nextgroupplot
					\addplot graphics [xmin=-4, xmax=0, ymin=-8, ymax=8] {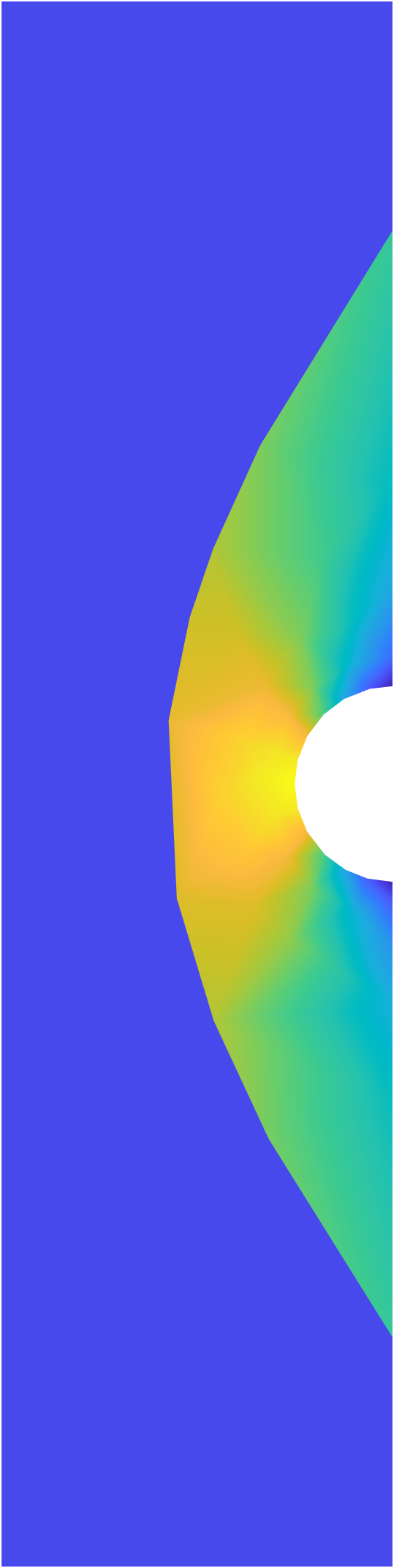};
				\end{groupplot}
			\end{tikzpicture}
			\colorbarMatlabParula{1.0}{1}{2}{3}{4.4}
				}
			\caption{HOIST iterations ($p=q=4$) projected to $p=0$, $q=1$ for $k=1$ (\textit{left}) and $k=100$ (\textit{middle}), $p=q=1$ for $k=100$ (\textit{right}) for the \ecyl~test case (density) shown with and without mesh edges.}
			\label{fig:euler:cyl:p1}
		\end{figure}

	Additionally, we investigate the impact of mesh refinement ($h$) and the approximation orders ($p$ and $q$) on the performance of the preconditioners. We apply the HOIST method to the bow shock problem, incorporating three additional refinement levels. Figure \ref{fig:euler:cyl:refinement} illustrates the density of converged solutions ($k=100$) and corresponding meshes resulting from this refinement process. To study the impact of polynomial degree, we first compute the HOIST solution on a grid with solution degree $p=4$ and mesh degree $q=4$. To obtain comparable lower order solutions, we restrict the $(p,q)=(4,4)$ solution to degrees $(p,q)=(3,3)$, $(p,q)=(2,2)$, $(p,q)=(1,1)$, and $(p,q)=(0,1)$ (Figure~\ref{fig:euler:cyl:p1}). We opt for this approach over computing a HOIST solution at a given polynomial degree to avoid situations where the HOIST iterations do not converge to a tracked configuration due to insufficient resolution.
	
	For all scenarios considered, the HOIST solver parameters \cite{huangRobustHighorderImplicit2022a} are set as follows
	\begin{itemize}
		\item Adaptive regularization $(\gamma_0,\gamma_\mathrm{min},\tau,\sigma_1,\sigma_2)=(10^{-2},10^{-2},10^{-1},10^{-2},10^{-1})$
		\item Adaptation of $\kappa$ $(\kappa_0,\kappa_\mathrm{min},\upsilon,\xi) = (1,0,1, 0.8)$
		\item Mesh operations $(c_1,c_2,c_3,c_4,c_4')=(0.025,10^{-10},5,0.025,10^{-2})$
		\item Reinitialization procedure $(c_5,c_6,c_7,c_8)=(0.5,10^{-2},0.5,10^{-2})$.
	\end{itemize}
	For brevity, these parameters have not been explicitly introduced in this manuscript. A complete description of all parameters and the overall algorithm can be found in \cite{huangRobustHighorderImplicit2022a}.

	\subsubsection{Supersonic flow over a two-dimensional diamond in a tunnel}
	\label{sec:numexp:euler:diamond}
\begin{figure}[!htbp]
	 		\centering
	 			 \ifbool{fastcompile}{}{
	 		\begin{tikzpicture}
	 			\begin{groupplot}[
	 				group style={
	 					group size=2 by 5,
	 					horizontal sep=0.4cm,
	 					vertical sep=0.4cm
	 				},
	 				width=0.57\textwidth,
	 				axis equal image,
	 				xticklabels={,,},
	 				yticklabels={,,},
	 				xmin=0, xmax=5,
	 				ymin=0, ymax=1.5
	 				]
	 				\nextgroupplot
	 				\addplot graphics [xmin=0, xmax=5, ymin=0, ymax=1.5] {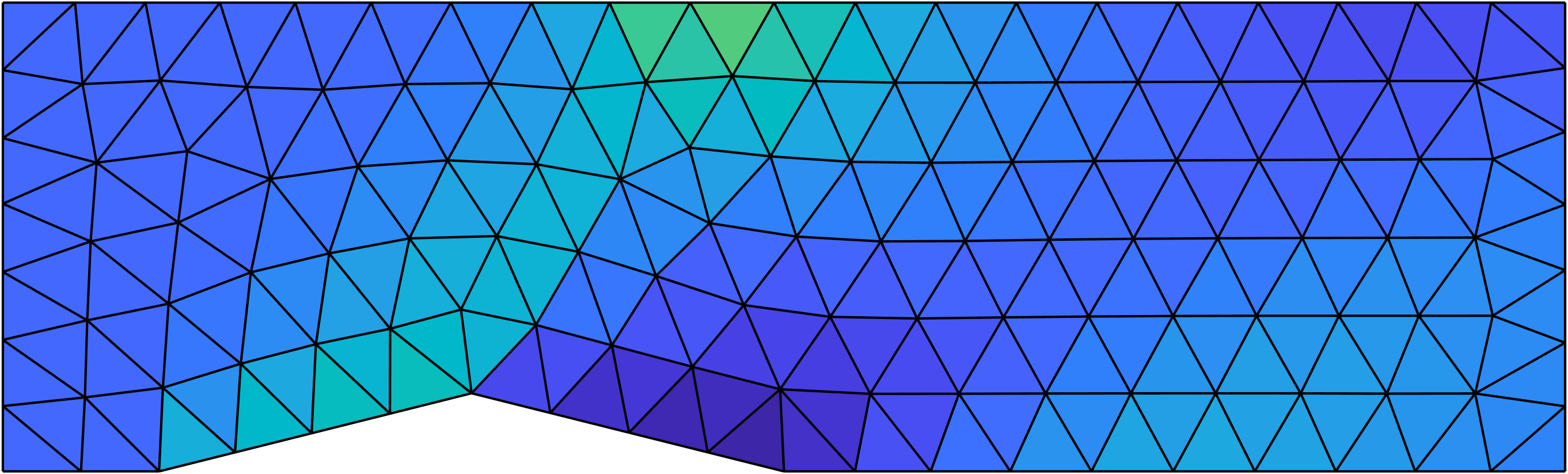};
	 				\nextgroupplot
	 				\addplot graphics [xmin=0, xmax=5, ymin=0, ymax=1.5] {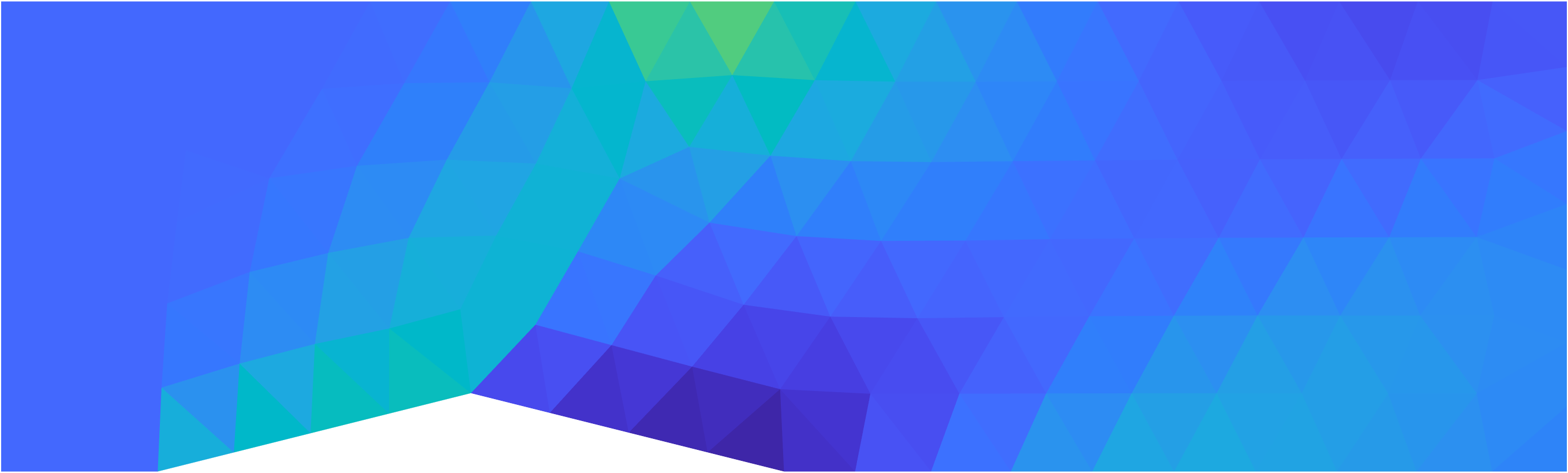};
	 				\nextgroupplot
	 				\addplot graphics [xmin=0, xmax=5, ymin=0, ymax=1.5] {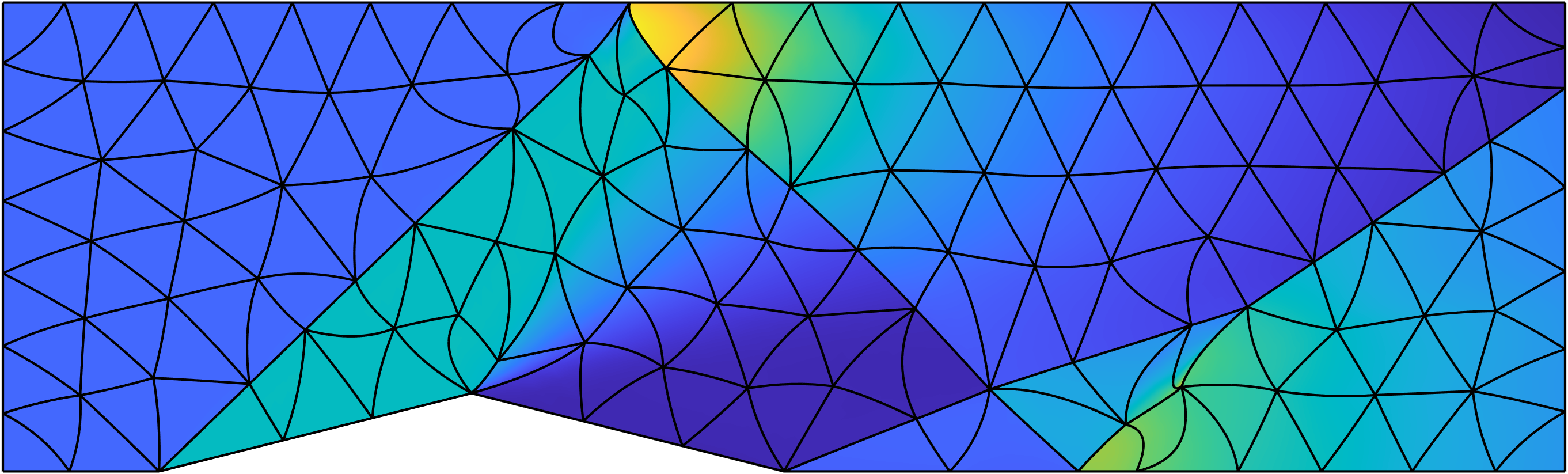};
	 				\nextgroupplot
	 				\addplot graphics [xmin=0, xmax=5, ymin=0, ymax=1.5] {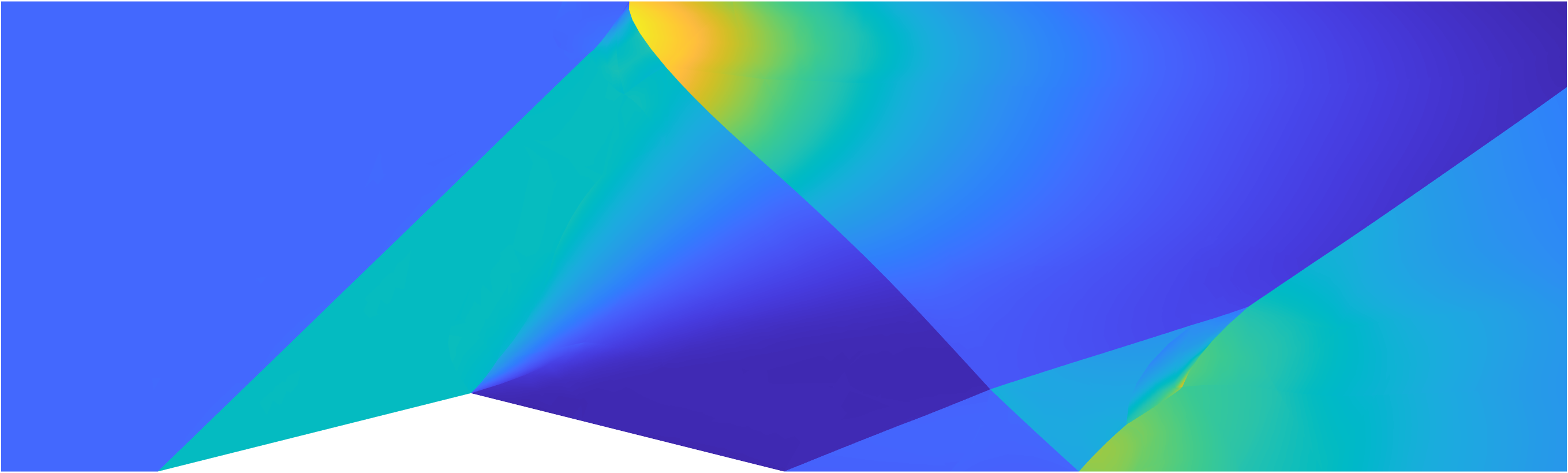};
	 				\nextgroupplot
	 				\addplot graphics [xmin=0, xmax=5, ymin=0, ymax=1.5] {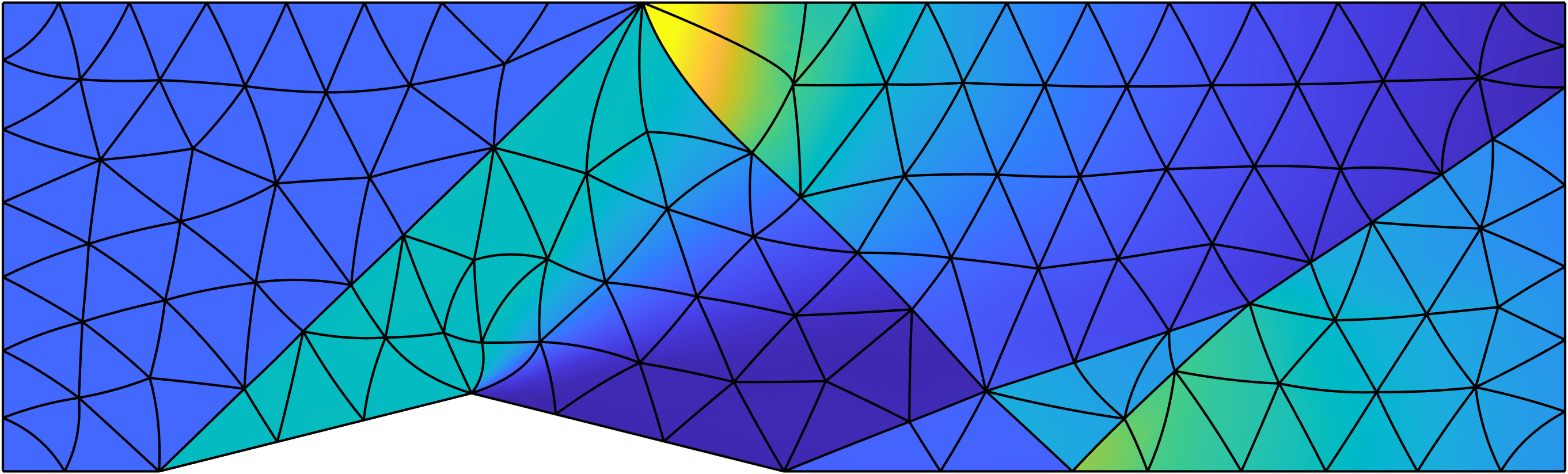};
	 				\nextgroupplot
	 				\addplot graphics [xmin=0, xmax=5, ymin=0, ymax=1.5] {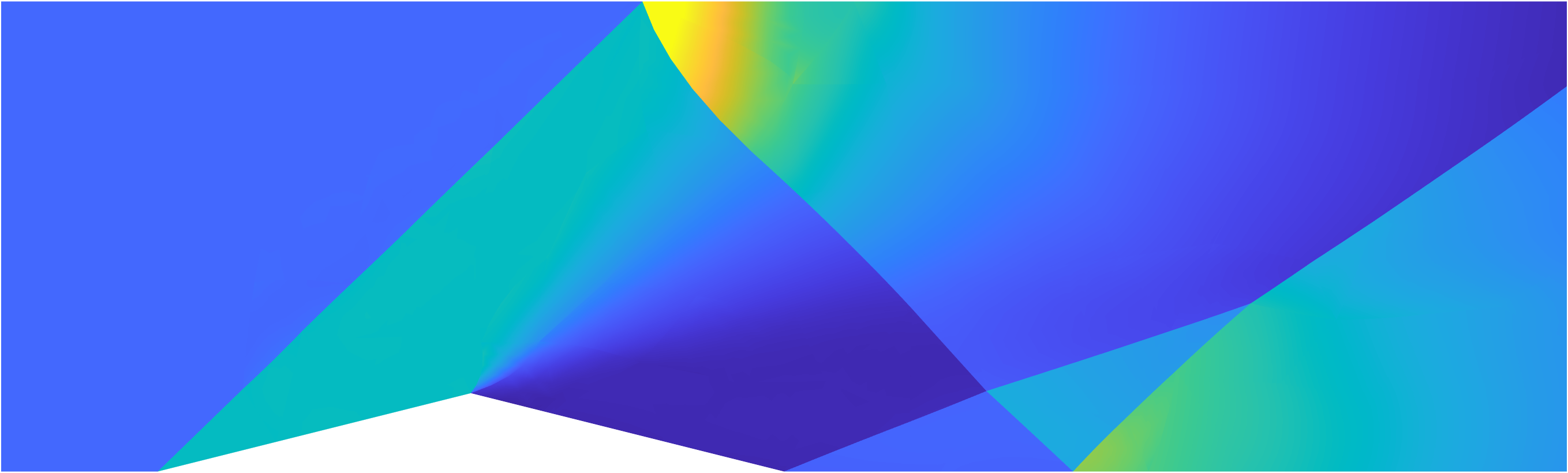};
	 			\end{groupplot}
	 		\end{tikzpicture}
	 		\colorbarMatlabParula{0.7}{1}{2}{3}{4.1}
	 		 }	
	 		\caption{Selected $p=q=2$ HOIST iterations $k\in \{1,150,300\}$ (\textit{top-to-bottom}) for the \edmnd~test case (density) with and without mesh edges.}
	 	\label{fig:euler:diamond:sltn}
	 \end{figure}

	Next, we study supersonic flow (Mach 2) passing over a two-dimensional diamond-shaped object within a tunnel (Figure~\ref{fig:euler:cyl:geom}). This test case, denoted as \edmnd, presents complex features such as reflecting, intersecting, and curved shocks. To discretize the domain, we generate an unstructured triangular mesh of 220 elements using DistMesh \cite{persson_simple_2004}. Similar to previous cases, a third-order approximation is employed for both the geometry and flow variables ($p=q=2$). Initializing the HOIST method with this unstructured mesh and the corresponding first-order finite volume solution, the method converges to a shock-aligned mesh that accurately represents all shocks and their intersections. Density plots of iterations at specific points in the optimization process ($k=1,150,300$) are highlighted in Figure \ref{fig:euler:diamond:sltn}. The corresponding states $\zbm_1,\zbm_{150}, \zbm_{300}$ will be utilized for the subsequent analysis.

	For all scenarios considered, the HOIST solver parameters \cite{huangRobustHighorderImplicit2022a} are set as follows
\begin{itemize}
	\item Adaptive regularization parameters $(\gamma_0,\gamma_\mathrm{min},\tau,\sigma_1,\sigma_2)=(10^{-2},10^{-2},10^{-1},10^{-2},10^{-1})$
	\item Mesh quality adaptation parameters $(\kappa_0,\kappa_\mathrm{min},\upsilon,\xi) = (1,0,1, 0.8)$
	\item Mesh operation parameters  $(c_1,c_2,c_3,c_4,c_4')=(0.25,10^{-10},4,0.3,10^{-2})$
	\item Reinitialization parameters $(c_5,c_6,c_7,c_8)=(0.5,10^{-2},0.5,10^{-2})$.
\end{itemize}  

\subsection{Results}
\label{sec:numexp:results}
In this section, we present and analyze the outcomes of the numerical experiments conducted in this study. These experiments gauge the efficiency of the preconditioners (Table \ref{tab:all_preconditioners}) designed for the HOIST method introduced in Section \ref{sec:linsolv:hoist:consideredprec}  and \ref{sec:mulgrd:hoist}, by comparing the required GMRES iterations to achieve the convergence criterion (\ref{eqn:gmres-convcrit-ex}). Due to the multitude of parameters potentially influencing the condition of the linear system, we conduct separate investigations for each relevant parameter. First, we explore the impact of the mesh adaptation parameter $\kappa$ on the GMRES iterations (Section \ref{sec:numexp:results:kappa}). Subsequently, we delve into the influence of variations in the state $\zbm_k$ (Section \ref{sec:numexp:results:state}) and regularization parameter $\gamma$ (Section \ref{sec:numexp:results:gamma}). Furthermore, we study the dependence on the polynomial degrees $p$, $q$ (Section \ref{sec:numexp:results:p}) and mesh size $\vert \Ecal_h \vert$. Finally, we illustrate the GMRES iterations required throughout the entire optimization history for both examples (Section~\ref{sec:numexp:results:evrystt}).

\subsubsection{Influence of mesh quality parameter $\kappa$}
\label{sec:numexp:results:kappa}
	\begin{figure}[!htbp]
		\centering
		\ifbool{fastcompile}{}{
		\input{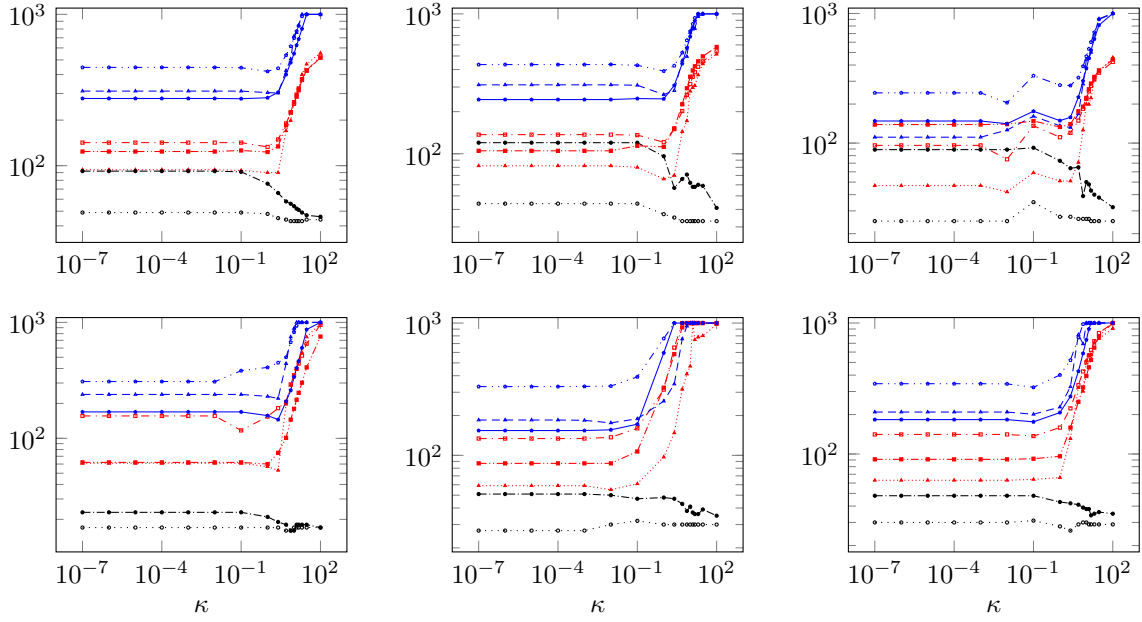}
	}
		\caption{(\# GMRES iterations) vs. mesh quality parameter $\kappa$ for both test cases (\textit{top}: \ecyl, \textit{bottom}: \edmnd) (legend in Table~\ref{tab:legend})	for polynomial degree $p=q=2$, regularization parameters $\gamma=0.1$, and states $z_k$ (\ecyl: $ k=1,50,100$, \edmnd: $ k=1,150,300$) (\textit{left-to-right}).}
		\label{fig:kppStdy}
	\end{figure}
In our first experiment, we aim to measure the impact of the choice of $\kappa$ by considering both test cases (\ecyl\ and \edmnd\ with $p=q=2$) with the linear system formed from specific states $\zbm_k$ corresponding to different HOIST iterations ($k=1,50,100$ for \ecyl~and $k=1,150,300$ for \edmnd ) with a single regularization parameter of $\gamma=0.1$. For these six configuration, we compute the HOIST matrix for $\kappa = \{10^{-10},10^{-9}, \hdots, 10^{2}\}$ and measure the needed GMRES iterations for each preconditioner tested (Table~\ref{tab:all_preconditioners}) (results in Figure~\ref{fig:kppStdy}).

Analyzing the results shown in Figure \ref{fig:kppStdy}, we make the following observations. Most preconditioners demonstrate significant deterioration after a critical $\kappa$, usually in the range $\kappa\in[0.1,1]$. As $\kappa$ increases beyond this range, the GMRES iterations of all preconditioners except \PrecV{0} and \PrecV{0p0} quickly increase, which suggests the approximations to $\Bbm_{\ybm \ybm}$ (the only block depending on $\kappa$) deteriorates as $\kappa$ rises beyond the critical value. Fortunately, such large $\kappa$ values rarely occur in the method, mitigating this sensitivity issue. Intriguingly, when $\kappa <0.1$, no sensitivity was observed for any preconditioner so we fix $\kappa=10^{-7}$ for the remaining studies in this work. Furthermore, we observe that as $\kappa$ rises beyond $0.1$, the benefit of $p$-multigrid diminishes.

For all values of $\kappa$, the expensive, best-case scenario $\PrecV{0}$ preconditioner outperforms all others across all states considered for both problems. It is also interesting to note that $p$-multigrid actually degrades the performance of the $\PrecV{0}$, in some cases making it worse than preconditioners that use approximate inverses. For the \ABILU\ preconditioners, the addition of $p$-multigrid and especially the inclusion of $\text{ilu}(\Bbm_{\ybm \ybm})$ as an approximation to $\Bbm_{\ybm \ybm}$ significantly enhances its performance (in some cases, reducing the GMRES iterations by a factor of two or more). The \ABJ\ preconditioners also benefit from both $p$-multigrid and the inclusion of $\text{ilu}(\Bbm_{\ybm \ybm})$ as an approximation to $\Bbm_{\ybm \ybm}$; however, in this case, $p$-multigrid provides the greater reduction in GMRES iterations. Finally, as expected, the \ABILU\ outperform the \ABJ\ preconditioners across test cases and states.


\subsubsection{Influence of linearization state $z_k$}
\label{sec:numexp:results:state}
\begin{figure}[!htbp]
	\centering
	\ifbool{fastcompile}{}{
		\input{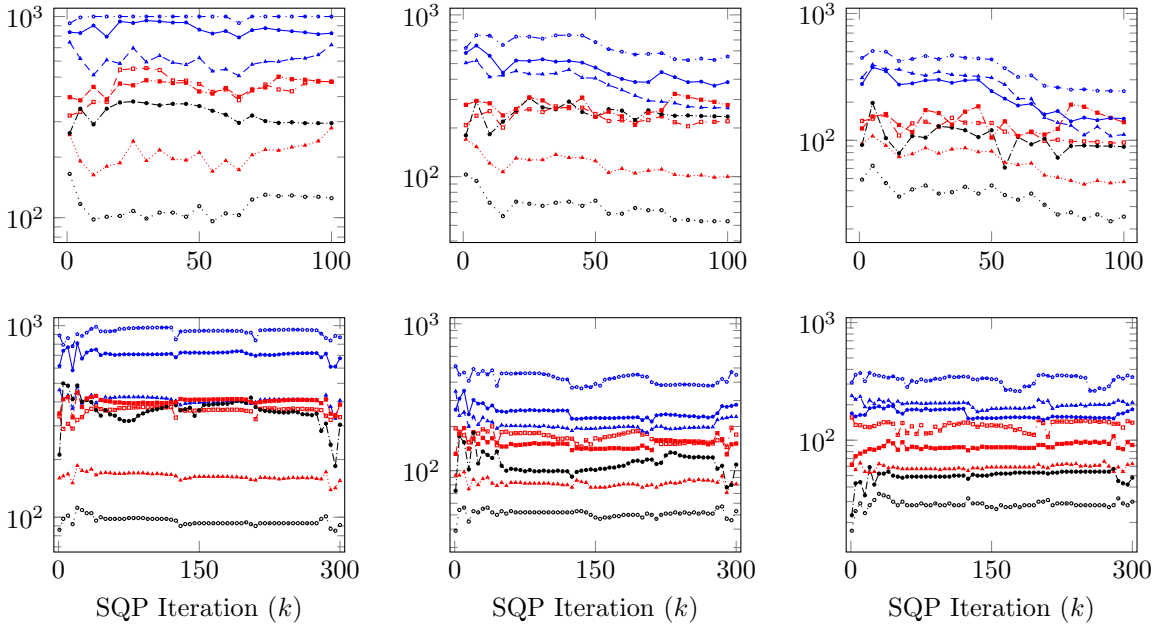}
	}
	\caption{(\# GMRES iterations) vs. state $\zbm_k$ for both test cases (\textit{top}: \ecyl, \textit{bottom}: \edmnd)  (legend in Table~\ref{tab:legend}) for polynomial degree $p=q=2$, mesh quality parameter $\kappa = 10^{-7}$, and different regularization parameters $\gamma=10^{-3},10^{-2},10^{-1}$ (\textit{left-to-right}).}
	\label{fig:stStdy}
\end{figure}

In our second experiment, we investigate the dependence of the preconditioner performance on the linearization state $\zbm_k$. We fix the mesh quality parameter $\kappa=10^{-7}$ and build six test cases from the two problems (\ecyl\ and \edmnd\ with $p=q=2$) and three choices for the regularization parameter $\gamma \in \{10^{-3},10^{-2},10^{-1}\}$. For each test case and preconditioner (Table~\ref{tab:all_preconditioners}), we record the number of GMRES iteration required to reach the convergence criteria \eqref{eqn:gmres-convcrit-ex} at every $5$th HOIST iteration, i.e., $k \in \{1,5,10,15,...,150\} $ for \ecyl~and $k \in \{1,5,10,15,...,300\} $ for \edmnd (Figure~\ref{fig:stStdy}).

Analyzing the results shown in Figure \ref{fig:stStdy}, we make the following observations. First, the linearization state has a modest impact on the GMRES iterations. In the \edmnd~ case, the iterations remain nearly constant, with minor fluctuations occurring due to abrupt state changes (e.g., solution reinitialization and element collapses). The \ecyl~ case exhibits a more pronounced state dependency, particularly for larger $\gamma$, where the iteration count tends to decrease as the final state is approached. Once again, the \PrecV{0} demonstrates superior performance across all scenarios and \ABILUilu\ is the best practical preconditioner (i.e., not involving the expensive $\Ju$, $\Ju^T$, and $\Bbm_{\ybm\ybm}$ inverses). Again, the $p$-multigrid counterpart of $\PrecV{0}$, $\PrecV{0p0}$, performs noticeably worse, in many cases requiring more iterations that \ABILUilu\ and often demonstrating similar performance to \ABILU~and \ABILUp, despite the use of exact inverses. Unlike the previous study, there is no clear conclusion regarding \ABILU~and \ABILUp. Finally, the addition of $p$-multigrid and the inclusion of $\text{ilu}(\Bbm_{\ybm \ybm})$ as an approximation to $\Bbm_{\ybm\ybm}$ enhance the performance of the \ABJ\ preconditioner with \ABJilu\ holding a clear advantage for smaller values of $\gamma$.


\subsubsection{Influence of regularization parameter $\gamma$}
\label{sec:numexp:results:gamma}
\begin{figure}[!htbp]
	\centering
	\ifbool{fastcompile}{}{
		\input{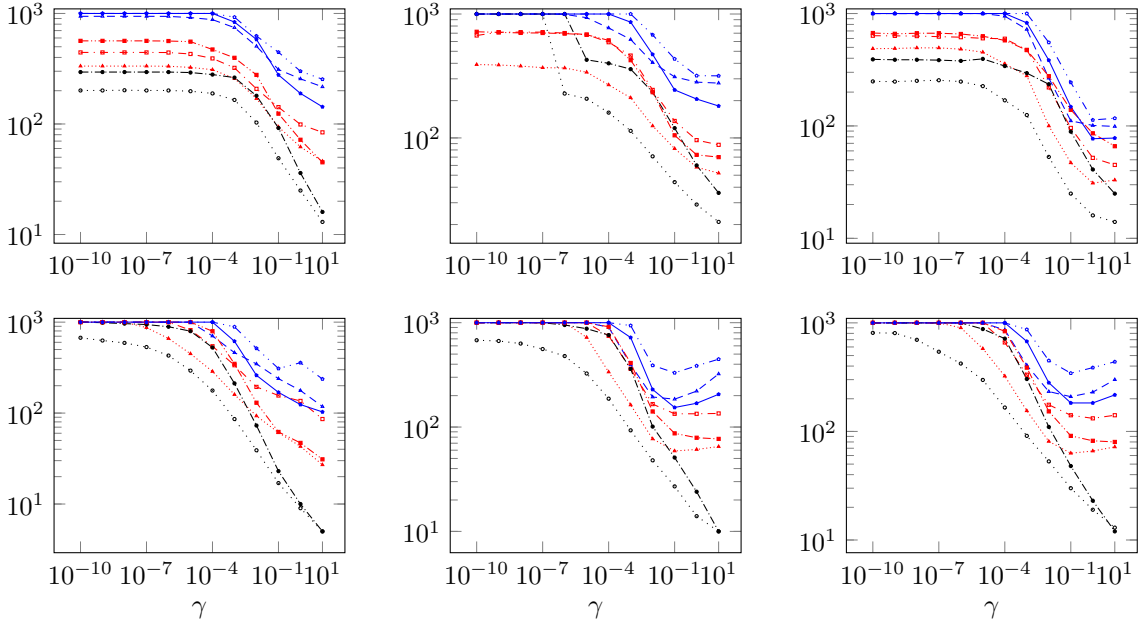}
	}
	\caption{(\# GMRES iterations) vs. regularization parameter $\gamma$ for both test cases (\textit{top}: \ecyl, \textit{bottom}: \edmnd) (legend in Table~\ref{tab:legend}) for polynomial degree $p=q=2$, mesh quality parameter $\kappa = 10^{-7}$, and different states $\zbm_k$  (\ecyl: $k=1,50,100$, \edmnd: $ k=1,150,300$) (\textit{left-to-right}).}
	\label{fig:lamStdy}
\end{figure}

In this experiment, we study the influence of the Hessian regularization parameter $\gamma$ on the performance of the preconditioners considered (Table~\ref{tab:all_preconditioners}). We build six test cases from the two problems (\ecyl\ and \edmnd\ with $p=q=2$) and three states ($k \in \{1,50,100\} $ for \ecyl\ and $k \in \{1,150,300\} $ for \edmnd). Furthermore, we fix the mesh quality parameter at $\kappa=10^{-7}$ and vary the regularization parameter $\gamma\in \{10^{-10},10^{-9},\hdots,10^{1}\}$. The resulting GMRES iterations needed to reach the convergence criteria \eqref{eqn:gmres-convcrit-ex} for each of these cases are shown in Figure~\ref{fig:lamStdy}.

Analyzing the results shown in Figure \ref{fig:lamStdy}, we make the following observations. Decreasing the parameter $\gamma$ reduces the regularization applied to the matrix $\Bbm_{\ybm \ybm}$, leaving the ill-conditioned (or singular) Gauss-Newton Hessian in the limit where $\gamma = 0$. As expected, this leads to a noticeable rise in the number of GMRES iterations, particularly evident in the case of \PrecV{0} and \PrecV{0p0}. However, for the other preconditioners, we observe a relative indifference to variations for $\gamma \in [10^{-1},10^{1}]$, particularly in later SQP iterations ($k>100$). This suggests that in this range, the loss in accuracy incurred by the approximations of $\Ju$ and $\Bbm_{\ybm\ybm}$ dominates ill-conditioning effects.

The results suggest the presence of a problem-dependent threshold value for $\gamma$ (\ecyl: $10^{-5}$, \edmnd: $10^{-9}$). Below this threshold, the number of iterations ceases to increase significantly. This phenomenon is especially prominent in the case of \ecyl, whereas for \edmnd, most preconditioners did not converge reaching the maximum number of iterations below $\gamma=10^{-6}$. Additionally, these findings imply the possibility of establishing a lower limit for the minimum regularization parameter $\gamma_{\text{min}}$ that should be set in the HOIST method. The results obtained for \edmnd~suggest that $\gamma_{\text{min}}$ should not be less than $10^{-4}$, as the iteration counts become impractical beyond this threshold. Considering the observed increase in iteration numbers with higher polynomial degrees (Section \ref{sec:numexp:results:p}) and finer meshes (Section \ref{sec:numexp:results:h}), setting a more conservative lower bound, for instance, $\gamma_{\text{min}}=10^{-2}$, is advisable.

For the BJ-based preconditioners (\ABJ, \ABJp, \ABJilu), the trends observed earlier remain evident: for large regularization parameters $\gamma \in [10^{-1},10^1]$, \ABJp~outperforms \ABJilu~while the opposite is true for $\gamma \leq 10^{-2}$. Both of these preconditioners perform favorly compared to \ABJ. The scenario is slightly different for the BILU-based preconditioners (\ABILU, \ABILUp, \ABILUilu). In this case, \ABILUilu~outperforms both \ABILU~and \ABILUp~across all cases, with the performance gap between \ABILUp~and \ABILUilu~widening for $\gamma \leq 10^{-2}$. The utilization of $p$-multigrid seems to add value only for $\gamma>10^{-2}$, as \ABILU~often exhibits similar or even better iteration counts than \ABILUp. 


\subsubsection{Influence of polynomial degrees $(p,q)$}
\label{sec:numexp:results:p}
\begin{figure}[!htbp]
	\centering
	\ifbool{fastcompile}{}{
		\input{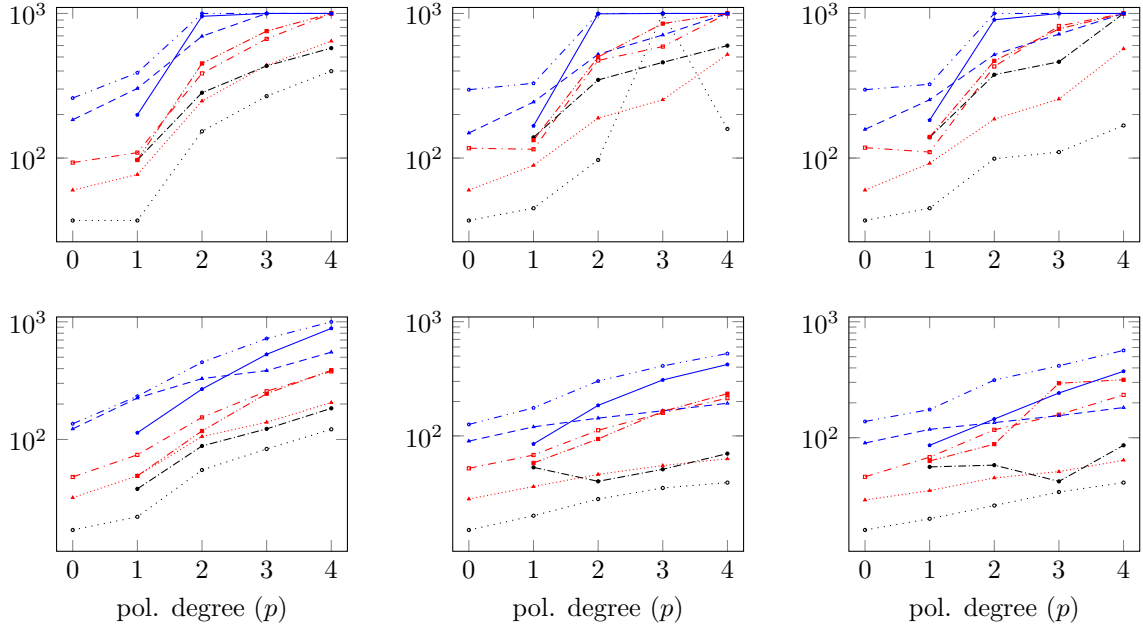}
	}
	\caption{(\# GMRES iterations) vs. polynomial degrees $(p,q)$ for different regularization parameters $\gamma=10^{-3}$ (\textit{top}) and $\gamma=10^{-1}$ (\textit{bottom}) (legend in Table~\ref{tab:legend}), mesh quality parameter $\kappa = 10^{-7}$, and different states $\zbm_k$, ($k=1,50,100$) for the \ecyl\ problem. For $p=0$, the coarse-scale updates from the $p$-multigrid preconditioners solve the problem directly, which only requires one GMRES iteration. These results are omitted for clarity.}
	\label{fig:pStdy}
\end{figure}

In this experiment, we study the effect of the polynomial degree $(p,q)$ on the GMRES iterations. We test each of our proposed preconditioners (Table~\ref{tab:all_preconditioners}) against six cases built from three states $\zbm_k$ for $k\in\{1,50,100\}$, two regularization parameters $\gamma \in \{10^{-3},10^{-1}\}$, and a fixed mesh quality parameter $\kappa=10^{-7}$ for the \ecyl\ problem. A $p=q=4$ HOIST simulation is used to compute the initial states ($\zbm_k$ for $k=1,50,100$), which are subsequently restricted to polynomial degrees $(p,q) \in \{(0,1),(1,1),(2,2),(3,3),(4,4)\}$. As discussed in Section~\ref{sec:numexp:euler:cylinder}, this approach is taken to yield a well-defined, systematic study and avoid HOIST convergence issues that can arise when the grid is sufficiently underresolved. The measured GMRES iterations required to achieve the convergence criteria \eqref{eqn:gmres-convcrit-ex} are depicted in Figure \ref{fig:pStdy}.

Analyzing the results shown in Figure \ref{fig:pStdy}, we make the following observations. First, increasing the polynomial degree on a fixed mesh results in a direct escalation of GMRES iterations for all preconditioners with more pronounced growth rate for the smaller regularization parameters $\gamma=10^{-3}$. The $p$-multigrid versions of the BJ (\ABJp) and BILU (\ABILUp) preconditioners are sensitive to the polynomial degree as their iteration count approaches that of the original BJ (\ABJ) and BILU (\ABILU) preconditioner as the polynomial degree increases. Both the original and $p$-multigrid version of the BJ and BILU are outperformed by inclusion of $\text{ilu}(\Bbm_{\ybm \ybm})$ as an approximation to $\Bbm_{\ybm\ybm}$, where $\ABJilu$ is the most effective BJ preconditioner and $\ABILUilu$ is the most effective BILU preconditioner. Furthermore, the $\ABILUilu$ preconditioner is the most effective practical preconditioner, only being outperformed by the best-case scenario $\PrecV{0}$ (and, in some cases, its $p$-multigrid variant). The $\ABILUilu$ preconditioner also exhibits the slowest iteration growth with polynomial degree, particularly for the larger regularization parameter $\gamma=10^{-1}$.

\subsubsection{Influence of number of mesh elements $\vert \Ecal_h \vert$}
\label{sec:numexp:results:h}
\begin{figure}[!htbp]
	\centering
	\ifbool{fastcompile}{}{
		\input{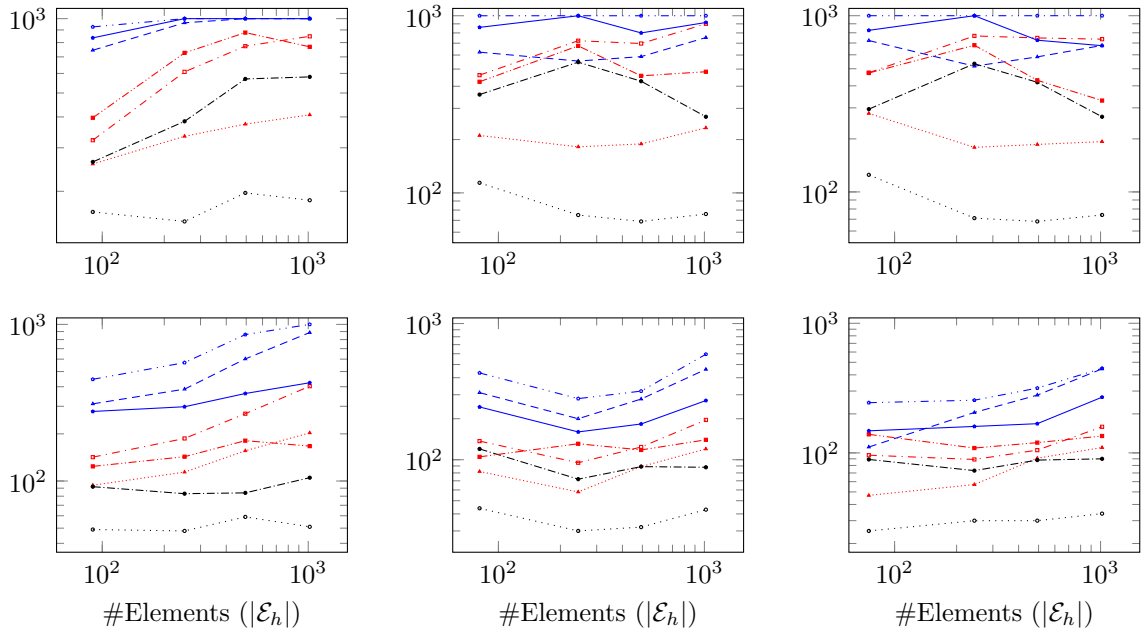}
	}
	\caption{(\# GMRES iterations) vs. (\# Elements $\vert \Ecal_h \vert $) for  different regularization parameters $\gamma=10^{-3}$ (\textit{top}) and $\gamma = 10^{-1}$ (\textit{bottom}) (legend in Table~\ref{tab:legend}), polynomial degree $p=q=2$, mesh quality parameter $\kappa = 10^{-7}$, and different states $\zbm_k$ (\ecyl: $ k=1,50,100$, \edmnd: $ k=1,150,300$) (\textit{left-to-right}).}
	\label{fig:grdStdy}
\end{figure}

In this experiment, we study the dependency of the GMRES iterations on the number of mesh elements. We test each of our proposed preconditioners (Table~\ref{tab:all_preconditioners}) against six cases built from three states $\zbm_k$ for $k\in\{1,50,100\}$, two regularization parameters $\gamma \in \{10^{-3},10^{-1}\}$, and a fixed mesh quality parameter $\kappa=10^{-7}$ for the \ecyl\ problem. For each of these cases, we consider four refinement levels (Figure~\ref{fig:euler:cyl:refinement}) at fixed polynomial degree $p=q=2$ with element count $\vert \Ecal_h \vert \in \{70,130,260,1000\}$. The measured GMRES iterations required to achieve the convergence criteria \eqref{eqn:gmres-convcrit-ex} are depicted in Figure \ref{fig:grdStdy}.

Analyzing the results shown in Figure \ref{fig:grdStdy}, we make the following observations. The exact preconditioner $\Ao$ demonstrates remarkable insensitivity to the number of elements. Its multigrid counterpart, \Aop, while less effective, exhibits a similar stable trend for the $\gamma=10^{-1}$ case. In some cases, the GMRES iteration count slightly decreases as the number of elements rises. For the BJ preconditioners, \ABJilu\ is most effective for the smaller $\gamma = 10^{-3}$ (the other BJ variants often reach the maximum iterations without convergence), whereas \ABJp\ is the most effective BJ preconditioner for $\gamma = 10^{-1}$ (although usually only slightly outperforms \ABJilu). Similarly, for the smaller $\gamma = 10^{-3}$, the \ABILUilu\ preconditioner is clearly superior to the other BILU variants and exhibits the slowest growth as the element count rises. The \ABILUilu\ is usually the best BILU preconditioner for the larger $\gamma=10^{-1}$, although the difference between the three BILU preconditioners is less dramatic for this scenario.


\subsubsection{Comparison of preconditioners for adaptive $\kappa_k,\gamma_k$}
\label{sec:numexp:results:evrystt}
\begin{figure}[!htbp]
	\centering
	\ifbool{fastcompile}{}{
		\input{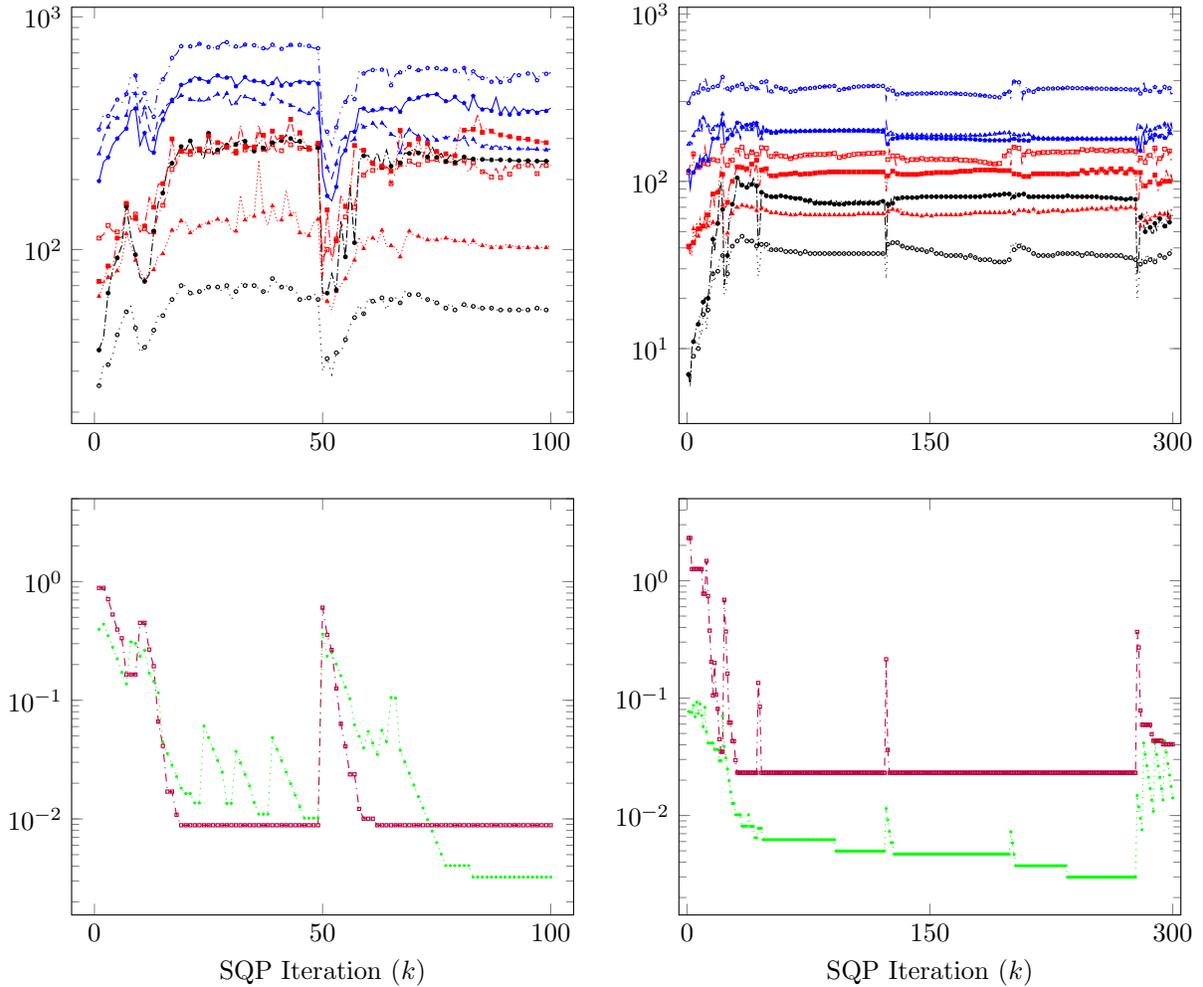}
	}
	\caption{\textit{Top}: (\# GMRES iterations) vs. states $\zbm_k$ using adaptive mesh/regularization parameters $\kappa_k$, $\gamma_k$ (legend in Table~\ref{tab:legend}) for both test cases with $p=q=2$ (\textit{left}: $\ecyl$, \textit{right}: $\edmnd$). \textit{Bottom}: History of the adaptive mesh quality parameter $\kappa_k$ (\ref{line:Kappa}) and regularization parameter $\gamma_k$ (\ref{line:Gamma}).}
	\label{fig:EStSStdy}
\end{figure}

In our final experiment, we investigate GMRES iterations across the entire optimization history for both problems (\ecyl\ and \edmnd\ with $p=q=2$). For this experiment, we use the adaptive mesh quality $\kappa_k$ and regularization parameters $\gamma_k$ from \cite{huangRobustHighorderImplicit2022a} with adaptation parameters in Sections~\ref{sec:numexp:euler:cylinder}-\ref{sec:numexp:euler:diamond}. The measured GMRES iterations required to achieve the convergence criteria \eqref{eqn:gmres-convcrit-ex} for each state $\zbm_k$ ($k \in \{1,2,\dots,100\}$ for \ecyl\ and $k \in \{1,2,\dots,300\}$ for \edmnd) encountered during the HOIST iterations are depicted in Figure \ref{fig:EStSStdy}. The evolution of $\gamma_k$ and $\kappa_k$ are also shown in this figure.

Analyzing the results shown in Figure \ref{fig:EStSStdy}, we observe the iteration count closely correlates to the $\gamma$ value, as expected from Section~\ref{sec:numexp:results:gamma}, in that the GMRES iterations rise as $\gamma$ decreases. However, extreme values of $\gamma$ are not encountered during the adaptation, which avoids excessive GMRES iteration counts. Abrupt changes in GMRES iterations are associated with abrupt alterations in $\gamma$ (e.g., in the \ecyl~case around $k=50$, where $\gamma$ is nearly equal to its initial value), which occur after elements are collapsed. Generally, larger $\gamma$ values tend to benefit the $p$-multigrid preconditioners the most, granting them an advantage over their counterparts (though \Ao~is an exception due to the overall poor performance of \Aop). For the BJ preconditioners, \ABJp~consistently outperforms \ABJ~and is on par with \ABJilu~for the \edmnd~cases. For the \ecyl~problem, \ABILUilu~performs better in the low $\gamma$ regime. Among the BILU preconditioners, \ABILUilu~consistently performs the best across all $k$, only matching \ABILUp~for high $\gamma$ values. As expected from the previous sections, \ABILUilu\ is the most effective practical preconditioner as it is only consistently outperformed by the best-case (but impractical) \Ao, making it our preferred preconditioner.


\section{Conclusion}
\label{sec:conclude}

In this work, we introduced matrix-based preconditioners for constrained high-order implicit shock tracking methods and thoroughly tested their performance across various critical parameters in the optimization solver. While we focused on the HOIST method that uses the enriched DG residual as the objective function, the preconditioners would apply to other constrained implicit shock tracking formulations.
By analyzing the block structure of the implicit shock tracking linear system and the sparsity of each block, we devised a family of approximate block anti-triangular preconditioners that integrate common DG preconditioners such as block Jacobi and block ILU0 with minimum discarded fill reordering. We also introduced a two-level $p$-multigrid scheme that can be combined with any of the proposed preconditioners.

All preconditioners were rigorously evaluated on two compressible inviscid flow problems, focusing on the number of GMRES iterations required to achieve a prescribed relative error norm. Our investigations revealed that the regularization parameter $\gamma$ has the most significant impact on GMRES iterations, with the proposed $p$-multigrid scheme offering added value only under high regularization conditions. The iteration count is sensitive to the polynomial degree of the solution and mesh, particularly in low regularization settings. Conversely, the number of mesh elements and the mesh quality parameter displayed relatively minor influence on the required GMRES iterations, with the latter being insignificant for parameters commonly used in practice. Overall, our findings highlighted the promising performance of BILU-based preconditioners across various problem and parameter configurations. The BILU variant that uses an ILU0 approximation to $\Bbm_{\ybm\ybm}$ emerged as the best and most reliable of all the practical preconditioners (i.e., those that did not require $\Ju$, $\Ju^T$, and $\Bbm_{\ybm\ybm}$ inverses). Our investigations also concluded that the two-level $p$-multigrid scheme did not yield sufficient advantages to justify its computational cost.

Interesting and relevant avenues of future research include the generalizability of the proposed preconditioners (and observations regarding their performance) to larger problems, the massively parallel setting, and viscous problems. In the parallel setting, it will be pertinent to study both CPU and wall time to reach convergence, in addition to the GMRES iteration count. It is well-known that the cost and memory requirements of direct solvers grow rapidly with problem size, and they scale poorly in parallel. On the other hand, this work showed the number of GMRES iterations was relatively insensitive to the problem size, which makes the proposed solvers a promising alternative to sparse direct solvers. A quantitative comparison to sparse direct solvers would be interesting to establish: (1) the iterative solver tolerances required to obtain robust shock tracking results (e.g., comparable to those obtained with a direct solver) and (2) the problem sizes at which the iterative solver becomes competitive with a direct solver. In the viscous setting, we expect the $p$-multigrid approach to provide additional benefit as the viscosity increases as observed in the DG context \cite{persson_newton-gmres_2008}.

\section*{Acknowledgments}
This work is supported by AFOSR award numbers FA9550-20-1-0236,
FA9550-22-1-0002, FA9550-22-1-0004, ONR award number
N00014-22-1-2299, and NSF award number CBET-2338843.
The first author is supported by the Graduate School CE within
Computational Engineering at Technische Universität Darmstadt.
The content of this publication does not necessarily reflect the position
or policy of any of these supporters, and no official endorsement should
be inferred.

\bibliographystyle{plain}
\bibliography{biblio}

\end{document}